\documentclass[11pt]{amsart}
\usepackage[hmargin=3cm,vmargin=3cm]{geometry}

\usepackage{combelow} 

\usepackage[latin1,utf8]{inputenc} 
\usepackage[english]{babel}
\usepackage[backend=bibtex,doi=false,url=false,isbn=false,style=alphabetic,maxnames=6]{biblatex}
\usepackage{etex}
\AtBeginBibliography{\small}

\makeatletter
\renewcommand\part{%
   \if@noskipsec \leavevmode \fi
   \par
   \addvspace{4ex}%
   \@afterindentfalse
   \secdef\@part\@spart}
   
\def\@part[#1]#2{%
    \ifnum \c@secnumdepth >\m@ne
      \refstepcounter{part}%
      \addcontentsline{toc}{part}{\thepart\hspace{1em}#1}%
    \else
      \addcontentsline{toc}{part}{#1}%
    \fi
    {\parindent \z@ \raggedright
     \interlinepenalty \@M
     \normalfont
     \ifnum \c@secnumdepth >\m@ne
       \Large\bfseries \partname\nobreakspace\thepart
       \par\nobreak
     \fi
     \huge \bfseries #2%
     \par}%
    \nobreak
    \vskip 3ex
    \@afterheading}
\def\@spart#1{%
    {\parindent \z@ \raggedright
     \interlinepenalty \@M
     \normalfont
     \huge \bfseries #1\par}%
     \nobreak
     \vskip 3ex
     \@afterheading}
\makeatother

\usepackage{xspace,amssymb,amsfonts,euscript}
\usepackage{amsthm,amsmath,amscd,stmaryrd,latexsym,youngtab,verbatim}
\usepackage[nohug]{diagrams}
\usepackage{palatino, euscript}

\usepackage{graphicx}
\usepackage[all]{xy}
\SelectTips{cm}{}

\usepackage{tikz, tikz-cd}
\usetikzlibrary{matrix, arrows, calc}
\tikzset{frontline/.style={preaction={draw=white,-,line width=6pt}},}  

\usepackage[dvips]{epsfig}
\usepackage{pinlabel}
\usepackage{psfrag}

\IfFileExists{comments.sty}{\usepackage{comments}}

\usepackage{enumitem}

\RequirePackage{color}
\definecolor{myred}{rgb}{0.75,0,0}
\definecolor{mygreen}{rgb}{0,0.5,0}
\definecolor{myblue}{rgb}{0,0.25,0.65}
\definecolor{references}{rgb}{0,0,1}

  \RequirePackage[pdftex,
   colorlinks = true,
   urlcolor = references, 
   citecolor = references, 
   linkcolor = references, 
]
{hyperref}

\newtheorem{thm}{Theorem}[section]
\newtheorem{lemma}[thm]{Lemma}
\newtheorem{theorem}[thm]{Theorem}

\newtheorem{prop}[thm]{Proposition}
\newtheorem{proposition}[thm]{Proposition}
\newtheorem{cor}[thm]{Corollary}
\newtheorem{corollary}[thm]{Corollary}

\newtheorem{conj}[thm]{Conjecture}
\newtheorem{conjecture}[thm]{Conjecture}

\newtheorem*{prop*}{Proposition}
\newtheorem*{lemma*}{Lemma}

\theoremstyle{definition}
\newtheorem{defn}[thm]{Definition}
\newtheorem{definition}[thm]{Definition}

\newtheorem{notation}[thm]{Notation}
\newtheorem{ex}[thm]{Example}
\newtheorem{example}[thm]{Example}

\theoremstyle{remark}
\newtheorem{remark}[thm]{Remark}

\numberwithin{equation}{section}

    \def\KM{{\mathbb{K}}}
    
  \def\mg{{\mathfrak m}}

    \def\PM{{\mathbb{P}}}
    \def\QM{{\mathbb{Q}}}
    
    \def\SM{{\mathbb{S}}}

    \def\ZM{{\mathbb{Z}}}

  \def\ab{{\mathbf a}}  \def\AC{{\mathcal{A}}}
    \def\BC{{\mathcal{B}}}
\def\CB{{\mathbf C}}  
\def\CC{{\mathcal{C}}}
\def\DB{{\mathbf D}}

\def\IB{{\mathbf I}}    
    \def\JC{{\mathcal{J}}}
\def\KB{{\mathbf K}}    \def\KC{{\mathcal{K}}} \def\K{{\KC}}
\def\LB{{\mathbf L}}    \def\LC{{\mathcal{L}}}
    \def\MC{{\mathcal{M}}}
    
    \def\OC{{\mathcal{O}}}
\def\PB{{\mathbf P}}    \def\PC{{\mathcal{P}}}
\def\QB{{\mathbf Q}}    
    
   \def\SC{{\mathcal{S}}} 
    \def\TC{{\mathcal{T}}}
    \def\UC{{\mathcal{U}}}
    \def\VC{{\mathcal{V}}}
    
    \def\XC{{\mathcal{X}}}
    \def\YC{{\mathcal{Y}}}
    \def\ZC{{\mathcal{Z}}}

\def\KS{{\EuScript K}}

\def\a{\alpha}
\def\b{\beta}
\def\g{\gamma}

\def\d{\delta}

\def\e{\varepsilon}
\def\k{\kappa}
\def\l{\lambda}
\def\L{\Lambda}

\let\phi=\varphi

\usepackage{bbm}
\def\C{{\mathbbm C}}

\def\Z{{\mathbbm Z}}
\def\Q{{\mathbbm Q}}
\def\1{\mathbbm{1}}
\newcommand{\one}{\1}

\newcommand{\firstMap}{\a}
\newcommand{\secondMap}{\b}
\newcommand{\firstProj}{\PB_\a}
\newcommand{\secondProj}{\PB_\b}

\newcommand{\firstShift}{\l}
\newcommand{\secondShift}{\mu}
\newcommand{\W}{\Omega}

\newcommand{\cidemp}{\CB}
\newcommand{\Coh}{\operatorname{Coh}}

\newcommand{\FHilb}{\operatorname{FHilb}}
\newcommand{\ft}{\operatorname{ft}}

\newcommand{\Cij}[2]{\CB_{#1,#2}}
\newcommand{\TL}{\operatorname{TL}}
\newcommand{\BN}{\TC\LC}

\newcommand{\idempotent}{\idemp}
\newcommand{\idemp}{\PB}
\newcommand{\otherIdemp}{\QB}

\newcommand{\Hecke}{\mathcal{H}}

\renewcommand{\setminus}{\smallsetminus}
\newcommand{\Homb}{\Hom^\bullet}
\newcommand{\Homg}{\operatorname{HOM}}
\newcommand{\Endg}{\operatorname{END}}
\newcommand{\smMatrix}[1]{\left[\begin{smallmatrix}#1\end{smallmatrix}\right]}

\newcommand{\im}{\operatorname{im}}
\newcommand{\sign}{\operatorname{sgn}}
\newcommand{\un}{\underline}

\newcommand{\ot}{\otimes}
\newcommand{\pa}{\partial}
\newcommand{\co}{\colon}

\renewcommand{\to}{\rightarrow}
\newcommand{\into}{\hookrightarrow}

\newcommand{\simto}{\xrightarrow{\sim}}

\newcommand{\define}{\stackrel{\mbox{\scriptsize{def}}}{=}}

\renewcommand{\sl}{\mathfrak{sl}}

\def\bimod{{\text{-bimod}}}

\def\gbimod{{\text{-gbimod}}}

\newcommand{\refequal}[1]{\xy {\ar@{=}^{#1}
(-1,0)*{};(1,0)*{}};
\endxy}

\newcommand{\Hom}{\operatorname{Hom}}

\newcommand{\End}{\operatorname{End}}

\newcommand{\Ext}{\operatorname{Ext}}

\newcommand{\Id}{\operatorname{Id}}
\newcommand{\Tot}{\operatorname{Tot}}

\newcommand{\Sh}{\operatorname{Sh}}
\newcommand{\inv}{^{-1}}

\newcommand{\Ch}{\operatorname{Ch}}
\newcommand{\Cone}{\operatorname{Cone}}

\newcommand{\Bim}{{\rm Bim }}

\newcommand{\SBim}{\SM\Bim}

\newcommand{\FT}{\operatorname{FT}}

\newcommand{\Br}{\operatorname{Br}}

\DeclareMathOperator{\Spec}{Spec}
\DeclareMathOperator{\XSpec}{[\Spec]}

\begin{document}

\begin{abstract} This paper lays the groundwork for the theory of categorical diagonalization. Given a diagonalizable operator, tools in linear algebra (such as Lagrange interpolation)
allow one to construct a collection of idempotents which project to each eigenspace. These idempotents are mutually orthogonal and sum to the identity. We categorify these tools. At the
categorical level, one has not only eigenobjects and eigenvalues but also \emph{eigenmaps}, which relate an endofunctor to its eigenvalues. Given an invertible endofunctor of a
triangulated category with a sufficiently nice collection of eigenmaps, we construct idempotent functors which project to eigencategories. These idempotent functors are mutually
orthogonal, and a convolution thereof is isomorphic to the identity functor.


In several sequels to this paper, we will use this technology to study the categorical representation theory of Hecke algebras. In particular, for Hecke algebras of type A, we will construct categorified Young symmetrizers by simultaneously diagonalizing certain functors associated to the full twist braids.

\end{abstract}

\title{Categorical diagonalization}

\author{Ben Elias}\address{University of Oregon, Eugene}\email{belias@uoregon.edu}\thanks{The first author was supported by NSF grant DMS-1553032, and by the Sloan Foundation}

\author{Matthew Hogancamp}\address{University of Southern California, Los Angeles}\email{hogancam@usc.edu}\thanks{The second author was supported by NSF grants DMS-1255334 and DMS-1702274}

\maketitle

\setcounter{tocdepth}{1}
\tableofcontents


In this paper we initiate the study of the diagonalization of functors. This theory belongs to the realm of \emph{categorification}, which upgrades concepts and constructions in linear algebra to a higher categorical level, replacing linear operators with functors. A key point in the philosophy of categorification (as developed by Chuang-Rouquier \cite{ChuRou}, Khovanov-Lauda \cite{KhoLau09, LauSurvey}, and Rouquier \cite{Rouq2KM-pp}) is that one must study the natural transformations between functors in order to define a theory with desireable
properties. For us this manifests itself in the notion of eigenmaps, which are certain natural transformations relating a functor $F$ with its ``eigenvalues.''

We believe that eigenmaps are a fundamental, yet until now unexplored, component of higher representation theory and categorification. The main goal of this paper is to illustrate this with examples and theory, culminating in our diagonalization theorem. Loosely stated, given an invertible endofunctor $F$ with sufficiently many eigenmaps, our diagonalization theorem
allows one to construct categorical idempotents which project to the eigencategories of $F$.

Let us give a few more details, and ignore any technical requirements. An eigenmap is a natural transformation $\a \co \l \to F$, whose source $\l$ is a ``scalar functor,'' a categorical analog of multiplication by a scalar $\k$. In a triangulated category one can consider $\Cone(\a)$, which is the categorical analog of the operator $(f - \k \cdot 1)$. An eigenobject is an object $M$ for which $\Cone(\a) \ot M \cong 0$, or equivalently, for which $\a \ot \Id_M$ is an isomorphism. (Here we are writing the action of a functor using tensor product notation, as though a monoidal category were acting on itself via the tensor product.) The full subcategory of eigenobjects for $\a$ is the $\a$-eigencategory.

A linear operator $f$ is diagonalizable if there are scalars $\k_i$, $i \in I$ such that \[ \prod_{i \in I} (f - \k_i \cdot 1) = 0. \] The categorical analog of this condition is a structure, that is, a chosen collection of eigenmaps $\a_i$, $i \in I$ such that \[ \bigotimes_{i \in I} \Cone(\a_i) \cong 0. \] Our main theorem, Theorem \ref{thm-firstEigenthm}, takes a nice enough functor $F$ with such a collection of eigenmaps, and explicitly constructs functors $\PB_i$ which project to each $\a_i$-eigencategory.



Our original motivation was to use the theory of diagonalization to develop the categorical representation theory of Coxeter groups and their Hecke algebras. Indeed, in a subsequent paper we categorify the Young symmetrizers in Hecke algebras of type $A$, by diagonalizing the action of Rouquier's full twist complex on the category of Soergel bimodules. The classical Young symmetrizers can be used to construct (and decompose) representations of the symmetric group, and the categorified Young symmetrizers play an analogous role for the categorification. In \S \ref{subsec:Hecke} and \S \ref{subsec:HeckeConj} we discuss the example of Hecke algebras in more detail, and state the main theorems and conjectures of these subsequent papers.

A similar goal would be to categorify then diagonalize the Casimir elements in the centers of quantum groups. The theory we develop in this paper does not apply directly to this case, since the Casimir operators are not invertible. We expect that our theory can be extended to non-invertible operators as well, but we relegate such investigations to future work. For
further discussion, see \S \ref{subsec:casimir}.

Categorified idempotents are also very useful in constructing colored link homology theories. The work presented here explains many observed properties of these idempotents, such as their periodic nature. Our theory of diagonalization produces not just categorified projection operators (which tend to be infinite chain complexes) but also finite relatives which are quasi-idempotent, that is, idempotent up to scalar multiples.  We expect these will yield colored link homology theories which are functorial with respect to link cobordisms. Functoriality gives link homology its 4-dimensional flavor, and is notoriously absent from existing colored homology theories. We mention the applications to link homology briefly in \S \ref{subsec:introTL}.


\subsection{Structure of the paper}

The paper is divided into two parts.  The first is an extended introduction in which we explain, in as elementary terms as possible, the main ideas of eigenmaps and categorical diagonalization, applications, and relations to other work.  The second part can be read independently from the first, and contains all relevant background material, detailed constructions, and proofs.

Part I begins with \S \ref{sec:extendedintro} which summarizes the main features of eigenmaps and categorical diagonalization and states our Diagonalization Theorem.  In \S \ref{sec:subtle} we discuss subtleties of categorical diagonalization which are not visible on the level of linear algebra.  Section \S \ref{sec:introapps} discusses applications of categorical diagonalization.

Part II begins with \S \ref{sec-homotopy}, which contains background on triangulated categories, specifically homotopy categories of complexes, which sets the stage for later constructions.  Then in \S \ref{sec-decompositions} we discuss idempotent decompositions of categories (experts in triangulated categories will recognize these as generalizations of semi-orthogonal decompositions).  Next, \S \ref{sec:prediag} concerns categorical pre-diagonalization and diagonalization in general.  In \S \ref{sec:interpolation} we define categorifications of the Lagrange intperpolation polynomials and some structural results regarding these.  In \S \ref{sec:diag} we prove the Diagonalization Theorem, which states that under certain hypotheses the categorified interpolation polynomials are the eigenprojections in a categorical diagonalization of a functor $F$.   In \S \ref{sec:generalizations} we consider potential generalizations of eigenmaps, which are relevant for the categorified Casimir operator of $\sl_2$.  Finally the appendix \S \ref{sec:appendix} proves some results on commutativity in triangulated monoidal categories.



\begin{remark} It should be mentioned that an earlier draft of this paper was shared with the authors of \cite{GNR}, a truly beautiful work which, among other things, takes our theory of categorical diagonalization and makes it heavily geometric. It is somewhat irresponsible to call upon their work in this paper, which is supposed to be a logical antecedent of theirs. However, they introduce an example (see \S\ref{subsec:cohp1}) which is utterly perfect for describing one of the thorniest issues in diagonalization, so we can not resist using it to make the exposition more clear. We give a brief introduction to parts of \cite{GNR} in \S \ref{subsec:introHilb}. \end{remark}

\subsection{Acknowledgements}
\label{subsec:ack}

The authors would like to thank E.~Gorsky, J.~Rasmussen, and A.~Negu\cb{t} for sharing their beautiful work. We would also like to thank M.~Khovanov for encouraging this collaboration at its inception.

\newpage

\part{Extended Introduction}

\section{Summary of the theory}
\label{sec:extendedintro}

\subsection{Diagonalization and interpolation polynomials in linear algebra}
\label{subsec-LAproj}

We begin by recalling some notions from linear algebra. Let $\KM$ be a commutative ring and let $A$ be a $\KM$-algebra, acting on a $\KM$-module $V$.

\begin{remark} The ring $A$ is where we assume that ``operators'' live. The reader can imagine that $A$ is the endomorphism ring of $V$. Most ideas pertaining to diagonalization can be
discussed within the ring $A$ itself, without mentioning $V$. \end{remark}

What does it mean for an operator $f \in A$ to be diagonalizable? We choose to distinguish between the two most common answers.

\begin{definition}\label{def:prediagLinAlg}
An element $f\in A$ is \emph{prediagonalizable} if there are distinct scalars $\{\k_i\}_{i \in I}$ such that
\begin{equation}\label{LA:minPoly}
\prod (f - \k_i 1) = 0.
\end{equation}
Moreover, we assume that the product of $(f - \k_i 1)$ over any proper subset of $I$ is nonzero, so that the product in \eqref{LA:minPoly} is the minimal polynomial of $f$.
\end{definition}

\begin{definition}\label{def:diagLinAlg}
An element $f\in A$ is \emph{diagonalizable} if there exist a finite set $I$, distinct scalars $\{\k_i\} \subset \KM$, and a complete collection of orthogonal idempotents $p_i\in A$ ($i\in I$) such that $(f-\k_i)p_i = 0=p_i(f-\k_i)$ for all $i\in I$.
If none of the idempotents $p_i$ is zero, then $\{\k_i\}_{i\in I}$ is called the \emph{spectrum} of $f$.
\end{definition}

It is an easy exercise to deduce from the definitions that if $f$ is diagonalizable with spectrum $\{\k_i\}$, then $f$ is prediagonalizable. Conversely, if $f$ is prediagonalizable and
$(\k_i-\k_j)$ is invertible in $\KM$ for all $i\neq j$, then it is also diagonalizable (we recall the proof in \S\ref{subsec:intro_catdiagYO}). Philosophically,
prediagonalizability is a condition, and an easy one to check, while diagonalizability is more structural (giving the operators $p_i$), and is the real payoff. The equivalence between
these notions therefore does a lot of work.

We now begin the process of categorifying these ideas.

\begin{remark} There is a third common notion: an operator $f \co V \to V$ is diagonalizable if $V$ has a basis of eigenvectors. Categorification, however, is not well-suited
to discussing bases in general. The Grothendieck group of a category does come equipped with certain natural bases, e.g. the symbols of the simple objects (or projective objects,
standard objects, etcetera), but these bases are typically not eigenbases. Consequently, the basis of eigenvectors is not the most natural concept to categorify. \end{remark}


\subsection{Categorifying rings, modules, and scalars}
\label{subsec-prelims}

Let $\VC$ be an additive category. Let $\AC$ be an additive monoidal category which acts on $\VC$, in the sense that there is a monoidal functor $\AC \rightarrow \End(\VC)$.
Alternatively, we may use abelian categories or triangulated categories instead of additive categories, in which case we assume that functors in the image of $\AC$ are exact. For more details on triangulated monoidal categories, see \S\ref{subsec:triang}.


\begin{remark} We think of $\AC$ as the world where operators live, and so we often refer to objects in $\AC$ as functors. Very often we will simply assume that $\VC = \AC$, acting on
itself on the left by the monoidal product. \end{remark}

\begin{definition}\label{def:homotopyMonoidalCat}
Let $(\AC,\otimes, \one_{\AC})$ be a triangulated monoidal category.  When we say $\AC$ is a \emph{monoidal homotopy category} we mean that there exists an additive monoidal category $(\BC,\otimes,\one_\BC)$ such that $\AC$ is equivalent to a full, triangulated subcategory of $\KC(\BC)$, with its induced tensor product.  We write this by $\AC \subset \KC(\BC)$. Here, $\KC(\BC)$ is the category of chain complexes in $\BC$, modulo homotopy.
\end{definition}

We are implicitly assuming that every object in $\AC\subset \KC(\BC)$ has appropriate boundedness conditions, so that $C\otimes D$ is guaranteed to exist for all $C,D\in \AC$.  The most common examples will be $\AC=\KC^{\circ}(\BC)$ for $\circ\in \{+,-,b\}$, with $\one_{\AC}=\one_{\BC}$, though not every example is of this form, see Lemma \ref{lem:HMCexample}.

We will restrict to monoidal homotopy categories in the bulk of the paper, to ensure that certain constructions are well-defined.\footnote{Many of our results do work for triangulated monoidal categories in additional generality. The interested reader should consult \S\ref{subsec-nothomotopy}.}

Before categorifying eigenvectors and eigenvalues, one should have a notion of scalars. For this purpose, we fix a monoidal subcategory $\KS \subset \AC$ which plays a role of the
subring $\KM \subset A$. The objects of $\KS$ will be called \emph{scalar objects}\footnote{When we return to scalar objects in \S \ref{subsec:scalars}, $\KS$ will not be a subcategory of $\AC$, but will be a braided monoidal category with a functor to the Drinfeld center $\ZC(\AC)$ of $\AC$.  Consequently, scalar objects enjoy a categorical analogue of the commutativity properties one expects from scalars inside a ring.}, and they give rise to \emph{scalar (endo)functors} of $\VC$. To make things interesting, let us assume
that $\VC$ and $\AC$ come equipped with various grading shift functors (e.g.~$\VC$ could be a graded triangulated category, having a grading shift $(1)$ and a homological shift $[1]$). A
simple but acceptable choice of $\KS$ would say that a scalar object is isomorphic to a direct sum of shifted copies of the monoidal identity $\one$. Scalar functors would then send an
object $M \in \VC$ to a direct sum of shifted copies of $M$.

\begin{remark}
If $\l$ is an invertible scalar functor, then a morphism $\l M\rightarrow N$ will also be regarded as a morphism $M\rightarrow N$ of degree $\l$.  
\end{remark}

\begin{remark} There are numerous contexts where one might want other functors to play the role of scalars (e.g. tensoring with a vector bundle in algebraic geometry). \end{remark}

Note that we have not made any concrete statements about the Grothendieck group. For example, we do not assert that the Grothendieck group of $\AC$ is the ring $A$, or the Grothendieck
group of $\KS$ is the subring $\KM$. This will be a theme in this paper: we are interested in making well-defined categorical constructions, not with ensuring that the Grothendieck
groups of our categories behave nicely. Everything that happens in linear algebra is a motivation for the categorical setting, but we will not try to make the connection precise. See
Remarks \ref{rmk:noGrothGpPlz1} and \ref{rmk:noGrothGpPlz2} for further discussion.

Notationally, we use capital or bold letters to indicate categorical analogs of their lower case friends in linear algebra. For example, a functor $F$ is the categorical analog of a linear operator $f$. The categorical analog of the scalar $\k_i$ will be the scalar object $\l_i$.

\subsection{Eigenobjects and eigenmaps}
\label{subsec:intro_catDiag1}

Now, for a fixed object $F \in \AC$, we ask the question: what might it mean for $F$ to have eigenvectors, or rather, \emph{eigenobjects} in $\VC$?  Here is a naive, and ultimately unsatisfactory answer.

\begin{defn} Let $\l \in \AC$ be a scalar object. An object $M \in \VC$ is a \emph{weak $\l$-eigenobject} for $F \in \AC$ if $F \ot M \cong \l \ot M$. \end{defn}

One reason this is too naive is that, when $\VC$ is triangulated, the full subcategory of weak $\l$-eigenobjects need not be triangulated (i.e. closed under mapping cones), see \S\ref{subsec:nonexamples}. Of course, a genuine categorical notion should involve morphisms; it should somehow fix the isomorphism $F \ot M \cong \l \ot M$ in a natural way. One possibility is below.

\begin{defn} \label{defn:eigenmap} Suppose $\l \in \KS$ is a scalar object and $\a \co \l \to F$ is a morphism in $\AC$. We call a nonzero object $M \in \VC$ an \emph{eigenobject of $F$ with eigenmap $\a$} or an \emph{$\a$-eigenobject} if \[\a \ot \Id_M \co \l \ot M \to F \ot M\] is an isomorphism in $\VC$. A morphism $\a$ as above which has an eigenobject we call a \emph{forward eigenmap}. One might also consider a \emph{backward eigenmap} $\b:F\rightarrow \lambda$, which is defined similarly.

If $\a$ is an eigenmap, define the $\a$-\emph{eigencategory} $\VC_{\a}$ to be the smallest full (additive) subcategory of $\VC$ containing the $\a$-eigenobjects. \end{defn}

Unlike the subcategory of weak eigenobjects, the $\a$-eigencategory is typically well-behaved. We now pass to the triangulated setting, so that the object $\Cone(\a)$ exists.

\begin{proposition}\label{prop:eigenCatsAreTriangulated}
Let $\lambda\in \AC$ be a scalar object, and $\a \co \lambda\to F$ (resp. $\a \co F \to \lambda$) a morphism in $\AC$. Then $M\in \VC$ is an $\a$-eigenobject if and only if $\Cone(\a)\otimes M\simeq 0$.  Consequently, the $\a$-eigencategory is a triangulated subcategory of $\VC$.
\end{proposition}

This is a categorical analogue of the fact that $m\in V$ is a $\k$-eigenvector for $f\in \End(V)$ if and only if $v$ is annihilated by $(f-\k)$.

\begin{proof}
Exactness of the tensor product implies that $\Cone(\a)\otimes M \cong \Cone(\a\otimes \Id_M)$.  It is a standard fact from triangulated categories that $\Cone(\a\otimes \Id_M)\simeq 0$ if and only if $\a\otimes \Id_M$ is an isomorphism in $\VC$.  It is an easy application of the 5-lemma that, for any  $C\in \AC$, the set of objects $M\in \VC$  such that $C\otimes M\simeq 0$ is closed under mapping cones and suspensions.
\end{proof}

Note that there is no reason why a weak eigenobject should admit an eigenmap, and often it does not. What is remarkable is that eigenmaps do occur frequently in nature! Starting in \S\ref{subsec:intro_example1} we give several examples and non-examples.

\begin{remark} Forward and backward eigenmaps are not sufficient for all purposes, e.g. for the diagonalization of the Casimir operator. We postpone this discussion until \S\ref{sec:generalizations}, and focus on forward eigenmaps until then. \end{remark}
	
%

\subsection{Prediagonalizability and commutativity}
\label{subsec:intro_catDiag2}

Henceforth, we restrict to the case where $\AC$ is a monoidal homotopy category, so that we can work with eigencones. Since $\Cone(\a)$ is the categorical analogue of $(f - \k)$, one
might expect the categorical analogue of prediagonalizability (cf. Definition \ref{def:prediagLinAlg}) to appear as follows.

\begin{defn}\label{def:prediagintro} (This definition assumes the vanishing of certain obstructions, c.f. Definition \ref{def:prediag} and see below.)

We say $F\in \AC$ is \emph{categorically prediagonalizable} if it is equipped with a finite set of maps $\{\a_i:\l_i\rightarrow F\}_{i\in I}$ (where $\l_i$ are scalar objects), such that
\begin{enumerate}
\item[(PD1)] $\bigotimes_{i\in I} \Cone(\a_i)\simeq 0$.
\item[(PD2)] the indexing set $I$ is minimal with respect to (PD1).
\end{enumerate}
In this case the collection $\{\a_i\}_{i \in I}$ is called a \emph{saturated collection of eigenmaps} or a \emph{prespectrum} for $F$. \end{defn}

The reason we call $\{\a_i\}$ a prespectrum, rather than the more definitive word ``spectrum," is because it is far from unique. See \S\ref{sec:subtle} for further discussion.

\begin{remark} If $\{\a_i:\l_i\rightarrow F\}_{i\in I}$ is a saturated collection of eigenmaps, then for any given $i \in I$, the tensor product $\bigotimes_{j\neq i}\Cone(\a_j)$ is
nonzero by (PD2), and is an $\a_i$-eigenobject. Thus the maps $\a_i$ are, in fact, eigenmaps. \end{remark}

This straightforward definition is masking a major technical issue. Downstairs, the subring $\KM\langle f \rangle \subset A$ is commutative. However, the full triangulated subcategory of
$\AC$ generated by $F$ and by the scalar functors $\KS$ is not tensor-commutative! In particular, it need not be the case that $\Cone(\a_i) \ot \Cone(\a_j) \cong \Cone(\a_j) \ot
\Cone(\a_i)$ for two eigenmaps $\a_i$ and $\a_j$. Thus in the most general settings the tensor product in (PD1) may depend on the order of the factors.

In \S\ref{subsec:commutativitystuff} we we introduce some homological obstructions to the commutation of cones, leaving proofs to the appendix. There is a primary obstruction, whose vanishing implies that $F \ot \Cone(\a_i) \cong \Cone(\a_i) \ot F$.  Assuming the vanishing of the primary obstructions, there is a secondary obstruction, whose vanishing implies that $\Cone(\a_i) \ot \Cone(\a_j) \cong \Cone(\a_j) \ot \Cone(\a_i)$. When all the secondary obstructions vanish, all eigencones commute, among other useful consequences. Then Definition \ref{def:prediagintro} is equivalent to Definition \ref{def:prediag}, which is our more general definition of prediagonalizability.  We know of no interesting examples where the obstructions do not vanish.

In any case, the proofs of the main theorems all go through without assuming the vanishing of these obstructions.  For more discussion, see \S\ref{subsec:donotcommute}.  For the rest of this introduction, we assume the vanishing of these homological obstructions for simplicity.

\subsection{Example: modules over a cyclic group algebra}
\label{subsec:intro_example1}

We now give an illustrative example, which we will follow throughout the rest of the introduction.

Set $A:=\Z[x]/(x^2-1)$, the group algebra of $\ZM/2\ZM$. Let $\BC$ denote the category of $A$-modules which are free as $\Z$-modules. This category is symmetric monoidal: if $M$, $N$ are
$A$-modules, then $M\otimes_\Z N$ is a $A$-module via $x\cdot(m\otimes n) = (xm)\otimes (xn)$. The monoidal identity is $\Z$, which is a $A$-module by letting $x$ act by 1. We let $\AC$
denote the bounded above homotopy category of $\BC$.

\begin{remark} The category $\BC$ is a toy version of the category $\SBim_2$ of Soergel bimodules in type $A_1$, so this example is actually rather significant. See \S \ref{subsec:typeA1} for the version using Soergel bimodules. \end{remark}

Consider the following complex of $A$-modules:
\[
F = \Big(0 \rightarrow \un{A} \buildrel x-1 \over \longrightarrow A \longrightarrow \Z \rightarrow 0 \Big).
\]
We have underlined the term in homological degree zero. Let us construct two forward eigenmaps for $F$.

The map $\Z \to A$ which sends $1 \mapsto 1 + x$ induces a quasi-isomorphism $\firstMap \co \Z \to F$. The cone of $\firstMap$ is the following acyclic (but non-contractible) complex:
\[ \Cone(\firstMap) = \Big(0 \rightarrow \Z \buildrel x+1 \over \longrightarrow \un{A} \buildrel x-1 \over \longrightarrow A \longrightarrow \Z \rightarrow 0 \Big). \]
We claim that $A$ is a $\firstMap$-eigenobject. Observe that $A \ot A \cong A \oplus A$ as $A$-modules. Consequently
\begin{equation}\label{eq:ACexample}
A\otimes_\Z \Cone(\firstMap) \ \cong \ \Big(0 \rightarrow A \rightarrow (A\oplus A) \rightarrow (A\oplus A) \rightarrow A \rightarrow 0\Big).
\end{equation}
This complex is acyclic since $\Cone(\firstMap)$ is, and consists only of free $A$-modules.  Thus, $A \ot \Cone(\firstMap)$ is contractible (any exact sequence of projective $A$-module is split exact).

The inclusion of the final term $\Z$ induces a chain map $\secondMap \co \Z[-2] \to F$. The complex 
\[
\Cone(\secondMap) \simeq \Big(0 \rightarrow \un{A}  \buildrel x-1 \over \longrightarrow A \rightarrow 0 \Big)
\]
is built entirely from the free module $A$ (as a ``convolution,'' an iterated cone). Consequently, $\Cone(\secondMap) \ot \Cone(\firstMap)$ is built entirely from $A \ot \Cone(\firstMap)$, which is zero. Thus $$\Cone(\secondMap) \ot \Cone(\firstMap) \cong 0.$$
As a consequence $F$ is categorically prediagonalizable with pre-spectrum $\{\firstMap, \secondMap\}$.

\begin{remark} \label{rmk:simsim} We use many arguments of this sort in this paper, where a tensor product of convolutions is simplified by simplifying each piece (in this case, replacing
it by zero) and then re-convolving. We call this \emph{simultaneous simplifications}, see Proposition \ref{prop:simplifications}. A more familiar example of simultaneous simplifications
says that, given a bicomplex where the columns are acyclic, the total complex is also acyclic. This statement is false for arbitrary bicomplexes, but true when the bicomplex satisfies
some boundedness condition. More generally, simultaneous simplifications applies only when certain boundedness conditions hold on the poset which governs the convolution. We provide a
thorough background chapter on convolutions in \S\ref{sec-homotopy} to justify simultaneous simplifications and other similar arguments, and point out where they break down.\end{remark}

\begin{remark} The above examples can be generalized to the group ring of $\ZM/ m \ZM$. One uses the ring $A = \ZM[x]/(x^m-1)$ instead, and replaces $x+1$ with $x^{m-1} + \ldots + x + 1$
everywhere above. Now $A \ot A \cong A^{\oplus m}$, but after minor modifications, one can repeat all the computations above. It is only the $m=2$ case which has a connection to Soergel
bimodules. \end{remark}

Note that $\Cone(\firstMap)$ is a $\secondMap$-eigenobject, though there are no $\secondMap$-eigenobjects which are actually $A$-modules (i.e. complexes supported in homological degree zero). This
motivates further that homotopy categories, rather than additive categories, are the correct setting for diagonalization.

Also note that the scalar objects $\Z$ and $\Z[-2]$ both induce the identity map on the Grothendieck group, so that their eigenspaces are indistinguishable in the Grothendieck group.
However, the $\firstMap$- and $\secondMap$-eigencategories are quite distinct.

\subsection{Related non-examples}
\label{subsec:nonexamples}

Now we discuss several related examples about weak eigenobjects which do not have eigenmaps.

Because $A \ot A \cong A \oplus A$, $A$ is a weak eigenobject for $A$ with eigenvalue $\Z \oplus \Z$. However, there is no forward or backward eigenmap which would realize this
isomorphism. After all, there is only one map $\Z \to A$ up to scalar, so no map $\Z \oplus \Z \to A$ could be injective.

Now consider the complex 
\[
G = \Big(0 \rightarrow \un{A} \longrightarrow \Z \rightarrow 0 \Big).
\]
One can compute that $F$ is homotopic to $G \ot G$, so that one might expect the eigenvalues of $G$ to be the square roots of the eigenvalues of $F$. Note that $G \ot A \simeq A$. Thus $A$ is a weak eigenobject of $G$ with eigenvalue $\Z$ (the square root of the eigenvalue for $\firstMap$). However, there is no corresponding eigenmap $\Z{} \to G$ or $G \to \Z{}$.

In fact, the category of complexes $C$ such that $G\otimes C\simeq C$ is not triangulated. The complex
\[ X = \Big(0 \rightarrow \un{A} \buildrel 1+x \over \longrightarrow A \rightarrow 0 \Big) \]
is sent by $G \ot (-)$ to
\[ Y = \Big(0 \rightarrow \un{A} \buildrel 1-x \over \longrightarrow A \rightarrow 0 \Big), \]
which is not homotopic to the original.

\begin{remark} \label{rem:preunconvincing} Let $\sign$ denote the $A$-module $\Z$, where $x$ acts by $-1$. There \emph{is} a chain map $\sign{}\rightarrow G$ which induces a homotopy
equivalence after tensoring with $A$. Thus, if one regards $\sign$ as a scalar object, then there does exist an eigenmap $\a:\sign\rightarrow G$ such that $A$ is an $\a$-eigenobject.
However, lest the reader become too optimistic, see Remark \ref{rem:unconvincing}. \end{remark}

We believe in the general principle that it is too easy to act diagonalizably on the Grothendieck group, while being categorically diagonalizable requires more structure and is far more
restrictive.

Finally, the complex $F$ is invertible with inverse 
\[
F\inv = \Big(0 \rightarrow \Z \longrightarrow A \buildrel x-1 \over \longrightarrow \un{A} \rightarrow 0 \Big).
\]
That is, $F \ot F\inv \simeq \Z$, the monoidal identity. Tensoring the eigenmaps $\firstMap$ and $\secondMap$ with $F\inv$ and shifting appropriately, one obtains two backward eigenmaps $\g_1 \co F\inv \to \Z$ and $\g_2 \co F\inv \to \Z[2]$.  Since $\Cone(\g_1) \cong \Cone(\firstMap)[2]$ and $\Cone(\g_2) \cong \Cone(\secondMap)[1]$, one still has $\Cone(\g_1) \ot \Cone(\g_2) \simeq 0$.

This example was raised to make the following point. The complex $F$ has two forward eigenmaps, and its inverse has two backward eigenmaps. However, $F$ has no backward eigenmaps, and $F\inv$ has no forward eigenmaps.

\subsection{Eigenmaps are not unique}
\label{subsec:warning}

We feel the need to interject with a significant warning: eigenmaps (for a given eigenobject, or even a given eigencategory) are not unique!

For example, let $\a$ be an eigenmap, and $c$ an invertible scalar in $\Bbbk$, where $\AC$ is a $\Bbbk$-linear category. Then $c\a$ is also an eigenmap, with the exact same
eigencategory. After all, $\Cone(\a) \cong \Cone(c \a)$. In many ways, it is the eigencone, rather than the eigenmap, which is the more canonical object, although the eigenmap is what
gives the precise relationship between $\Cone(\a)$ and the original object $F$.

However, even in important situations, rescaling is not the only way eigenmaps can fail to be unique.  In \S \ref{sec:subtle} we discuss many subtleties of eigenmaps, including examples of non-colinear eigenmaps with the same eigencategory.

\subsection{The goal of categorical diagonalization}
\label{subsec:intro_catdiaggoal}
Let us return to the linear algebra setup from \S\ref{subsec-LAproj}. There we stated that an operator $f \in A$ is prediagonalizable (the product $\prod (f - \k_i)$ vanishes) if and
only if it is diagonalizable (there exist projections $p_i$ to the $\k_i$-eigenspace). More precisely, one should construct $p_i$ such that
\begin{subequations} \label{whatdiagis}
\begin{equation} p_i p_i = p_i, \end{equation}
\begin{equation} p_i p_j = 0, \textrm{ for } i \ne j, \end{equation}
\begin{equation} \sum p_i = 1, \end{equation}
\begin{equation} f p_i = \k_i p_i. \end{equation}
\end{subequations}
Moreover, for any $A$-module $V$ and $v \in V$, $f v = \k_i v$ if and only if $p_i v = v$.

This motivates our definition of a categorically diagonalized functor.

\begin{definition}\label{def:catdiagzed}
Let $F\in \AC$ be an object of a homotopy monoidal category.  Let $I$ be a finite poset, and suppose we are given scalar objects $\l_i\in \AC$, maps $\a_i:\l_i\rightarrow F$, and nonzero objects $\PB_i\in \AC$, indexed by $i\in I$.  We say that $\{(\PB_i,\a_i)\}_{i\in I}$ is a \emph{diagonalization of $F$} if
\begin{itemize}\setlength{\itemsep}{2pt}
\item The $\PB_i$ are mutually orthogonal, in the sense that
\begin{equation} \label{orthogonal} \PB_i \otimes  \PB_j \simeq 0 \textrm{ for } i \ne j. \end{equation} 
\item They are each idempotent, in the sense that there is a homotopy equivalence
\begin{equation} \label{idempotent} \PB_i \otimes \PB_i \simeq \PB_i. \end{equation}
\item There is an idempotent decomposition of identity (see Definition \ref{def-resOfId_intro} below)
\begin{equation} \label{resofid}
\one \simeq \left(\bigoplus_{i\in I} \idemp_i \  ,\  d\right)
\end{equation}
\item $\Cone(\a_i)\otimes \PB_i\simeq 0 \simeq \PB_i\otimes \Cone(\a_i)$.
\end{itemize} 
We say that the diagonalization is \emph{tight} if, whenever $\AC$ acts on a triangulated category $\VC$ and $M$ is an object of $\VC$, one has $\Cone(\a_i) \ot M \simeq 0$ if and only if $\PB_i \ot M \simeq M$.
\end{definition}

The tightness property says that the idempotent $\PB_i$ projects to the $\a_i$-eigencategory; in general it need only project to a subcategory of the $\a_i$-eigencategory. Tightness fails
in important examples, see \S\ref{sec:subtle}.

\begin{remark}
In \S \ref{subsec-diagonalizable} we show that a diagonalized functor satisfies (PD1) in the definition of categorically pre-diagonalizable (Definition \ref{def:prediagintro}), while a tightly diagonalized functor also satisfies (PD2), hence is categorically pre-diagonalizable.
\end{remark}

\begin{defn} \label{def-resOfId_intro} An \emph{idempotent decomposition of} $\one\in \AC$, is a collection of complexes $\idemp_i \in \AC$, indexed by a finite poset $(I,\leq)$,  such that $\idemp_i^{\otimes 2}\simeq \idemp_i$, $\idemp_i\otimes \idemp_j\cong 0$ for $i\neq j$,  and
\begin{equation}\label{eq:introDecomp}
\one \simeq \left(\bigoplus_{i\in I} \idemp_i \  ,\  d\right),
\end{equation}
where $d$ is some lower triangular differential. More precisely, $d=\sum_{i>j} d_{ij}$ where $d_{ij}$ is a degree 1 linear map $\PB_j\rightarrow \PB_i$ and $d_{ii}=d_{\PB_i}$.  This last property describes the monoidal identity $\one$ as homotopy equivalent to a filtered complex whose subquotients are the $\idemp_i$ (more accurately, $\one$ is a convolution with layers $\idemp_i$, as defined in \S\ref{sec-homotopy}). \end{defn}

The main implication of an idempotent decomposition of identity is a canonical filtration of the category, whose subquotients are the images of the idempotents. That is, given a categorically diagonalized functor, one obtains a filtration (more precisely a \emph{semi-orthogonal decomposition}; see \S \ref{subsec:generalDecomp}) of the category whose subquotients are eigencategories (in the tight case). This is in contrast to linear algebra, where one obtains a splitting by eigenspaces rather than a filtration.

The theory of idempotent decompositions of identity is developed at length in work of the second author \cite{Hog17a}, and is recalled in \S\ref{subsec:generalDecomp}.

One key thing to observe is that categorical diagonalization works with a poset of eigenmaps, rather than just a set. This partial order, which governs the filtration on the category
induced by the idempotent decomposition of identity, is a new feature of categorical diagonalization not appearing in ordinary linear algebra; we suggest an interpretation of this
partial order in \S\ref{subsec:spectrum}.

We give some examples which illustrate various subtleties of this definition in \S\ref{sec:subtle}.

\subsection{Categorical diagonalization}
\label{subsec:intro_catdiagYO}

In linear algebra, the eigenprojections $p_i$ can be constructed via an explicit formula which is polynomial in $f$ and rational in the $\k_i$, using Lagrange interpolation. We temporarily treat $f$ as a formal variable, and fix arbitrary scalars $\{\k_i\}$ with $\k_i - \k_j$ invertible for all $i \ne j$. Suppose that $I$ has size $r+1$, and define the $i$-th \emph{Lagrange interpolating polynomial} by
\begin{equation} \label{eq:interPoly}
p_i(f) \define \prod_{j \ne i} \frac{f-\k_j}{\k_i-\k_j}.
\end{equation}
Clearly $p_i(f)\in \KM[f]$ is a polynomial in $f$ of degree $r$.  Furthermore, $p_i(\k_j)=0$ for $j\neq i$, and $p_i(\k_i)=1$.  It follows that for any scalars $\{a_i\}_{i \in I}$ in $\KM$, the expression
\[
L(f):=\sum_i a_i p_i(f)
\]
is the unique polynomial of degree $\leq r$ in $f$ such that $L(\k_i)=a_i$ for all $i$.  In particular if $Q(f)$ is a polynomial of degree $\leq r$, then
\begin{equation}\label{eq:interpolation}
Q(f) = \sum_i Q(\k_i)p_i(f).
\end{equation}
Choosing $Q$ to be the constant polynomial $Q(f)=1$, we see that $1=\sum_i p_i(f) $.

In the concrete setting where $f$ is prediagonalizable and $\prod_i (f-\k_i)=0$, it is an easy exercise to deduce that $p_i:=p_i(f)$ are the eigenprojections for $f$.

Our goal is to take a categorically prediagonalizable functor $F \in \AC$, equipped with eigenmaps $\a_i \co \l_i \to F$, and to categorify \eqref{eq:interPoly} in order to construct
objects $\idemp_i \in \AC$ which project to eigencategories. Without further ado, let us state our first main theorem, which achieves this under some restrictions.

We say that a scalar object $\l$ is \emph{invertible} if there is a scalar object $\l\inv$ with $\l \ot \l\inv \cong \one\cong \l\inv\otimes \l $. We say it is \emph{small} if the infinite direct sum
$\oplus_{n \ge 0} \l^{\ot n}$ exists in $\AC$, and is isomorphic to the infinite direct product. For example, when $\AC$ is the bounded above homotopy category, the homologically-shifted monoidal identity $\one[d]$ is small if and only if $d > 0$. For more discussion of smallness, see \S\ref{subsec-completed}.

\begin{thm}[Diagonalization Theorem]\label{thm:introEigenthm} Let $\AC$ be a monoidal homotopy category (Definition \ref{def:homotopyMonoidalCat}) with category of scalars $\KS$. Fix $F \in \AC$ which is categorically prediagonalizable with pre-spectrum $\{\a_i \co \l_i \to F\}_{i \in I}$. Moreover, assume that: \begin{itemize}
\item Each scalar object $\l_i$ is invertible.
\item The set $I$ is a finite totally-ordered set, which we identify with $\{0,1,\ldots,r\}$. Whenever $i > j$, the scalar object $\l_i \l_j\inv$ is small.
\end{itemize}
Then one can construct objects $\PB_i$ as in Definition \ref{defn:intro_interp} below, called the \emph{categorified interpolation polynomials}, such that $\{(\PB_i,\a_i)\}_{i\in I}$ is a tight diagonalization of $F$ as in Definition \ref{def:catdiagzed}. \end{thm}

Our construction of the categorified interpolation polynomials will mimic \eqref{eq:interPoly}, and can be found in the next section. Note that the complexes involved will be infinite, even when the complex $F$ is finite; this is the expected behavior.

\begin{remark} It would be very interesting to categorify other applications of Lagrange interpolation; see Conjecture \ref{conj:categorifiedInterpolation}. \end{remark}

When $F$ is an invertible complex in $\AC$, its eigenvalues must be indecomposable, and typically the indecomposable scalar objects are invertible. If the eigenvalues have distinct
homological shifts, we can sort them by the homological shift, and apply Theorem \ref{thm:introEigenthm} within the bounded above homotopy category. Thus this theorem applies to
invertible complexes with homologically distinct invertible eigenvalues, such as $F$ from \S\ref{subsec:intro_example1}.

There are many categorically prediagonalizable (invertible) functors to which this theorem does not apply, for instance because the scalar objects $\l_i$ may not be distinct.
Consequently, it is not always possible to place a total order on the indexing set $I$ as in Theorem \ref{thm:introEigenthm}. Examples of this can be found later in this introduction,
such as $\OC(1)$ in \S\ref{subsec:cohp1}, or $\FT_n$ from \S\ref{subsec:HeckeConj} for $n \ge 6$. This situation is extremely interesting, and we discuss it further in \S\ref{sec:subtle}.

\begin{remark} \label{rmk:eigenoverlap} For example, suppose that $\l_i \cong \l_j$. It is possible that there is an object $M$ which is both an $\a_i$-eigenobject and an
$\a_j$-eigenobject. In this case, it is impossible that a diagonalization is tight. This would imply that $\PB_i \ot M \cong M$ and $\PB_j \ot M \cong M$, but then $M \cong \PB_i \ot
\PB_j \ot M \cong 0$. See \S\ref{subsec:cohp1} for an interesting example. \end{remark}

It would be exciting to categorify other aspects of Lagrange interpolation. We conclude this section with an intriguing conjecture.

\begin{conj}[Categorified Lagrange interpolation]\label{conj:categorifiedInterpolation}
Let $\a_i:\l_i\rightarrow F$, $i \in \{0, \ldots, r\}$, be maps such that $\l_i$ is an invertible scalar object and $\l_i\l_j\inv$ is small whenever $i>j$. Let $\PB_i(F)$ denote the interpolation complexes.  Then for each $k\in \Z$ there is a convolution
\[
L_k(F) = \Tot^\oplus\Big(
\begin{tikzpicture}[baseline=-.12cm]
\node (aa) at (0,0){$\l_0^k \PB_0(F)$};
\node (ab) at (3.3,0){ $\l_{1}^k \PB_{1}(F) $};
\node (ac) at (6.3,0){$\cdots$};
\node (ad) at (8.5,0){$\l_r^k \PB_r(F) $ };
\path[->,>=stealth',shorten >=1pt,auto,node distance=1.8cm,font=\small]
(aa) edge node {$[1]$} (ab)
(ab) edge node {$[1]$} (ac)
(ac) edge node {$[1]$} (ad);
\end{tikzpicture}
\Big)
\]
so that  $L_k\simeq F^{\otimes k}$ for $0\leq k\leq r$.
\end{conj}

Our main theorem proves the case $k=0$ when the $\a_i$ are a saturated collection of eigenmaps.

\subsection{Categorified quasi-idempotents}
\label{subsec:introQuasiidemp}

Before summarizing the construction of the projectors in Theorem \ref{thm:introEigenthm}, we discuss a recurring theme in the categorification of idempotents.  

Recall that, decategorified, the eigenprojections of a diagonalizable operator $f$ with eigenvalues $\k_1,\ldots,\k_r$ are given by $p_i=\prod_{j\neq i}\frac{f-\k_j}{\k_i-\k_j}$.  Categorifying denominators is tricky business, so let us consider for the moment only the numerator $k_i:=\prod_{j\neq i} (f-\k_j)$.  Being a scalar multiple of an idempotent, this element is \emph{quasi-idempotent} in the sense that $k_i^2$ is a scalar multiple of $k_i$.

Now we consider a categorical analaogue of this.  Suppose that $\{\a_i:\l_i\rightarrow F\}_{i\in I}$ is a saturated collection of eigenmaps for $F$. Set
\begin{equation} \label{eq:Ki}
\KB_i := \bigotimes_{j \ne i} \Cone(\a_j),
\end{equation}
which is a categorical analogue of $k_i$ considered above.  In \S \ref{subsec:quasiIdemp1} we give a precise statement and prove the following.

\begin{theorem}
Assuming certain obstructions vanish, $\KB_i$ is a quasi-idempotent:
\[
\KB_i\otimes \KB_i \cong \bigotimes_{j\neq i}(\l_i\oplus \l_j[1])\KB_i.
\]
Further, $\KB_i\otimes \Cone(\a_i)\cong 0\cong \Cone(\a_i)\otimes \KB_i$.
\end{theorem}

Given $F$, we thus factor the problem of categorical diagonalization of $F$ into two halves:
\begin{enumerate}\setlength{\itemsep}{3pt}
\item construct a saturated collection of eigenmaps $\a_i:\l_i\rightarrow F$ ($i\in I$).
\item construct the eigenprojections $\PB_i$ via a categorical analogue of dividing $\KB_i$ by $\prod_{j\neq i}(\k_i-\k_j)$.
\end{enumerate}
In this paper we solve (2), under the hypotheses of Theorem \ref{thm:introEigenthm}.  Essentially, the relationship between $\PB_i$ and $\KB_i$ is a dg version of the Koszul duality between polynomial and exterior algebras, which we discuss in subsequent sections.


\subsection{Interpolation complexes}
\label{subsec:interpolation}

Now we discuss the details of our construction of the categorical eigenprojections.  We have categorical analogs of $(f - \k_i)$, namely $\Cone(\a_i)$. What we need next is a categorical analog of $$p_i = \prod_{j \ne i} \frac{f-\k_j}{\k_i - \k_j}.$$ Let \begin{equation} \label{eqn:cji} c_{j,i} = \frac{f-\k_j}{\k_i - \k_j}.\end{equation} We will define complexes $\Cij{j}{i}$ which are categorical analogs of $c_{j,i}$, and define $\idemp_i$ by
\begin{equation} \label{eq:defpi} \idemp_i = \bigotimes_{j \ne i} \Cij{j}{i}. \end{equation} The fact that $c_{j,i}$ is a multiple of $(f - \k_j)$ indicates that $\Cij{j}{i}$ should be ``built from'' copies of $\Cone(\a_j)$, but the crucial idea comes from understanding the denominator of $c_{j,i}$.
	
Let us expand $c_{j,i}$ as a power series in $\k_j \k_i\inv$. We have
\begin{equation} c_{j,i} = \k_i\inv (f - \k_j) (1 + \k_j \k_i\inv + (\k_j \k_i\inv)^2 + \ldots). \end{equation}
This power series would naively categorify to the infinite direct sum $\bigoplus_{n \ge 0} (\l_j \l_i\inv)^{\ot n}$. Note that this infinite direct sum behaves well when $j > i$, by the smallness assumption in Theorem \ref{thm:introEigenthm}. So we might naively construct $\Cij{j}{i}$, when $j > i$, as the infinite direct sum $\bigoplus_{n \ge 0} \l_i\inv \Cone(\a_j) (\l_j \l_i\inv)^{\ot n}$. Similarly, when $i > j$, we should expand $c_{j,i}$ as a power series in $\k_i \k_j\inv$.
\begin{equation} c_{j,i} = -\k_j\inv (f - \k_j) (1 + \k_i \k_j\inv + (\k_i \k_j\inv)^2 + \ldots). \end{equation}
Then, since $\l_i \l_j\inv$ is small, we might naively construct $\Cij{j}{i}$ as the corresponding infinite direct sum, with a homological shift to account for the minus sign.

In truth, $\Cij{j}{i}$ will have the same underlying objects (in each homological degree) as the infinite direct sum, but will have a more interesting differential.

\begin{defn} \label{def:CijIntro} Let $\a \co \l \to F$ and $\b \co \mu \to F$ be maps from scalar objects $\l$ and $\mu$, and assume that $\l \mu\inv$ is small. Define $\Cij{\a}{\b}$ as the total complex of the following convolution diagram.
\begin{equation}\label{eq:introCab} \Cij{\a}{\b} = \left(
\begin{tikzpicture}[baseline=1em]
\node(aa) at (0,0) {$\displaystyle \frac{1}{\mu} F$};
\node(ba) at (2.5,1) {$\displaystyle \frac{\l}{\mu} $};
\node(ca) at (5,0) {$\displaystyle \frac{\l}{\mu^2} F$};
\node(da) at (7.5,1) {$\displaystyle \frac{\l^2}{\mu^2} $};
\node(ea) at (10,0) {$\cdots$};
\tikzstyle{every node}=[font=\small]
\path[->,>=stealth',shorten >=1pt,auto,node distance=1.8cm]
(ba) edge node[above] {$ \frac{1}{\mu} \a$} (aa)
(ba) edge node[above] {$ -\frac{\l}{\mu^2} \b $} (ca)
(da) edge node[above] {$ \frac{\l}{\mu^2}\a$} (ca)
(da) edge node[above] {$- \frac{\l^2}{\mu^3}\b$} (ea);
\end{tikzpicture} \right).
\end{equation}
Similarly, $ \Cij{\b}{\a}(F)[1]$ (note the homological shift) is the total complex of the following convolution diagram:
\begin{equation}\label{eq:introCba} \Cij{\b}{\a} = \left(
\begin{tikzpicture}[baseline=1em]
\node(za) at (-2.5,1) {$\one $};
\node(aa) at (0,0) {$\displaystyle \frac{1}{\mu} F$};
\node(ba) at (2.5,1) {$\displaystyle \frac{\l}{\mu} $};
\node(ca) at (5,0) {$\displaystyle \frac{\l}{\mu^2} F$};
\node(da) at (7.5,1) {$\displaystyle \frac{\l^2}{\mu^2} $};
\node(ea) at (10,0) {$\cdots$};
\tikzstyle{every node}=[font=\small]
\path[->,>=stealth',shorten >=1pt,auto,node distance=1.8cm]
(za) edge node[above] {$-\frac{1}{\mu} \b$} (aa)
(ba) edge node[above] {$ \frac{1}{\mu} \a$} (aa)
(ba) edge node[above] {$ -\frac{\l}{\mu^2} \b $} (ca)
(da) edge node[above] {$ \frac{\l}{\mu^2}\a$} (ca)
(da) edge node[above] {$ -\frac{\l^2}{\mu^3}\b$} (ea);
\end{tikzpicture} \right).
\end{equation}
\end{defn}

For example, omitting the southeast-pointing arrows from \eqref{eq:introCab}, one obtains the naive direct sum $\bigoplus_{n \ge 0} \mu\inv \Cone(\a) (\l \mu\inv)^{\ot n}$.

\begin{remark} \label{rmk:signsIntro} The signs on the southwest-pointing arrows in the definition above are conventional. The complexes $\Cij{\a}{\b}$ and $\Cij{c\a}{d\b}$ are isomorphic, for any invertible scalars $c$ and $d$, so these signs can be removed up to isomorphism. The utility of this particular sign convention will become clear in Remark \ref{rmk:signsIntro2}. \end{remark}

\begin{defn} \label{defn:intro_interp} Set $\Cij{j}{i} = \Cij{\a_j}{\a_i}$, which is defined using \eqref{eq:introCab} when $j > i$ and \eqref{eq:introCba} when $i > j$. Define $\idemp_i$ using \eqref{eq:defpi}. \end{defn}

The proof that these complexes $\idemp_i$ satisfy the properties of Theorem \ref{thm:introEigenthm} is the topic of \S\ref{sec:diag}. In \S\ref{subsec:intro_Koszul} we motivate this
definition of $\Cij{\a}{\b}$ using Koszul duality.

\begin{remark} Many readers may be unfamiliar with convolutions and their total complexes, like the one in \eqref{eq:introCab} or the one in \eqref{resofid}. Given an appropriate diagram of
complexes, one can take an iterated cone to obtain its total complex. Again, a familiar example is a bicomplex (which could be thought of as a complex of complexes) and its total
complex. Convolution diagrams are indexed by a partially ordered set (e.g. $\Z$ for a bicomplex), though the literature tends to focus on totally ordered sets or on finite sets. We have
included \S\ref{sec-homotopy} to give a thorough background on convolutions, including those governed by more interesting posets. \end{remark}

\begin{remark} \label{rmk:zigzag} In \S\ref{subsec:intro_catDiag2} we mentioned some homological obstructions to the commutation of eigencones discussed in the appendix. A convolution
like \eqref{eq:introCab} is called a \emph{zigzag diagram}. We prove in the appendix that the same obstructions govern whether zigzag diagrams tensor-commute. Thus, under some technical
assumptions, the complexes $\Cij{j}{i}$ and $\idemp_i$ will also commute up to homotopy equivalence. Hence \eqref{eq:defpi} is independent of the order of the tensor product. \end{remark}

\subsection{Interpolation complexes for the cyclic group algebra}
\label{subsec:interpolation_example}

Let us continue the notation of \S\ref{subsec:intro_example1}. Our eigenmaps were called $\firstMap$ and $\secondMap$, with associated scalar objects $\firstShift = \Z$ and $\secondShift = \Z[-2]$. Consequently, $\firstShift \secondShift\inv$ is small in the bounded above homotopy category. Note that there are only two eigenvalues, so that $\Cij{j}{i} = \idemp_i$.

Recall that
\[ \Cone(\firstMap) = \Big(0 \rightarrow \Z \buildrel x+1 \over \longrightarrow \un{A} \buildrel x-1 \over \longrightarrow A \longrightarrow \Z \rightarrow 0 \Big). \]
The complex $\secondProj = \Cij{\firstMap}{\secondMap}$ is built using \eqref{eq:introCab} by gluing together infinitely many copies of $\Cone(\firstMap)$, where the initial $\Z$ in one copy is sent to the final $\Z$ in the next copy.
\[
\begin{tikzpicture}[baseline=1em]
\node(aa) at (0,0) {$\Z$};
\node(ba) at (-2,0) {$A$};
\node(ca) at (-4,0) {$A$};
\node(da) at (-6,0) {$\Z$};
\node(ab) at (-4,-1.5) {$\Z$};
\node(bb) at (-6,-1.5) {$A$};
\node(cb) at (-8,-1.5) {$A$};
\node(db) at (-10,-1.5) {$\Z$};
\node(ac) at (-8,-3) {$\Z$};
\node(bc) at (-10,-3) {$A$};
\node(cc) at (-12,-3) {$\cdots$};
\node at (-4,-.75) {$\oplus$};
\node at (-6,-.75) {$\oplus$};
\node at (-8,-2.25) {$\oplus$};
\node at (-10,-2.25) {$\oplus$};
\tikzstyle{every node}=[font=\small]
\path[->,>=stealth',shorten >=1pt,auto,node distance=1.8cm]
(da) edge node[above] {} (ca)
(ca) edge node[above] {} (ba)
(ba) edge node[above] {} (aa)
(db) edge node[above] {} (cb)
(cb) edge node[above] {} (bb)
(bb) edge node[above] {} (ab)
(cc) edge node[above] {} (bc)
(bc) edge node[above] {} (ac)
(da) edge node[above=4pt] {$-1$}  (ab)
(db) edge node[above=4pt] {$-1$}  (ac);
\end{tikzpicture}
\]
Applying Gaussian elimination to remove the contractible copies of $\Z\buildrel -1\over\longrightarrow \Z$, one obtains the following complex:
\begin{equation}
\secondProj = \Big(\cdots \buildrel x-1\over \longrightarrow A \buildrel x+1 \over\longrightarrow A  \buildrel x-1 \over \longrightarrow A\rightarrow \underline{\Z} \rightarrow 0\Big),
\end{equation}
in which which we have underlined the term in homological degree zero.  A similar analysis computes $\firstProj = \Cij{\secondMap}{\firstMap}$ built using \eqref{eq:Cba} to be the following:
\begin{equation}
\firstProj = \Big(\cdots \buildrel x-1\over \longrightarrow A \buildrel x+1 \over\longrightarrow A  \buildrel x-1 \over \longrightarrow \underline{A} \rightarrow 0\Big).
\end{equation}

Let us argue that $\firstProj \ot \secondProj \cong 0$. Recall that $A \ot \Cone(\firstMap) \cong 0$ from \S\ref{subsec:intro_example1}. Applying simultaneous simplifications (see Remark
\ref{rmk:simsim}) we deduce that $A \ot \secondProj \cong 0$, because $\secondProj$ is built from copies of $\Cone(\firstMap)$. Now $\firstProj$ is built from copies of $A$, so we again
apply simultaneous simplifications and deduce that $\firstProj \ot \secondProj \cong 0$.

\begin{remark} The convolution building $\firstProj$ out of $A$ is governed by a poset with the ascending chain condition (ACC), while the one building $\secondProj$ out of $\Cone(\firstMap)$ has the DCC. It is our smallness hypothesis which allows one to use simultaneous simplifications in both situations. See \S\ref{sec-homotopy} for more details. \end{remark}

It is easy to see that there is a chain map $\secondProj[-1] \to \firstProj$, whose cone is ${\Z}$, the monoidal identity. This description of ${\Z}$ as a cone is precisely the
(idempotent) decomposition of identity from \eqref{resofid} in Theorem \ref{thm:introEigenthm}.

The fact that $\secondProj \ot \secondProj \cong \secondProj$ can now be deduced in several ways.  For instance, one can use simultaneous simplifications and the fact that $A \ot \secondProj \cong 0$ to see that the only
surviving term in the tensor product $\secondProj \ot \secondProj$ is ${\Z} \ot \secondProj$. Equivalently, one has $$\secondProj \cong \secondProj \ot {\Z} \cong \secondProj \ot \Cone(\secondProj[-1] \to
\firstProj) \cong \secondProj \ot \secondProj,$$ where the last equality follows because $\secondProj \ot \firstProj \cong 0$. This line of argument also proves easily that $\firstProj \ot \firstProj \cong \firstProj$.

It is easy to complete the proof of Theorem \ref{thm:introEigenthm} in this example.

\begin{remark}
We refer to $\firstProj$ as a \emph{counital idempotent}, since it is equipped with a chain map $\firstProj\rightarrow \one$ which becomes a homotopy equivalence after tensoring on the right or left with $\firstProj$.  Similarly, $\secondProj$ is a \emph{unital idempotent}. 
\end{remark}

\begin{remark} Our arguments above are fairly general. Under the hypotheses of Theorem \ref{thm:introEigenthm}, it is fairly straightforward to prove that $\PB_i \ot \PB_j \cong 0$ for
$i \ne j$ using (PD1) and simultaneous simpifications. It follows easily that any convolution $\one\simeq \Tot\{\PB_i,d_{ij}\}$ is automatically an idempotent decomposition of identity.
Consequently, the most interesting part of the proof of Theorem \ref{thm:introEigenthm} is the construction of such a convolution. \end{remark}

\begin{remark} If instead we work in the bounded below homotopy category, it is $\secondShift \firstShift\inv$ which is small (so the eigenmaps should now be ordered so that $\secondMap>\firstMap$). One must now define $\Cij{\firstMap}{\secondMap}$ using \eqref{eq:introCba} instead of \eqref{eq:introCab}. We call the corresponding idempotent complexes $\firstProj^\vee$ and $\secondProj^\vee$.
\[
\secondProj^\vee \ = \ \Big(0\rightarrow \un{\Z} \rightarrow A \buildrel x-1 \over  \longrightarrow A \buildrel x+1 \over  \longrightarrow A \buildrel x-1 \over  \longrightarrow \cdots \Big).
\]
\[
\firstProj^\vee \ = \ \Big(0\rightarrow \un{A} \buildrel x-1 \over  \longrightarrow A \buildrel x+1 \over  \longrightarrow A \buildrel x-1 \over  \longrightarrow \cdots \Big).
\]
To prove that $\firstProj^\vee \ot \secondProj^\vee \cong 0$, we use simultaneous simpifications as above, this time relying on the fact that the convolutions satisfy a descending chain condition. Now we have an isomorphism ${\Z} \cong \Cone(\firstProj^\vee[-1] \to \secondProj^\vee)$, as desired. \end{remark}

%

\subsection{Koszul duality and interpolation complexes}
\label{subsec:intro_Koszul}

We now wish support our claim that the categorified interpolation complexes $\PB_i$ and the finite complexes $\KB_i:=\bigotimes_{j\neq i} \Cone(\a_j)$ are related by Koszul duality.  We first recall the classical Koszul duality between polynomial and exterior algebras.

Let $R = \C[u]$ denote the polynomial ring in one variable, $\L = \Lambda^\ast(\pa)$ the exterior algebra in one variable, and $\C$ the one-dimensional graded module over either ring.  Let $\l=(1)$ be the grading shift functor.  We say
that the degree of $u$ is $\l^{-1}=(-1)$, meaning that the map $R(-1) \to  R$ given by multiplication by $u$ is homogeneous of degree zero.

Let $\CC$ be the category of finitely generated, graded $R$-modules.  Let $V'=K_0(\CC)$ denote the (abelian) Grothendieck group of $\CC$.  Note that $V'$ is a $\Z[\k,\k\inv]$-module, via $\k[M] = [M(-1)]$.  Let $V:=\Z[\k,\k\inv, (1-\k)\inv]\otimes_{\Z[\k,\k\inv]} V'$.

Observe that there is a short exact sequence
\[
0\rightarrow R(-1)\buildrel u\over \rightarrow R \rightarrow \C\rightarrow 0,
\]
which categorifies the identity $[\C]=(1-\k)[R]$ in $V$.  Consequently, $[R] = \frac{1}{1-\k}[\C]$ in $V$.  This  latter identity is categorified by the fact that $R$ admits a filtration with associated graded $\bigoplus_{k\geq 0} \C(-k)$.

It is instructive to observe that, in the derived category of $\CC$, the short exact sequence $0\rightarrow R(-1)\rightarrow R\rightarrow \C\rightarrow 0$ corresponds to a quasi-isomorphism
\begin{equation}\label{eq:CintermsOfR}
\Cone(R(-1)\buildrel u\over \rightarrow R) \cong \C.
\end{equation}
Heuristically, one regards $R$ as built from copies of $\C$ glued in a ``polynomial algebra configuration,'' while $\C$ is built from copies of $R$ glued
in an ``exterior algebra configuration.''

Similarly, working inside the category of graded $\L$-modules, we see that $\L$ has a filtration by copies of $\C$, which categorifies the statement $[\L] = (1 + \k) [\C]$.  Conversely, the identity $[\C]\simeq \frac{1}{1+\k}[\L]$ is categorified by the fact that $\C$ is quasi-isomorphic to the complex
 \begin{equation}\label{CintermsofLambda}
\C \cong  \Big( \ldots \to \L(-2) \to \L(-1) \to \underline{\L} \Big),
 \end{equation}
 where the maps are multiplication by $\pa$.  Heuristically, one regards $\C$ as built from copies of $\L$ glued in a polynomial algebra configuration, while $\L$ is built from copies of $\C$ glued in an
exterior algebra configuration.  In certain contexts, $\partial$ should be thought of as carrying some homological degree, in which case the expressions $1+\k$ above get replaced by $1-\k$.

Koszul duality relates these two observations.  Applying the functor $\Ext_R^*(-,\C)$ sends $\C$ to $\L$. The element $\partial \in \Ext_R^1(\C,\C)$ is the evident endomorphism of the left-hand side of \eqref{eq:CintermsOfR}, which has degree $(1)[-1]$.  Also, $\Ext_R^\ast(-,\C)$ sends $R$ to $\C$, and sends any complex of graded $R$-modules to a complex of graded $\L$-modules. Conversely $\Ext^*_\L(-,\C)$ sends $\L$ to $\C$, $\C$ to $R$, and any complex of graded $\L$-modules gets sent to a complex of graded $R$-modules.   Informally, Koszul duality lifts mutiplication and division by $1-\k$ to a pair of inverse autoequivalences relating the derived categories of $\L$ and $R$.

Now let us return to the case at hand.  Let $\l,\mu$ be invertible scalar objects such that $\l\mu\inv$ is small, and let $\a:\l\rightarrow F$ and $\b:\mu\rightarrow F$ be given.  We will analyze the relationship between $\Cij{\a}{\b}$ and $\Cone(\a)$ with a lens similar to the above.  As $\Cij{\a}{\b}$ is built from copies of $\Cone(\a)$ glued in a polynomial algebra configuration, we expect $\Cij{\a}{\b}$ to be analogous to $R$ inside $R$-modules or $\C$ inside $\L$-modules, while we expect $\Cone(\a)$ to be analogous to $\C$ inside $R$-modules or
$\L$ inside $\L$-modules.  Indeed, we claim the following:
\begin{itemize}
\item $\Cij{\a}{\b}$ has an endomorphism $u$ of degree $\l \mu\inv$ (the ``periodicity map''), analogous to multiplication by $u$ in $R$.
\item Using $u$, we may construct a complex of the form
\[
\Big(\begin{diagram}\frac{\l}{\mu}\Cij{\a}{\b}[1] &\rTo^{u} &\Cij{\a}{\b}\end{diagram}\Big)
\]
which is homotopy equivalent to $\Cone(\b)$, in analogy with (\ref{eq:CintermsOfR}).
\item $\Cone(\a)$ has an endomorphism $\pa$ of degree $\l\inv \mu[-1]$, satisfying $\pa^2 = 0$, analogous to multiplication by $\pa$ in $\L$.
\item Using $\partial$, we may construct a complex of the form
\[
\Big( \begin{diagram}\Cone(\a) &\rTo^{[1]}& \l\mu\inv \Cone(\a) &\rTo^{[1]}&  \l^2\mu^{-2}\Cone(\a) &\rTo^{[1]}&  \cdots \end{diagram} \Big)
\]
which is homotopy equivalent to $\Cij{\a}{\b}$, in analogy with (an alternate version of) (\ref{CintermsofLambda}).
\end{itemize}

The map $\partial$ is constructed using the map $\b$:
\[
\begin{tikzpicture}[baseline=-4em]
\node at (4.9,0) {$\Big)$};
\node at (.7,0) {$\Big($};
\node at (1.3,-2) {$\Big)$};
\node at (-2.4,-2) {$\Big($};
\node at (-3.5,0) {$=$};
\node at (-3.5,-2) {$=$};
\node(f) at (-5.7,0) {$\displaystyle \frac{\mu}{\l}[-1] \Cone(\a)$};
\node(h) at (-5.7,-2) {$\Cone(\a)$};
\node(a) at (4,0) {$\displaystyle\frac{\mu}{\l} F[-1]$};
\node(b) at (1,0) {$\mu $};
\node(c) at (1,-2) {$F$};
\node(d) at (-2,-2) {$\l[1] $};
\tikzstyle{every node}=[font=\small]
\path[->,>=stealth',shorten >=1pt,auto,node distance=1.8cm]
(b) edge node[above] {$\a $} (a)
(b) edge node[right=3pt] {$-\b$} (c)
(d) edge node[above] {$\a$} (c)
(f) edge node[right] {$\partial$} (h);
\end{tikzpicture}.
\]
The construction of $\Cij{\a}{\b}$ using $\Cone(\a)$ and $\partial$ is immediate from the construction in diagram \eqref{eq:introCab}.  Finally, the periodicity map $u$ is simply the chain map $\frac{\l}{\mu}\Cij{\a}{\b}\rightarrow \Cij{\a}{\b}$ which shifts the zigzag diagram \eqref{eq:introCab} to the right.

There are similar observations concerning $\Cij{\b}{\a}$ and its relationship with $\Cone(\b)$.  One difference is that the periodicity map $\frac{\mu}{\l}\Cij{\b}{\a}$ shifts the zigzag diagram \eqref{eq:introCba} to the left.

\begin{example} We illustrate this in the example from \S\ref{subsec:interpolation_example}.

Let us temporarily denote the periodicity map of $\firstProj$ by $u_{\secondMap,\firstMap}$, and the periodicity map of $\secondProj$ by $u_{\firstMap,\secondMap}$. The degree of $u_{\secondMap,\firstMap}$ is $\one[-2]$, and that of $u_{\firstMap,\secondMap}$ is $\one[2]$. The reader should confirm that these maps look as follows.
\[
\begin{tikzpicture}[baseline=-.12cm]
\node (aa) at (0,0){$\firstProj$};
\node at (1.5,0) {$=$};
\node (ab) at (3,0){ $\cdots$};
\node (ac) at (5,0){$A$};
\node (ad) at (7,0){$A$};
\node (ae) at  (9,0){$A$};
\node (af) at (11,0){$0$};
\node (ag) at (13,0){$0$};
\node (ba) at (0,2){$\firstProj[-2]$};
\node at (1.5,2) {$=$};
\node (bb) at (3,2){ $\cdots$};
\node (bc) at (5,2){$A$};
\node (bd) at (7,2){$A$};
\node (be) at  (9,2){$A$};
\node (bf) at (11,2){$A$};
\node (bg) at (13,2){$A$};
\path[->,>=stealth',shorten >=1pt,auto,node distance=1.8cm,font=\small]
(ab) edge node {$x-1$} (ac)
(ac) edge node {$x+1$} (ad)
(ad) edge node {$x-1$} (ae)
(ae) edge node {} (af)
(af) edge node {} (ag)
(bb) edge node {$x-1$} (bc)
(bc) edge node {$x+1$} (bd)
(bd) edge node {$x-1$} (be)
(be) edge node {$x+1$} (bf)
(bf) edge node {$x-1$} (bg)
(ba) edge node {$u_{\secondMap,\firstMap}$} (aa)
(bc) edge node {$\Id$} (ac)
(bd) edge node {$\Id$} (ad)
(be) edge node {$\Id$} (ae)
(bf) edge node {} (af)
(bg) edge node {} (ag);
\end{tikzpicture}
\]
\[
\begin{tikzpicture}[baseline=-.12cm]
\node (aa) at (0,0){$\secondProj[2]$};
\node at (1.5,0) {$=$};
\node (ab) at (3,0){ $\cdots$};
\node (ac) at (5,0){$A$};
\node (ad) at (7,0){$A$};
\node (ae) at  (9,0){$\Z$};
\node (af) at (11,0){$0$};
\node (ag) at (13,0){$0$};
\node (ba) at (0,-2){$\secondProj$};
\node at (1.5,-2) {$=$};
\node (bb) at (3,-2){ $\cdots$};
\node (bc) at (5,-2){$A$};
\node (bd) at (7,-2){$A$};
\node (be) at  (9,-2){$A$};
\node (bf) at (11,-2){$A$};
\node (bg) at (13,-2){$\Z$};
\path[->,>=stealth',shorten >=1pt,auto,node distance=1.8cm,font=\small]
(ab) edge node {$x+1$} (ac)
(ac) edge node {$x-1$} (ad)
(ad) edge node {} (ae)
(ae) edge node {} (af)
(af) edge node {} (ag)
(bb) edge node {$x+1$} (bc)
(bc) edge node {$x-1$} (bd)
(bd) edge node {$x+1$} (be)
(be) edge node {$x-1$} (bf)
(bf) edge node {} (bg)
(aa) edge node {$u_{\firstMap,\secondMap}$} (ba)
(ac) edge node {$\Id$} (bc)
(ad) edge node {$\Id$} (bd)
(ae) edge node {} (be)
(af) edge node {} (bf)
(ag) edge node {} (bg);
\end{tikzpicture}
\]
It is straightforward to check that $\Cone(u_{\firstMap,\secondMap})\simeq \Cone(\firstMap)$ and $\Cone(u_{\secondMap,\firstMap})\simeq \Cone(\secondMap)$, up to shifts.

The endomorphism $\pa_{\firstMap,\secondMap}$ of $\Cone(\firstMap)$ is
\[
\pa_{\firstMap,\secondMap} = \quad
\begin{diagram}[small]
0 & \rTo & 0 & \rTo & 0 & \rTo & \Z & \rTo^{x+1} & A & \rTo^{x-1} & A & \rTo & \Z  \\
 \dTo &&\dTo && \dTo && \dTo^{\Id} && \dTo && \dTo &&\dTo \\
\Z & \rTo^{x+1} & A & \rTo^{x-1} & A & \rTo & \Z & \rTo & 0 & \rTo & 0 & \rTo & 0 .
\end{diagram}
\]
The endomorphism $\pa_{\secondMap,\firstMap}$ of $\Cone(\a_2)$ is
\[
\pa_{\secondMap,\firstMap} = \quad
\begin{diagram}[small]
0 & \rTo & A & \rTo^{x-1} & A &\rTo & 0 &\rTo & 0 \\
\dTo &&\dTo && \dTo^{x+1} && \dTo &&\dTo  \\
0 & \rTo & 0 & \rTo & A & \rTo^{x-1}& A &\rTo & 0.
\end{diagram}
\]
Here, multiplication by $x+1$ arises from the composition of $A \to \Z$ with $\Z \to A$.

\end{example}

Now, suppose $I=\{0,1,\ldots,r\}$, $\l_i$ are invertible scalars with $\l_i \l_j\inv$ small for $i > j$, and $\a_i:\l_i\rightarrow F$ are maps.  Set $\KB_i = \bigotimes_{j \ne i} \Cone(\a_j)$ as in \eqref{eq:Ki}.  The relation between $\Cij{j}{i}$ and $\Cone(\a_j)$  is essentially a dg version of Koszul duality between polynomial and exterior algebras in one variable, as discussed above.  Taking appropriate tensor products, one sees that the relation between $\PB_i$ and $\KB_i$ is Koszul duality between polynomial and exterior algebras in $r$ variables.  More precisely, each $\Cij{j}{i}$ has an endomorphism (the ``periodicity map'') of degree $\l_j/\l_i$, which we will denote suggestively by $u_j/u_i$, where $u_k$ are formal indeterminates of degree $\l_k$, $0\leq k\leq r$.  Consequently, the tensor product $\PB_i$ has commuting endomorphisms $u_j/u_i$ for each $j \ne i$, which we also call periodicity maps.  Using these endomorphisms, one can form an $r$-dimensional cube-like complex (the usual Koszul complex), in which the edges are $\pm u_j/u_i$, the total complex of which is homotopy equivalent to $\KB_i$.

Conversely, $\KB_i$ is has an action of an exterior algebra in $r$ variables, and $\PB_i$ can be reconstructed as the total complex of a certain infinite, multiply periodic complex constructed from shifted copies of $\KB_i$. We will return to this theme in \S\ref{subsec:compact}.



\subsection{Periodicity maps}
\label{subsec:intro_period}

We have used the notation $u_j/u_i$ for a periodicity map on the projector $\PB_i$. This notation is somewhat misleading, as we have not defined morphisms $u_j$ or $u_i$ which one ``divides'' to obtain $u_j/u_i$. However, under certain assumptions, such morphisms $u_i$ exist: they are the eigenmaps.

Suppose that $M$ is an $\a_i$-eigenobject for $F$. Then there is an isomorphism $(\a_i \ot \Id_M)^{-1} \co F \ot M \to \l_i M$, which we denote $\a_i^{-1}(M)$. In particular, $\PB_i$ is an $\a_i$-eigenobject, so $\a_i^{-1}(\PB_i)$ exists. The eigenmap $\a_j$ also induces a map $\a_j(\PB_i) \co \l_j \PB_i \to F \ot \PB_i$. Composing these maps and shifting appropriately we obtain the map
\begin{equation} \label{eq:ajoverai} \frac{\a_j}{\a_i}(\PB_i) \co \l_j \l_i\inv \PB_i \to \PB_i, \end{equation}
which has the same degree as the periodicity map $u_j/u_i$.

In Proposition \ref{prop:eigenaction} we prove that
\begin{equation} \label{eq:ujuiajai} \frac{u_j}{u_i} = \frac{\a_j}{\a_i}(\PB_i) \end{equation}
when certain homological obstructions vanish. In fact, these are the same homological obstructions discussed earlier (and defined in \S\ref{subsec:commutativitystuff}), which governed the commuting properties of eigencones. Our proof of Proposition \ref{prop:eigenaction} resembles many of the proofs in the Appendix.

\begin{remark} \label{rmk:signsIntro2} The equality \eqref{eq:ujuiajai} relies on our sign convention in \eqref{eq:introCab}, see Remark \ref{rmk:signsIntro}. \end{remark}

\subsection{Simultaneous diagonalization}
\label{subsec:simultaneous}

When $f$ and $g$ are commuting operators, and both are diagonalizable, then they can be simultaneously diagonalized. We now explore what can be done to categorify this fact.

Suppose that $\{\PB_i\}_{i \in I}$ and $\{\QB_j\}_{j \in J}$ are two collections of idempotent functors (indexed by posets $I$ and $J$) which can be arranged into idempotent decompositions of identity, as
in \eqref{def-resOfId_intro}. Suppose that $\PB_i \ot \QB_j \cong \QB_j \ot \PB_i$ for each $i \in I$ and $j \in J$. Then it is relatively straightforward (see Proposition \ref{prop:tensoringResolutions}) to prove
that $\PB_i \ot \QB_j$ is an idempotent, and that tensoring together the two idempotent decompositions of identity yields a third one $\{\PB_i \ot \QB_j\}_{(i,j) \in I \times J}$.

This is an example of refining an idempotent decomposition of identity (it is the mutual refinement of $\{\PB_i\}$ and $\{\QB_j\}$). See \S\ref{subsec-idempotentreassoc} for further discussion of reassociation of idempotent decompositions and refinement.

\begin{remark} Of course, some of the idempotents $\PB_i \ot \QB_j$ may be zero, if the corresponding ``joint eigencategory'' is zero. \end{remark}

Let $F$ and $G$ be two functors in the same monoidal homotopy category $\AC$, which are both categorically diagonalized in the sense of Definition \ref{def:catdiagzed}, with projectors $\PB_i$ and $\QB_j$ respectively.   The complicating issue is that, even if one assumes that $F \ot G \cong G \ot F$, there is no guarantee that $\PB_i \ot \QB_j \cong \QB_j \ot \PB_i$. As usual, commutativity is more complicated on the categorical level.

However, if Theorem \ref{thm:introEigenthm} can be applied to both $F$ and $G$, then we have explicit constructions of $\PB_i$ and $\QB_j$. Using the same results about zigzag complexes mentioned in Remark \ref{rmk:zigzag}, we are able to provide homological obstructions whose vanishing implies that $\PB_i$ and $\QB_j$ commute. This is proven in the appendix, Proposition \ref{prop:zigzagsCommute}. This enables the simultaneous diagonalization of nicely-commuting functors.  We give a more detailed statement in the sequel, since simultaneous diagonalization plays an important role in the applications to Soergel bimodules.

\subsection{Cell filtrations and diagonalization}
\label{subsec:cells}

Let us conclude with a brief description of some technology to be developed in the sequel.

For any diagonalizable functor $F:\VC\rightarrow \VC$ there is a filtration of $\VC$ coming from the idempotent resolution of identity (in the context of triangulated categories this is
a semi-orthogonal decomposition; see Remark \ref{rmk:semiorthogonal}). When $\VC$ is the homotopy category of an additive category $\CC$, one might ask whether the filtration of $\VC$ is
inherited from a filtration on $\CC$. This is often the case: for the Hecke category (see \S\ref{subsec:Hecke}), this filtration is the cell filtration on the monoidal category.

Let us continue the setup of \S\ref{subsec:intro_example1}. That $F$ was categorically prediagonalizable was roughly proven as follows: $\Cone(\firstMap)$ kills $A$ so it kills any complex built out of the object $A$, while $\Cone(\secondMap)$ is built from the object $A$. Let us rephrase this argument.  Inside the category of $A$-modules, we have the full subcategory $\JC$ consisting of direct sums of $A$. This category is closed under tensor product with any object; it is a \emph{two-sided tensor ideal}. Tensoring with $F$ acts as the identity on $\JC$; more accurately every complex in $\KC(\JC)$ is a $\firstMap$-eigenobject.  On the other hand, $\Cone(\secondMap)$ is homotopic to a complex in $\KC(\JC)$, so that tensoring with it takes any complex into this tensor ideal. This behavior is summarized by the following diagram:
\[
\begin{tikzpicture}[baseline=1em]
\node(a) at (0,0) {$\KC^b(A\text{-mod})$};
\node(b) at (4,0) {$\KC^b(\JC) $};
\node(c) at (8,0) {$0$};
\tikzstyle{every node}=[font=\small]
\path[->,>=stealth',shorten >=1pt,auto,node distance=1.8cm]
(a) edge node[above] {$\Cone(\secondMap)\otimes(-)$} (b)
(b) edge node[above] {$ \Cone(\firstMap)\otimes(-)$} (c);
\end{tikzpicture}
\]
Rephrasing slightly, one might say that $\Cone(\secondMap)$ annihilates the monoidal identity, modulo $\JC$, and $\Cone(\firstMap)$ annihilates $\JC$ modulo zero. Heuristically, $F$ acts
by an ``upper-triangular matrix'' with respect to the filtration $0\subset \KC^b(\JC)\subset \KC^b(A\text{-mod})$. The eigenmaps $\firstMap$, $\secondMap$ represent the relationship
between this ``matrix'' and its ``diagonal'' entries, corresponding to the eigenvalues.

To prove that $\{\firstMap, \secondMap\}$ is a saturated set of eigenmaps, it is easier to first prove that their cones annihilate the different layers of the filtration. In the sequel
we develop this trick as a tool to prove that, given a functor which acts in upper-triangular fashion on a filtered additive category, and given maps which are eigenmaps on the
filtration subquotients, then the functor is prediagonalizable.

\section{Subtleties of categorical diagonalization}
\label{sec:subtle}

Despite the many comparisons drawn between categorical diagonalization and its counterpart in linear algebra, there are many situations which arise in the categorical situation which have no analog in linear algebra.  Most notably, categorical diagonalization involves a choice of structure, and there are often many inequivalent choices one could make.  To illustrate this further, let us present a more ``structure-oriented'' perspective on diagonalization in linear algebra.

Let $A$ be a $\KM$-algebra and $f\in A$ an element.  Choose a finite set $I$, scalars $\k_i\in \KM$ and elements $p_i\in A$ ($i\in I$).  We say that $\{(p_i,\k_i)\}_{i\in I}$ is a \emph{diagonalization of $f$} if the $p_i$ are a complete collection of mutually orthogonal idempotents such that $(f-\k_i)p_i = 0 = p_i(f-\k_i)$ and none of the $p_i$ is zero.

Note that we allow repetitions among the scalars, so that the idempotent decomposition $1=\sum_i p_i$ is potentially finer than the decomposition into eigenspaces.  Additionally, the eigenvalue $\k_i$ is determined by the projection $p_i$, so the ``choice'' of $\k_i$ is not much of a choice at all.  Furthermore, for a given $\k\in \KM$, the projection onto the $\k$-eigenspace is simply the sum of $p_i$ over all indices such that $\k_i=\k$.

In the categorified situation, a diagonalization of $F\in \AC$ in a triangulated monoidal category is a choice $\{(\PB_i,\a_i)\}_{i\in I}$ where $I$ is now a finite poset, the $\PB_i$ are nonzero mutually orthogonal idempotents such that $\one\simeq\Tot\{\PB_i,d_{ij}\}$ and $\Cone(\a_i)\otimes \PB_i\simeq 0 \simeq \PB_i\otimes \Cone(\a_i)$.  In analogy with the decategorified situation, $\PB_i$ may project onto a proper subcategory of the $\a_i$-eigencategory.

One can now ask: ``is there some combination of idempotents which projects onto the $\a_i$-eigencategory?''  The answer to this is certainly no, as we shall see.  Morally speaking, two eigenmaps with the same eigenvalue may interact in interesting ways which are not visible at the decategorified level.

Before continuing, recall the notion of a \emph{tight diagonalization of $F$}, $\{(\PB_i,\a_i)\}_{i\in I}$.  In any tight diagonalization, $\PB_i$ annihilates every $\a_j$-eigenobject for $j\neq i$.  Thus, there are no nonzero objects which are both an $\a_i$-eigenobject and an $\a_j$-eigenobject for $i\neq j$.  In some sense, tight diagonalizations satisfy many properties which are intuitive from the perspective of linear algebra.  Non-tight examples abound, however, and in this situation there is a lot of misleading intuition from linear algebra that we intend to nip off at the bud.

Foremost amongst the potential confusions is this. In linear algebra, we can place a partial order on the diagonalizations of $f$ using
the notion of refinement of an idempotent decomposition. There is a unique maximal (coursest) diagonalization, where each $p_i$ is projection to an entire eigenspace. In the categorified
situation, there is also a partial order on diagonalizations of $F$, using refinement (see \S\ref{subsec-idempotentreassoc}). Any tight diagonalization is maximal in this partial order.
However, there is no unique maximal diagonalization in general, so there is no ``best" way to diagonalize a functor.

\begin{remark} Even in linear algebra, non-maximal diagonalizations are useful. When one simultaneously diagonalizes two commuting operators $f$ and $g$, the resulting diagonalization of $f$ is typically non-maximal. \end{remark}




The first two examples below illustrate how eigencategories for different eigenmaps may overlap in interesting ways.  These examples are not obviously diagonalizable, so we follow them with some similar looking examples which are diagonalizable, but somewhat contrived.  There are less contrived examples coming from categorical Hecke theory, which are discussed in \S\ref{subsec:Hecke} and \S\ref{subsec:HeckeConj}.

\subsection{Example: Projective space}
\label{subsec:cohp1}

We learned of this example from Negu\cb{t}, and it is discussed at length in \cite{GNR}. Though it appeared chronologically after the initial drafts of this
paper, we have included it here because it gives an excellent illustration of some of the major subtleties of categorical eigentheory.

Let $\AC$ denote the derived category of coherent sheaves on $\PM^r$, with homogeneous coordinates $[x_0:x_1:\ldots:x_r]$. Let $\OC_p$ denote the structure sheaf of the point $p \in \PM^r$. The sheaves $\OC_p$ are nice ``test objects'' for the application of functors. Let $\LC$ be any line bundle. It is a standard fact that $\LC \ot \OC_p \cong \OC_p$, so that $\OC_p$ is a weak eigenobject for $\LC$ with eigenvalue $\one = \OC$.

Let $F = \OC(1)$, and consider the scalar object $\one=\OC$.  The maps $\one\rightarrow F$ (i.e.~sections  of $\OC(1)$) are given by linear polynomials in $\C[x_0,\ldots,x_r]$. If $f$ is a linear polynomial, we let $\a_f$ denote $f$, thought of as a section. In particular we have $\a_i:\one\rightarrow F$ corresponding to the monomials $x_i$.  Observe that $\OC_p$ is an $\a_f$-eigenobject precisely when $p$ does not lie on the hyperplane $f=0$.   In particular, for any two maps $\a_f$,  $\a_{g}$ there exists an object $\OC_p$ which is simultaneously an $\a_f$ and an $\a_g$ eigenobject.  Consequently, no two such eigenmaps can appear in a tight diagonalization.

\begin{remark} It should be noted that every nonzero element of $\Hom(\1,F)$ is an eigenmap, and almost every element is an eigenmap for $\OC_p$ for any given point $p$. On the other hand, every eigenmap is uniquely determined (up to scalar) by its eigencategory. \end{remark}

Nonetheless, $F$ is pre-diagonalizable.  If $f_0,\ldots,f_r$ is any basis of linear polynomials, then $\bigotimes_{i=0}^r \Cone(\a_{f_i})$ is essentially the Koszul complex which resolves the trivial representation of $\C[f_0,f_1,\ldots,f_r]$, where all the $f_i$ vanish.  It is supported on the subset of $\PM^r$ where $f_i = 0$ for all $i$. That is to say, it is supported on the empty set, so it is zero as desired. Hence the ways in which $F$ can be prediagonalized are in bijection with the ways to choose a basis for linear polynomials.

Note that we can not apply Theorem \ref{thm:introEigenthm}, because $\OC$ is not small (even in the derived category of quasi-coherent sheaves, a countable direct sum $\bigoplus \OC$ is not isomorphic to the countable direct product). 

In \cite{GNR} the situation is rigidified by placing a torus action on $\PM^r$ and considering equivariant sheaves. The structure sheaves $\OC_p$ only admit an equivariant structure if
$p$ is a fixed point for the torus action, which means there is an $i$ such that $p=\bigcap_{j\neq i}\{x_j=0\}$. Now, this gives a reason to prefer the eigenmaps $\a_i$, since no two of
these maps have a common nonzero eigenobject. The eigenmaps $\a_i$ map from $\l_i$ to $F$, where $\l_i$ is $\OC$ with a twisted equivariant structure. If the $\l_i$ are all distinct,
then it may be possible to possible to modify Theorem \ref{thm:introEigenthm} to apply, giving a tight diagonalization of $F$. In any case \cite{GNR} imitates our construction of the
projectors $\PB_i$, proves that they are categorical idempotents, and relates them to constructions in algebraic geometry.

\subsection{Example: The integers}
\label{subsec:integers}

Let $\AC$ denote the category of all abelian groups. Then $\AC$ is additive and symmetric monoidal, with tensor product $\otimes_\Z$ and monoidal identity $\one=\Z$. Then clearly $F=\Z$
is diagonalizable with eigenmap $\Id:\Z \rightarrow \Z$ and eigenprojection $\PB = F = \one = \Z$. However, there are many other eigenmaps for $\Z$. Indeed, for each nonzero $k\in \Z$
multiplication by $k$ gives an eigenmap $\a_k:\Z\rightarrow \Z$. The $\a_k$-eigencategory consists of all $\Z$-modules on which $k$ acts invertibly.

\begin{remark} To make $F$ categorically prediagonalizable, one should restrict to an appropriate full subcategory of $\AC$.\end{remark}

This example shares many of the features of the example of coherent sheaves on $\PM^r$, as one may expect from algebraic-geometric considerations. For example, $\Q$ is an
$\a_k$-eigenobject for each nonzero $k$. There are additional phenomena, such as the fact that the $\a_6$-eigencategory is the intersection of the $\a_2$-eigencategory and the
$\a_3$-eigencategory.

\begin{remark} A completely thorough examination of the theory of categorical diagonalization is likely to contain large portions of algebraic geometry as a subtheory. Algebraic geometry
is interesting, so this is motivation to study categorical diagonalization further. Algebraic geometry is also difficult, so this is also motivation to restrict one's attention to more tractable scenarios. \end{remark}

\subsection{Diagonalizations of $\one$}

Suppose $\AC$ is a monoidal homotopy category and $\one\simeq \Tot\{\PB_i,d_{ij}\}_{i\in I}$ is an idempotent decomposition of identity. Set $F = \one$ and $\a_i:=\Id_{\one}$ for all $i$.  Then $\{(\PB_i,\a_i)\}_{i\in I}$ is a non-tight diagonalization of $\one$.  The next example concerns a special case of this, in which there are many possibilities for the eigenmaps.

\subsection{Example: Semisimple categories}
\label{subsec:semisimple}

This next example shares some superficial similarities with the coherent sheaves example above, but has a ``finite spectrum" and consequently has drastically different behavior.

Let $\BC$ be an additive monoidal category in which $\one$ decomposes as a direct sum $\one\cong e\oplus f$, where $e\otimes f\cong f\otimes e$ are orthogonal idempotents.   An example would be  the category of endofunctors of a semisimple category with two simple objects $L_1$ and $L_2$, where $e = \Hom(L_1,-) \boxtimes L_1$ and $f = \Hom(L_2,-) \boxtimes L_2$.

Obviously the identity map $\Id_\one:\one\rightarrow \one$ is an eigenmap for $\one$, and $\{\one,\Id_{\one}\}$ is a diagonalization of $F:=\one$ (with one projector and one eigenmap). 

However, this is not the only diagonalization of $F$. Let $\a_e,\a_f\in \Hom(\one,F)=\End(\one)$ denote the projections onto $e$ and $f$, respectively.  Then $\a_e$ is an eigenmap with eigencategory consisting of all objects such that $e\otimes M\cong M$, and similarly for $\a_f$. Moreover, $\Big\{\{e,f\},\{\a_e,\a_f\}\Big\}$ is a tight diagonalization of $F$, with indexing set $I=\{e,f\}$ and the trivial partial order.

In fact any nonzero linear combination $c_e \a_e+c_f \a_f$ is an eigenmap, with eigencategory equal to  $\im e$ if $c_f=0$, $\im f$ if $c_e=0$, and all of $\AC$ if both $c_f,c_e$ are nonzero.

More generally, suppose that $\one$ decomposes as a direct sum $\one \cong \bigoplus_{i=0}^r e_i$ of orthogonal idempotent functors. Let $F = \one$ and
$\a_i$ be the projection to $e_i$. Then any nonzero linear combination $\sum c_i \a_i$ is an eigenmap, and its eigencategory is the image of $\bigoplus_{j} e_j$ over the subset for which $c_j \ne 0$.

Let $\a = \sum c_i \a_i$ with all $c_i \ne 0$. Then the collection $\{e_i,\a\}$ is a diagonalization of $F$, where every eigenmap is equal to $\a$. It is not tight.

This example illustrates that a given object can be diagonalized in a variety of ways, with potentially different numbers of eigenmaps and projections.  Generically, one needs only one eigenmap to diagonalize $F$, but there are special eigenmaps which resolve the big eigencategory $\AC$ into smaller subcategories.

\subsection{Fission and fusion}
\label{subsec:fusion}

An important phenomenon appeared in the previous example. Suppose that two eigenmaps $\a_i \co \l \to F$ (for $i=1, 2$) have the same scalar functor $\l$ as their source; we call this a
case of \emph{coincident eigenvalues}. We do not use the term ``repeated eigenvalues,'' because to do so would suggest (incorrectly) the failure of diagonalizability, and also because the $\a_i$-eigencategories may be quite distinct.

We say that that $\a_1$ and $\a_2$ can be \emph{fused} if there exists an eigenmap $\b:\l\rightarrow F$ whose eigencategory contains both the eigencategory of $\a_1$ and that of $\a_2$.  We also say that the $\b$-eigencategory fuses the $\a_1$- and $\a_2$- eigencategories.

\begin{remark} In the (non-equivariant) coherent sheaves example of \S\ref{subsec:cohp1}, all eigenmaps have the same eigenvalue $\OC$ but distinct eigencategories.  Two different eigenmaps cannot be fused because no eigencategory properly contains any other. \end{remark}

\begin{remark}
In \S\ref{subsec:fusionGuarantee} we discuss a common scenario in which every pair of coincident eigenmaps can be fused.
\end{remark}

Fusion of eigenmaps can be a useful tool in applying the diagonalization theorem.  For instance, suppose there exist invertible scalars $\l_i$ and a saturated collection of eigenmaps $\{\a_i:\l_i\rightarrow F\}$.  It may occur that some eigenvalues coincide, in which case Theorem \ref{thm:introEigenthm} is not applicable.  However, if $\b$ fuses together two eigenmaps $\a_i$, $\a_j$, then both factors $\Cone(\a_i)$, $\Cone(\a_j)$ may be replaced by a single factor $\Cone(\b)$.  If all coincident eigenmaps can be fused in this way, then one is a good position to apply the Diagonalization Theorem.  It would then remain to \emph{fission} the resulting diagonalization into a diagonalization which distinguishes the fused eigenmaps.

\begin{remark}
This technique will play an important role in a sequel, in which we diagonalize the full twist in the Hecke category (see \S \ref{subsec:HeckeConj}).  
\end{remark}

\begin{remark}
The potential of fission is another major difference between linear algebra and its categorification. In linear algebra two vectors in the same eigenspace for $f$ cannot be distinguished using $f$ alone, because any polynomial expression $p(f)$ acts as a constant on each $f$-eigenspace. The categorical analogue of this is false, as we have seen.  Even if $M_1,M_2$ are both $\b$-eigenobjects for some $\b:\mu\rightarrow F$, it may be possible to find maps $\a_i:\mu\rightarrow F$ such that $M_1$ is an $\a_1$-eigenobject but not a $\a_2$-eigenobject, and vice versa. \end{remark}



Let us give a more precise illustration of what we mean by fission and fusion in the context of diagonalization.  Suppose $F\in \AC$ is an object of a monoidal homotopy category, and suppose $\{(\PB_i,\a_i)\}_{i\in I}$ is a diagonalization of $F$.  Suppose $i_1$ and $i_2$ are incomparable indices with the same eigenvalue $\mu$, and suppose there is a map $\b:\mu\rightarrow F$ such that the $\b$ eigencategory contains both the $\a_{i_1}$- and $\a_{i_2}$-eigencategories.  We define a coarser diagonalization of $F$ in two steps.

First, define a new family of eigenmaps $\a_i':\l_i\rightarrow F$ by $\a_i'=\a_i$ if $i\not\in\{i_1,i_2\}$ and $\a_i=\b$ otherwise.  Then $\{(\PB_i,\a_i')\}_{i\in I}$ is a diagonalization of $F$ which is ``less tight'' than the original diagonalization.  

Second, we merge the idempotents $\PB_{i_1}$ and $\PB_{i_2}$ together.  Define a new poset $J$ in which $i_1$ and $i_2$ have been fused.  Formally, let $J:=I/\sim$ where $\sim$ is the equivalence relation generated by $i_1 \sim i_2$.  Since $i_1$ and $i_2$ are incomparable, $J$ inherits a partial order from $I$. We define a new, $J$-indexed diagonalization of $F$ by
\begin{enumerate}
\item $(\QB_j,\b_j) = (\PB_j,\a_j)$ if $j\not\in\{i_1,i_2\}$, and
\item $ (\QB_j,\b_j) = (\PB_{i_1}\oplus \PB_{i_2},\b)$ if $j=\{i_1,i_2\}$ is the non-singleton equivalence class.
\end{enumerate}

It is an exercise to check that $\{(\QB_j,\b_j)\}_{j\in J}$ is a diagonalization of $F$.  We say that $\{(\QB_j,\b_j)\}_{j\in J}$ is obtained from $\{(\PB_i,\a_i\}_{i\in I}$ by \emph{fusing} the $i_1$ and $i_2$ terms.

For example in \S\ref{subsec:semisimple} we essentially showed that the fusion of eigenmaps converts the diagonalization $\{(e,\a_e),(f,\a_f)\}$ into a coarser diagonalization $\{\one, \b\}$, where almost every linear combination $c_e \a_e + c_f \a_f$ is a valid choice for $\b$.


In the next section we discuss a situation in which coincident eigenmaps can always be fused.

\subsection{Guaranteed fusion}
\label{subsec:fusionGuarantee}

Suppose that $\VC$ is an additive Krull-Schmidt category over a field $\Bbbk$, and $X$ is an indecomposable object. Then $\End(X)$ is local with maximal ideal $\mg$ consisting of all the non-isomorphisms. If $F X \cong \l X$ and $\l$ is invertible then, by fixing any identification $F X \to \l X$, we may identify $\Hom(\l X, F X)$ with $\End(\l X)$ and thereby with $\End(X)$.

 For an object $X \in \VC$ with $F X \cong \l X$, let $U_X \subset \Hom(\l,F)$ denote the subset of eigenmaps for $X$, i.e.~the set of maps $\a$ such that $X$ is an $\a$-eigenobject. The following lemma is obvious.

\begin{lemma} \label{lem:decomposeU}  $U_{X \oplus Y} = U_X \cap U_Y$. \end{lemma}

The set $U_X$ is precisely the complement of the preimage of zero under the composition \begin{equation} \label{eq:iseigen}
\Hom(\l,F)
\to \Hom(\l X, F X) \simto \End(X) \to \End(X)/\mg.
\end{equation} The first map is tensoring with $\Id_X$, the second is our chosen identification, and the third is the quotient map.

Let us equip the morphism spaces in $\VC$ (and in $\End(\VC)$) with a topology.

\begin{defn} Let $V$ be a finite dimensional vector space over $\Bbbk$. The \emph{linear topology} on $V$ (and on the projective space $\PM V$) is the topology where the closed sets are
finite unions of linear subspaces. \end{defn}

In this topology, every nonempty open set is dense, because every proper closed set has strictly smaller ``dimension.'' This is the coarsest topology on the collection of all vector
spaces where linear maps are continuous and where $\{0\}$ is a closed set.

Now, the maps of \eqref{eq:iseigen} are continuous, and consequently $U_X$ is an open set, which is dense whenever it is nonempty. Using Lemma \ref{lem:decomposeU}, the set $U_X$ for any
finite direct sum of indecomposable objects is also a dense open set.

Consequently, if the eigencategory of $\a_1$ has finitely many indecomposable objects, then any eigenmap for their direct sum is an eigenmap for the entire eigencategory, and the subset
$U_1$ of $\Hom(\l,F)$ consisting of eigenmaps for the entire eigencategory of $\a_1$ is also a dense open set. Suppose $\a_2$ is another eigenmap with eigenvalue $\l$. As the
intersection $U_1 \cap U_2$ of two dense open sets is another dense open set, one can deduce that almost every linear combination of $\a_1$ and $\a_2$ yields a fused eigenmap. Of course,
one can fuse multiple coincident eigenmaps in this way.

\begin{remark} This argument does not apply when there are infinitely many indecomposable objects, as the infinite intersection of open sets need not be open. See the example of coherent sheaves on $\PM^r$ in \S\ref{subsec:cohp1}, where coincident eigenmaps can not be fused.   \end{remark}

We are more interested in actions on triangulated categories than on additive categories.  In such categories there are typically infinitely many indecomposables, so a direct application of the above argument breaks down, unfortunately.  However, if one has a particularly nice, finite set of objects which generate (in the sense of triangulated categories) then there is a modification of this argument which goes through intact.  We will return to this in a sequel.

This discussion suggests that in the Krull-Schmidt case, coincident eigenmaps can be fused when the category being acted upon satisfies some finiteness condition, but not in general. The converse question is more difficult: given an eigenmap $\b$, when do there exist coincident eigenmaps $\a_1$ and $\a_2$ which fission its eigencategory? We know of no reasonable criterion as of yet.

\subsection{What is the spectrum of an endofunctor?}
\label{subsec:spectrum}

The above examples should make it clear that this is a complicated question. Should one keep track of the scalar functors, or the eigenmaps, or the possible saturated sets of eigenmaps,
or the possible eigencategories? If two eigencategories can undergo fusion, how should this relationship be encoded? We do not know a satisfactory answer. Let us make some basic definitions which will help future discussion, at least.

\begin{defn} Let $F$ be a functor acting on a triangulated category $\VC$. The \emph{spectrum of forward eigenmaps}, or simply \emph{forward spectrum} $\Spec(\l,F)$ (or $\Spec_{\to}(\l,F)$) is the set of all eigenmaps $\a:\l\rightarrow F$. We place a
topology on $\Spec(\l,F)$ where a base for the open sets is given by $U_X$ for $X \in \VC$.  Define $\Spec(F)$ to be the union $\bigsqcup_{\l} \Spec(\l,F)$ where $\l$ ranges over a set of representatives of isomorphism classes of scalar objects.  Two eigenmaps in $\Spec(F)$ are \emph{equivalent} if they have the same
eigencategory. Taking the quotient of $\Spec(F)$ by this equivalence relation, we obtain the topological space $\XSpec(F)$.  One may also consider the similar notion of the \emph{backward spectrum} $\Spec_{\leftarrow}(F)$.\end{defn}

Why the set $U_X$ should be open has been discussed in the previous section. Note that the quotient map $\Spec_{\to} \to \XSpec_{\to}$ factors through $\PM \Spec_{\to}$, the projectivization of the forward spectrum.

\begin{example} Consider the general example from \S\ref{subsec:semisimple}, where $\one = \bigoplus_{i=0}^r e_i$, and the eigenmaps are the sums $\sum c_i \a_i$ where $c_i \ne 0$ for some $i$. The set $\XSpec_{\to}(\1,\1)$ has $2^{r+1}-1$ elements, each eigenmap being identified with the nonempty subset of $\{0, \ldots, r\}$ recording which coefficients $c_i$ are nonzero. The topology is the subset topology. \end{example}

\begin{remark} Note that these definitions cater only to the notion of forward or backward eigenmaps, and ignore more generalized notions of eigenmaps, such as those mentioned in
\S\ref{sec:generalizations}. \end{remark}

\begin{remark}
Many important properties of eigenmaps (for instance whether two eigenmaps can be fused) can be phrased entirely in terms of the corresponding eigencategories, which is one reason to prefer $\XSpec(F)$ over $\Spec(F)$.
\end{remark}



Finally, let us remark on an additional structure on $\XSpec(F)$ which we believe plays an important role.  In order to describe it, we will utilize the language developed in \S \ref{sec-homotopy} and \S \ref{sec-decompositions}.  
 For $\PB\in \AC$ a counital idempotent, let us say $\a\in \XSpec(F)$ is \emph{$\PB$-compatible} if $\PB$ fixes each $\a$-eigenobject up to isomorphism, and $\a$ is \emph{$\PB$-orthogonal} if $\PB$ annihilates each $\a$-eigenobject.  We say that $\PB$ \emph{separates} $\a$ from $\b$ if $\a$ is $\PB$ compatible and $\b$ is $\PB$-orthogonal.  Finally, we say $\a$ is separated from $\b$ if there exists a $\PB$ which separates $\a$ from $\b$.

The relation of separation is intimately connected with the partial order on eigenmaps which appear in a diagonalization.  For instance if $\a_0,\a_1,\ldots,\a_r$ are eigenmaps to which Theorem \ref{thm:introEigenthm} applies, then $\a_i$ is separated from $\a_j$ if $i>j$.

In linear algebra one can deduce the spectrum of a linear operator (even a non-diagonalizable one) by examining the characteristic polynomial. However, we do not know of any analogue of characteristic polynomial in the categorified setting, and the problem of computing the spectra of endofunctors seems quite difficult in general. In this paper and related work we are content with finding examples of eigenmaps (enough to prove categorical diagonalizability), and we ignore the problem of computing the entire forward spectrum.

\section{Applications of categorical diagonalization}
\label{sec:introapps}

In recent decades mathematicians have unearthed a rich new world of categorical representation theory. Diagonalization is a powerful tool in representation theory, and we intend that
categorical diagonalization should become an important tool in categorical representation theory.  The applications we currently have in mind are for categorified quantum groups and Hecke algebras. We imagine the average reader is more familiar with the quantum group story, so we begin there.

\subsection{Structural and constructive results}
\label{subsec:CRT1}
Traditional representation theory is filled with classifications and constructions. For example, consider the complex semisimple lie algebra $\sl_2$, and its category of representations. Here are some of the key results: \begin{itemize}
\item Representations are semisimple; they split canonically into isotypic components.
\item Each isotypic component $V$ of type $\l$ (a dominant integral weight) can be described canonically as the external tensor product $V = V_\l \boxtimes V(\l)$ of a minimal (i.e. simple) representation $V_\l$ and a multiplicity vector space $V(\l)$ (the highest weight space).
\item Each minimal representation $V_\l$ has a universal description (in this case, as a simple quotient of a Verma module).
\item Each minimal representation $V_\l$ also has an explicit construction (using direct formulas, or homogeneous polynomials and differential operators, etcetera), allowing one to compute with them.
\end{itemize}
Combining these points, one has a relatively complete understanding of every $\sl_2$-representation.

In their seminal paper \cite{ChuRou}, Chuang and Rouquier defined categorical representations of $\sl_2$, and proved ``lifts" of these classification results. \begin{itemize} \item
Categorical representations have filtrations, whose subquotients are categorifications of isotypic components. \item Each isotypic component $\VC$ of type $\l$ can be described
canonically as an exterior tensor product $\VC = \VC_\l \boxtimes \VC(\l)$ of a minimal categorical representation $\VC_\l$ and a multiplicity category $\VC(\l)$. \item Each minimal
categorical representation $\VC_\l$ has a universal description. \item The minimal categorical representation $\VC_\l$ has an explicit construction. \end{itemize} Later, as categorical
representations of other complex semisimple lie algebras were being developed by Rouquier \cite{Rouq2KM-pp} and Khovanov-Lauda \cite{KhoLau09, KhoLau10, KhoLau11}, minimal categorical
representations $\VC_\l$ were described using the so-called ``cyclotomic quotient," effectively allowing one to compute (see e.g. \cite{LauVaz}) how the categorification $\UC(\sl_2)$ of
the enveloping algebra $U(\sl_2)$ acts on these minimal categorical representations.

Unlike the traditional case, categorical representations are not semisimple, so reducible representations can have interesting new structure. For example, Webster \cite{WebsKIHRT}
constructed natural categorifications of tensor product representations, and Webster and Losev \cite{LosWeb} described the structural properties of minimal categorifications of tensor
products: each weight category has the structure of a standardly stratified category, or in the simpler case of a tensor product of minuscule representations, a highest weight category.
The filtration by isotypic components descends to a highest weight filtration on each weight category.

In an ongoing series of delightful papers \cite{MazMie11, MazMie14, MazMie16, MazMie16-2}, Mazorchuk and Miemietz have attempted to place structural results like those of Chuang and
Rouquier in a natural and general setting. They describe the representation theory of \emph{fiat $2$-categories}, which is a natural categorical lift of the idea of an algebra with a
distinguished antiinvolution. One key observation which underlies many of their results is that one is always categorifying, not just an algebra or representation, but an algebra or
representation with a fixed basis, given by the symbols of the indecomposable (projective) objects. As such, one should study this algebra and its categorification using the theory of
\emph{cells}. One might hope for a Chuang-Rouquier-style result like: every categorical representation has a filtration by isotypic copies of cell modules, and isotypic copies of cell
modules are obtained from a universally-described cell module by inflation against a multiplicity category. This result fails in general, although it does hold in important special
cases, like the Hecke category in type $A$ in characteristic zero, see \cite{MazMie16}.


\subsection{The spectral approach}
\label{subsec:CRT2}

There is also a spectral approach to the representation theory of $\sl_2$. The universal enveloping algebra $U(\sl_2)$ has an element $c$, the \emph{Casimir element}, which generates
the center $Z(U(\sl_2))$. The decomposition of a representation into its isotypic components is the same as its decomposition into eigenspaces for $c$. On the categorical side,
Beliakova-Khovanov-Lauda \cite{BKL} have constructed a complex $C$ in the homotopy category of $\UC(\sl_2)$, which lives in its Drinfeld center, and descends to $c$ in the Grothendieck
group.

Thus the natural question is whether one can use the techniques of categorical diagonalization, applied to the complex $C$, to construct categorical projections to isotypic components.
If so, the idempotent resolution of identity equips any categorical representation with a canonical filtration whose subquotients are isotypic components. One expects this will give a
new perspective on the filtration constructed by Chuang and Rouquier.

\begin{remark} Chuang and Rouquier work with abelian categorifications, while we work with homotopy categories (of projective objects). One advantage of our perspective would be that it
makes homological algebra easier, because it works directly with projective resolutions rather than simple objects. \end{remark}

\begin{remark} The interactions between categorical diagonalization and highest weight structures, such as the ones discussed in the previous section, are completely unexplored and
presumably quite interesting. \end{remark}

\begin{remark} Optimistically, one should also be able to reproduce the constructions of minimal categorical representations using projections to eigencategories, though the details are
mysterious. Let us contrast this hypothetical construction to the usual cyclotomic quotient construction. The latter imposes certain cyclotomic relations upon the Khovanov-Lauda-Rouquier
algebra. The former would instead identify certain infinite complexes over $\UC(\sl_2)$, for which the cyclotomic relations would act as nulhomotopic endomorphisms. Instead of formally
imposing these relations, one finds a natural setting where they already act trivially. A careful examination of the endomorphism rings (modulo homotopy) of the projection operators and
their images under the action of $\UC(\sl_2)$ would hopefully permit one to recover these complexes and the cyclotomic Khovanov-Lauda-Rouquier algebra in a natural way. \end{remark}

The (weak) categorical eigenvalues of $C$ are not invertible, and $C$ does not admit any forward eigenmaps. Consequently, none of the theorems in this paper can be applied to diagonalize
$C$. In \S\ref{sec:generalizations} we discuss generalizations of eigenmaps, focusing on the mixed eigenmaps introduced in \S\ref{subsec:mixed}. In \S\ref{subsec:casimir} we conjecture that $C$ has an interesting mixed eigenmap in a cyclotomic quotient of $\UC(\sl_2)$. Alas, there do not appear to be enough mixed eigenmaps to prediagonalize $C$ outside of special cases. Moreover, we do not know how to construct enough the complexes $\Cij{j}{i}$ or the projectors $\idemp_i$ for non-invertible eigenvalues (the appropriate analog of Koszul
duality is not known).

We hope to develop our technology further to address this natural application.

Let us now segue to a field where we have interesting positive results.

\subsection{Categorical Hecke theory}
\label{subsec:Hecke}

As mentioned above, Mazorchuk and Miemietz have already begun to study the structure results in categorical Hecke theory from a cellular perspective. In this chapter we discuss the spectral approach. A
thorough introduction to this topic and its applications will be found in the sequel to this paper, so we only give a brief version here.

To explain, consider the example of the symmetric group $S_n$. Traditionally, one constructs the Specht modules of $S_n$ by explicitly constructing a family of \emph{Young idempotents}
$p_T\in \C[S_n]$, indexed by the set of standard tableaux on $n$ boxes, such that the module $\C[S_n]p_T$ is irreducible and depends only on the shape $\Sh(T)$ of $T$ up to
isomorphism. For this reason, many references only construct the idempotents $p_T$ for a special tableau of each shape. However, in the decomposition of the regular representation of
$\C[S_n]$, all tableaux are required. There is also a central idempotent $p_{\l}$ for each partition $\l$ of $n$, which projects to the entire isotypic component; one has $p_\l = \sum_{\Sh(T) = \l} p_T$.

The idempotents $p_T$ can be realized in terms of diagonalization. The Young-Jucys-Murphy operators $0=y_1,y_2,\ldots,y_n$ generate a maximal commutative subalgebra
$Y \subset \C[S_n]$. The action of $Y$ on $\C[S_n]$ by right multiplication is simultaneously diagonalizable, with joint eigenspaces indexed by standard Young tableau $T$.
Decomposing $1\in\C[S_n]$ into simultaneous eigenvectors gives $1= \sum_T p_T$, where $p_T$ are the Young idempotents. Then the decomposition of $\C[S_n]$ into its simultaneous
eigenspaces $\C[S_n]\cong\bigoplus_T \C[S_n]p_T$ gives a decomposition of $\C[S_n]$ into irreducible summands. Thus, simultaneous diagonalization of the Jucys-Murphy operators is key to
many of the classical and more modern approaches to the representation theory of symmetric groups, such as \cite{OkounkovVershik}.

There is a similar story for Hecke algebras in type $A$. The Hecke algebra $\Hecke_n$ of type $A_{n-1}$ is a $\Z[v,v\inv]$-deformation of $\Z[S_n]$. There are Jucys-Murphy operators $1
=j_1, j_2, \ldots, j_n$ in this setting as well\footnote{The operators $j_i$ all come from pure braids, so after specializing to $v=1$ they all yield the identity element of $\Z[S_n]$.
Instead, $j_i$ deform $y_i$ in a more subtle way.}, whose simultaneous eigenvalues are also described by tableaux. There is also a connection between the modules $\Hecke_n p_T$ and the
cell modules constructed by Kazhdan-Lusztig \cite{KL79}, but we refer to the sequel for more details.

We prefer different diagonalizable operators which generate the same commutative subalgebra. There is a map $\rho \co \Br_n\rightarrow \Hecke_n^\times$ from the braid group on $n$
strands. Let $\ft_k$ denote full-twist braid on the first $k$ strands. Then $1=\ft_1,\ft_2,\ldots,\ft_n$ generate a maximal commutative subgroup of the braid group, whose images in
$\Hecke_n$ generate a maximal commutative subalgebra, which agrees with the algebra generated by $j_i$\footnote{Indeed, $j_k = \ft_k \ft_{k-1}^{-1}$, and $j_1 = \ft_1 =1$.}. One advantage
of full twists over Jucys-Murphy operators is that the complete full twist $\ft_n$ is central in $\Hecke_n$.

The Hecke algebra $\Hecke(W)$ (for any Coxeter group) is categorified by the monoidal category of Soergel bimodules, or by its homotopy category. The braid group representation
$\rho:\Br(W) \rightarrow \Hecke(W)$ was categorified by Rouquier, who produced a complex of Soergel bimodules for each braid, well-defined up to homotopy. If $W$ is finite, the longest
element lifts to an element of $\Br(W)$ called the half-twist, and its square is called the full-twist $\FT(W)$. We write $\FT_n$ for $\FT(S_n)$.

\subsection{Conjectures and results}
\label{subsec:HeckeConj}

We conjecture that the Rouquier complex of the full-twist braid is categorically diagonalizable, for any finite Coxeter group, and that for a chain of finite Coxeter groups, the full
twists of each are simultaneously diagonalizable. We prove this conjecture for type $A$ in a sequel, and we prove it for dihedral types in another sequel. The example of type $A_1$ is in
\S\ref{subsec:typeA1}.

Let us split these conjectures and results into two halves, one dealing with prediagonalizability, one with diagonalizability. Because $\l$ typically denotes a partition, rather than a
scalar object, we write $\g$ for a scalar object in this section.

\begin{thm} \label{thm:typeAprediag} For each partition $\l\vdash n$ there is a shift $\g_\l$ and a chain map $\a_\l : \g_\l\one\rightarrow \FT_n$. This equips $\FT_n$ with a
prediagonalizable structure: one has \[\bigotimes_{\l \vdash n} \Cone(\a_\l) \cong 0. \] \end{thm}

\begin{conj}\label{conj:otherprediag} Let $(W,S)$ be a finite Coxeter system. For each 2-sided cell $\l$ of $W$ there is a shift $\g_\l$ and a chain map $\a_\l : \g_\l\one\rightarrow
\FT(W)$. This equips $\FT(W)$ with a prediagonalizable structure. \end{conj}

\begin{thm}\label{thm:dihedralprediag} Conjecture \ref{conj:otherprediag} holds for finite dihedral groups. \end{thm}

\begin{remark} The shifts $\g_\l$ are determined from Lusztig's $\ab$-function for the cell $\l$ and the opposite cell $w_0 \cdot \l$. In type $A$, $\g_\l$ is determined by the total
content (i.e. the sum of the contents of each box) of the partition, and the total column number of the partition. Let us defer the details to the sequel. \end{remark}

\begin{remark} When the Hecke category is defined in finite characteristic, the cell theory of the Hecke algebra changes, and new cells may appear. Amazingly, it seems that new eigenmaps
appear as well! In the sequel we also examine the case of type $B_2$ in characterstic 2, and prove an analog of Theorem \ref{thm:dihedralprediag}. \end{remark}

Now we hope to construct idempotents $\PB_\l$ which project to the $\a_\l$-eigenspaces.

\begin{conj} \label{conj:otherdiag} Let $(W,S)$ be a Coxeter system. There are idempotents $\PB_\l$ for each two-sided cell such that $\{\PB_\l,\a_\l\}$ is a diagonalization of
$\FT(W)$, respecting the usual partial order on two-sided cells. \end{conj}

For dihedral groups, Theorem \ref{thm:introEigenthm} applies as stated, yielding the following theorem.

\begin{thm} \label{thm:dihedraldiag} Conjecture \ref{conj:otherdiag} holds for finite dihedral groups. \end{thm}

In type $A$, starting with $S_6$, there are coincident eigenvalues: for example, $\g_{(3,1,1,1)} \cong \g_{(2,2,2)}$. Consequently, we can not apply Theorem \ref{thm:introEigenthm}
directly. Our approach to bypass this problem combines the fusion of eigenmaps with simultaneous diagonalization. We prove that eigenmaps $\l$ with coincident eigenvalues can be fused,
and that we can apply Theorem \ref{thm:introEigenthm} to the fused prespectrum to construct idempotents $\PB_{[\l]}$ indexed by certain equivalence classes of partitions. Then, by
mutually refining these idempotent decompositions for each full twist $\FT_1$ through $\FT_n$, we can construct projections $\PB_T$ for each standard tableau $T$. Reassociating the
tableaux of a given shape together, we obtain idempotents $\PB_{\l}$ which project to the $\a_\l$ eigenspaces, and thus we successfully fission $\PB_{[\l]}$ into its more interesting
pieces.

\begin{remark} In \S\ref{subsec-relative} we state and prove a generalization of Theorem \ref{thm:introEigenthm} which we term the Relative Diagonalization Theorem.  It is actually this theorem that  allows one to diagonalize $\FT_n$ ``relative to'' the diagonalizations of $\FT_1,\ldots,\FT_{n-1}$, which one constructs inductively.\end{remark}

Consequently, we have the following theorem.

\begin{thm} \label{thm:typeAdiag} There are central idempotents $\PB_\l$ for each $\l \vdash n$ such that $\{\PB_\l,\a_\l\}$ is a diagonalization of $\FT_n$, respecting the dominance
order on partitions. Moreover, there are (non-central) idempotents $\PB_T$ for each standard tableau $T$, which simultaneously diagonalize $\FT_1$ through $\FT_n$, respecting the
dominance order on tableaux. These idempotents categorify the primitive Young idempotents in the Hecke algebra. \end{thm}

In further work, we aim to reproduce many of the arguments and constructions in \cite{OkounkovVershik} on a categorical level.

The theory presented here also implies the existence of quasi-idempotent complexes $\KB_T$.   In addition to its representation-theoretic interest, this theorem has many applications to knot theory, especially to the triply graded knot homology of torus links, and to definitions of colored link homologies. Again, more details on these applications are in the sequel.

\begin{remark} It is clear from this type $A$ example that something more sophisticated must be done to diagonalize full twists for other Coxeter groups. Even the spectral theory of
chains of full twists is not well-studied in the literature. \end{remark}

\begin{remark} For certain infinite Coxeter groups, it appears that there are infinite complexes of Soergel bimodules which behave as though they were the Rouquier complex attached to
the longest element, even though no such longest element exists. These infinite complexes are idempotent, and serve as both the half-twist and the full-twist. There is an analog of
Conjecture \ref{conj:otherdiag} for such complexes, which we prove for the infinite dihedral group in the aforementioned sequel. \end{remark}

\subsection{Example: type $A_1$}
\label{subsec:typeA1}

The examples from \S\ref{sec:extendedintro} involving $A$-modules (starting in \S\ref{subsec:intro_example1}) can be lifted to examples involving Soergel bimodules in type $A_1$. Soergel
bimodules are graded, unlike $A$-modules, but otherwise the examples are almost identical.

Let $R = \QM[x_1, x_2]$ be the polynomial ring with $\deg x_i = 2$. Let $s$ denote the simple reflection of the symmetric group $S_2$, which acts on $R$. Soergel bimodules are a
particular kind of graded $R$-bimodule, from which it inherits its monoidal structure. Let $B_s$ denote the Soergel bimodule $R \ot_{R^s} R(1)$, where $(1)$ denotes the grading shift, so
that the element $1 \ot 1$ lives in degree $-1$. Let $\a_s = x_1 - x_2 \in R$. Then $\a_s^2 \in R^s$, so that $\a_s^2 \ot 1 - 1 \ot \a_s^2 = 0$ inside $B_s$. Note that $B_s$ is also a
ring.

Comparing this to \S\ref{subsec:intro_example1}, one should replace $A$ with $B_s$ and $\Z$ with $R$. Note that $B_s \ot B_s \cong B_s(1) \oplus B_s(-1)$. Note also that $(\a_s \ot 1 + 1 \ot \a_s)(\a_s \ot 1 - 1 \ot \a_s) = 0$ inside $B_s$, analogous to the fact that $(x+1)(x-1)=0$ inside $A$.

The full twist complex $\FT_2$ is homotopy equivalent to the following complex of Soergel bimodules:
\[
\FT_2 \simeq \Big(\begin{diagram}\underline{B_s}(-1) & \rTo & B_s(1) &\rTo & R(2)  \end{diagram}\Big).
\]
We have underlined the term in homological degree zero. The first map is multiplication by $(\a_s\otimes 1 -1\otimes \a_s)$. The second map is multiplication.

There are two eigenmaps here. The map $R(-2)\rightarrow B_s(-1)$ sending $1 \mapsto (\a_s \ot 1 + 1 \ot \a_s)$ lifts to a chain map $\firstMap:R(-2){}\rightarrow \FT_2$. The
inclusion of the final term $R(2)$ induces a chain map $\secondMap:R(2)[-2]\rightarrow \FT_2$. It is a straightforward computation to confirm that $\Cone(\firstMap)\otimes
\Cone(\secondMap) \simeq \Cone(\secondMap) \ot \Cone(\firstMap) \simeq 0$, hence $\FT_2$ is categorically prediagonalizable.

One can also confirm that $B_s$ is a $\firstMap$-eigenobject. Like before, $\Cone(\firstMap)$ is a $\secondMap$-eigenobject, though there are no $\secondMap$-eigenobjects which are
concentrated in a single homological degree. Unlike in \S\ref{subsec:intro_example1}, the two eigenvalues of $\FT_2$ can be distinguished on the Grothendieck group, as they are $v^2$ and
$v^{-2}$ (where $v$ is the decategorification of the grading shift).

We encourage the reader to compute the complexes $\firstProj$ and $\secondProj$, analogous to those in \S\ref{subsec:interpolation_example}. The non-examples of \S\ref{subsec:nonexamples} can also be reproduced in this setting.

%
%
%

\begin{remark} \label{rem:unconvincing} Let $G = (0\rightarrow \underline{B_s}\rightarrow R(1)\rightarrow 0)$. As in \S\ref{subsec:nonexamples}, this complex categorifies a diagonalizable
operator but is not itself diagonalizable. Analogous to Remark \ref{rem:preunconvincing}, there is a map $R_s(-1) \to G$ which seems like an eigenmap for the eigenobject $B_s$. Here,
$R_s$ is the analog of $\sign$; it is the $R$-bimodule which is isomorphic to $R$ as a left module, and for which the right action of $R$ is twisted by the simple reflection $s$. However,
unlike $\sign$, $R_s$ is not in the Drinfeld center of the category of $R$-bimodules (and it is not even in the category of Soergel bimodules), so it is not permissible as a scalar
object. \end{remark}

\begin{remark} \label{rmk:fortopologists}  Eigenmaps appear in topological situations as well.  For instance, let $\TC\LC_{n,m}$ denote Bar-Natan's category of tangles and cobordisms associated to a rectangle with $n$ marked points on the top and $m$ marked points on the bottom of the rectangle \cite{B-N05}.  Associated to any $(n,m)$ tangle diagram $T$ there is a complex $\llbracket T\rrbracket \in \KC^b(\TC\LC_{n,m})$.  Let $X,X\inv\in \KC^b(\TC\LC_{2,2})$  be the complexes associated to the positive and negative crossings, and let $C\in \KC^b(\TC\LC_{2,0})$ be the (planar) the cup diagram, so that $C\otimes X \simeq C\otimes X\inv \simeq C$ up to a shift.  Then there is a map $\a:X\inv \rightarrow X$ such that $\Id_C\otimes \a$ is a homotopy equivalence. Tensoring on the right with $X$ gives an eigenmap $\a':\one\rightarrow X^{\otimes 2}$ for $C$.  Similar examples can be found for other link homology theories. \end{remark}

\subsection{Addendum: connections to coherent sheaves}
\label{subsec:introHilb}

After preliminary versions of this work were made available, Gorsky, Negut, and Rasmussen took our idea of categorical diagonalization and ran with it, beating us to the presses with their beautiful paper \cite{GNR}. We want to avoid making too many comments about \cite{GNR}, because we do not wish to create a dangerous timeloop. Nonetheless, we will say a few words
in this addendum.

In \cite{GNR} it is conjectured that there exists a very deep relationship between Soergel bimodules in type $A$ and the \emph{flag Hilbert scheme} $\FHilb_n(\C^2)$ (we refer to \emph{loc.~cit.} for details).  Precisely, they conjecture the existence of a pair of adjoint functors
\[ \iota_\ast : \KC^b(\SBim_n)\leftrightarrow D^b_{\C^\times\times \C^\times}(\FHilb_n(\C^2)):\iota^\ast, \]
where this latter category is the derived category of sheaves which equivariant with respect to the action of $\C^\times \times \C^\times$. These functors categorify the projection to and inclusion of the Jucys-Murphy subalgebra of the Hecke algebra, and are expected to enjoy a number of remarkable properties. Categorical diagonalization plays a key role in formulating and understanding their conjectures. A discussion of their conjecture belongs more appropriately in the sequel on type $A$.

Also in \cite{GNR}, they make the observation that giving the data of a prediagonalizable functor $F$ in a triangulated monoidal category $\AC$ is the same as giving a monoidal functor $\Coh(\PM^r) \to \AC$. Here, $\Coh(\PM^r)$ is the category of coherent sheaves on $\PM^r$.  The monoidal functor sends $\OC(1)$ to $F$, and the sections $x_i$ of $\OC(1)$ to the eigenmaps of $F$. Additional structures on $\AC$ (e.g. grading shifts, homological shifts) are accounted for by placing a torus action on $\PM^r$ and considering equivariant coherent sheaves.

For example, eigenmaps induce an action of the polynomial ring $\Z[\frac{\a_0}{\a_i}, \ldots, \hat{\frac{\a_i}{\a_i}}, \ldots, \frac{\a_r}{\a_i}]$ on the projector $\PB_i$, and observations made in \cite{GNR} explain this as an action of the coordinate ring of the hyperplane complement $\{x_i\neq 0\}$ in $\PM^r$. Moreover, our diagonalization theorem is interpreted to give a categorification of the Thomason localization formula in equivariant $K$-theory.

\subsection{Connections to previous work on idempotents}
\label{subsec:introTL}

In order to put the present work in perspective, we provide some history on related categorified idempotents.

In \S \ref{subsec:Hecke} we discussed a family of primitive idempotents $p_T\in \Hecke_n$ in the Hecke algebra associated to $S_n$, indexed by standard Young tableaux $T$ on $n$ boxes.  The $\sl_N$ quotient of the Hecke algebra is by definition the quotient of $\Hecke_n$  by the two-sided ideal generated by the $p_T$ whose tableaux have more than $N$ rows. This quotient is the image of $\Hecke_n$ inside the endomorphism ring of the $n$-th tensor power of the standard representation of $\sl_N$, appearing in Schur-Weyl duality. The $\sl_2$ quotient is isomorphic to the Temperley-Lieb algebra.  The one-row idempotent ($q$-Young symmetrizer) $p_{(n)}$ survives in each $\sl_N$ quotient, and we refer to it as the Jones-Wenzl projector \cite{Wenzl87}.

Below we will focus on the $\sl_2$ quotient, since it is the most familiar (the $\sl_2$ quotient was also the first to be considered in the context of categorification in representation theory and topology \cite{BFK,Kh00}).  In 2010 Cooper-Krushal \cite{CK12a}, Rozansky \cite{Roz10a}, and Frenkel-Stroppel-Sussan \cite{FSS} categorified the Jones-Wenzl projectors.   The Cooper-Krushkal and Rozansky projectors both live in the same category $\KC^-(\BN_n)$and are isomorphic.  Here $\BN_n$ is Bar-Natan's category of tangles of cobordisms. The Frenkel-Stroppel-Sussan projector lives in a category of endofunctors of some category $\OC$, and is in a sense Koszul dual to the other two \cite{SS12}.   Similar constructions can be made in categories of $\sl_2$ matrix factorizations as defined in \cite{KhoRoz08}, or directly with the quotient of the Soergel category considered in \cite{ETL}. 

The categorified Jones-Wenzl projectors were later extended to a categorification of all the primitive idempotents $p_T\in \TL_n$ by the second author and Ben Cooper \cite{CH12}, where $T$ ranges over the SYT with two-row shapes.  Reassociation (see \S \ref{subsec-idempotentreassoc}) yields complexes $\PB_{n,k}\in \KC^-(\BN_n)$ which categorify central idempotents in $\TL_n$; these are naturally indexed by pairs of integers $0\leq k\leq n$ with the same parity, or equivalently two-row partitions $(\frac{n+k}{2},\frac{n-k}{2})$.  As a special case, one recovers both the categorified Jones-Wenzl idempotents $\PB_{n,n}$ of Cooper-Krushkal-Rozansky, and also the ``minimal weight'' idempotents $\PB_{2n,0}$ constructed earlier by Rozansky \cite{Roz10b}.

Let us discuss the idempotents $\PB_{n,0}$ and $\PB_{n,n}$ further.  To distinguish them we will refer to $\PB_{n,n}$ as the ``top projector'' and $\PB_{n,0}$ (or $\PB_{n,1}$ if $n$ is odd) as the ``bottom projector.''   Rozansky constructs these idempotents as infinite limits of powers of the full twist complex in the setting of Khovanov homology (in the case of the bottom projector, Rozansky considers only even $n$).  

In general, the $\PB_{n,k}$ are conjectured in \cite{CH12} to be eigenprojections for the full twist $\FT_n$ acting on $\BN_n$.  Given this, the existence of Rozansky's two limits could then be explained by the fact that the top and bottom projectors correspond to the maximal and minimal eigenvalues of $\FT_n$ in a certain sense.  It is more-or-less obvious from the construction that for each $T$, the idempotent $\idemp_T$ in \cite{CH12} is constructed as an infinite convolution constructed from shifted copies of a complex $\QB_T$ which itself is an eigenobject for the full twist, in the weak sense.  However, this is not sufficient to conclude that the total complex $\idemp_T$ is an eigenobject, since the category of weak eigenobjects with a given eigenvalue need not be triangulated.

What all these previous works lacked is eigenmaps, which give the true connection between the idempotents and the full twist.   If one could show that $\QB_T$ is an eigenobject for some eigenmap $\a_{n,k}:\l_{n,k}\rightarrow \FT_n$ (where $\l_{n,k}$ is the appropriately shifted copy of the identity object), then this would prove that the idempotents in \cite{CH12} define a diagonalization of $\FT_n$.    As it turns out, when constructing eigenmaps it is easier to work directly with the category of Soergel bimodules than with the Temperley-Lieb ``quotient'' $\BN_n$.  Fortunately, our results for $\SBim_n$ will carry over to the $\sl_2$ (in fact $\sl_N$) specializations, thanks to the existence of a monoidal functor relating these categories.  We remark that the eigenmaps $\a_{n,n}$ and $\a_{n,0}$ can be constructed by hand, and so the top and bottom projectors are already known to be eigenobjects for $\FT_n$.  

Finally, let us remark how eigenmaps shed light on Rozansky's construction.  First, suppose $\a:\one \rightarrow \FT_n$ is a map of degree $\l$, where $\FT_n\in \KC^b(\BN_n)$.  Then one can form a directed system (omitting shifts)
\[
\one\rightarrow \FT_n\rightarrow \FT_n^{\otimes 2}\rightarrow \cdots
\]
in which all the maps are $\a\otimes \Id^{\otimes k}$.  There is an abstract process (homotopy colimit, or mapping telescope) by which one can take the limit of such a sequence.  This colimit always exists, but lives usually in some formal extension of the original category.  If $\a=\a_{n,n}$ is the ``top'' eigenmap, then the homotopy colimit is the categorified Jones-Wenzl projector $\PB_{n,n}$, which lives in $\KC^-(\BN_n)$.  At the other end, if $\a=\a_{n,0}$ is the ``bottom'' eigenmap, then the homotopy colimit is Rozansky's bottom projector, which lives in $\KC^+(\BN_n)$.  We denote this by $\PB_{n,0}^\vee$, since we reserve the notation $\PB_{n,0}$ for its dual, which lives in $\KC^-(\BN_n)$. 

Conjecturally, the homotopy colimit above is nice for any of the eigenmaps.

\begin{conjecture}
If $\a:\one\rightarrow \FT_n$ is an eigenmap then the homotopy colimit $\LB_\a$ of $\one\rightarrow \FT_n\rightarrow \FT_n^{\otimes 2}\rightarrow \cdots$ is isomorphic to the internal endomorphism ring of the corresponding projector.  In particular $\Homg(\one,\LB_\a)\cong \Endg(\PB_\a)$.
\end{conjecture}

This conjecture is true in case of the maximal and minimal eigenvalues.

\begin{remark}
There are many more categorified idempotents \cite{Rose14,CautisClasp,Hog15,AbHog15} that can be constructed via analagues of Rozansky's argument, and one expects each of these to fit into a similar story as the above.   In particular all of these idempotents should be extended to give a complete collection of eigenprojections for the appropriate full twist complex in each context.  We take this up in a sequel.
\end{remark}

The theory of diagonalization produces also finite, quasi-idempotent complexes which categorify the Jones-Wenzl idempotents in which certain denominators have been cleared (see \S \ref{subsec:quasiIdemp1}).  These complexes are much more subtle than their idempotent cousins, as they cannot be constructed by the infinite full twist construction.  Rather, they owe their existence to the presence of a saturated collection of eigenmaps.  These were constructed in special cases already by the second author in \cite{H14a,Hog15}.  We expect that these complexes give rise to a functorial categorification of the Reshetikhin-Turaev link invariants.  Such a theory will be essential in any eventual applications of link homology to 4-dimensional topology.  We consider this in future work.

\part{Categorical Diagonalization}

\section{Convolutions and homotopy categories}
\label{sec-homotopy}

We state our conventions for complexes in \S\ref{subsec-conventions}. In \S\ref{subsec-convolutions} through \S\ref{subsec-convolutions3}, we recall the standard notion of convolution in
a homotopy category, and state their major properties, omitting proofs. We use some particularly complicated manipulations of infinite complexes in this paper, so we have included this
extremely careful background section to set the reader at ease.

In \S\ref{subsec:triang} and \S\ref{subsec-nothomotopy}, we discuss to what degree the results in this paper will remain true in more general triangulated categories. The reader willing
to stick with complexes is welcome to skip these sections.

Throughout this paper, we use the following convention: $\BC$ will always represent an additive monoidal category, and $\VC$ an additive category on which $\BC$ acts. Analogously, $\AC$ will always represent a triangulated monoidal category, and $\TC$ a triangulated category on which $\AC$ acts. If we want to represent a monoidal category which is either additive or triangulated, we use $\CC$.

\subsection{Conventions for complexes}
\label{subsec-conventions}

Let $\VC$ be an additive category, and $\Ch(\VC)$ the category of complexes over $\VC$. The differential in a complex will always increase the homological degree by 1.  If $A,B$ are complexes, we let $\Hom^k(A,B)$ denote the space of linear maps which are homogeneous of homological degree $k$.  Then $\Homb(A,B)=\bigoplus_{k\in\Z}\Hom^k(A,B)$ is a chain complex with differential
\[
[d,f]:=d_B\circ f - (-1)^{|f|} f\circ d_A,
\] 
where $|f|=k$ is the homological degree of $f\in \Hom^k(A,B)$.
A degree zero cycle in $\Homb(A,B)$ is precisely a chain map. Two chain maps $f,g\in \Homb(A,B)$ are \emph{homotopic}, written $f\simeq g$, if $f-g$ is a boundary.  The \emph{homotopy category} $\KC(\VC)$ is the category with the same objects as $\Ch(\VC)$, whose morphism spaces are the degree zero chain maps modulo homotopy.  Let $\KC^{+}(\VC)\subset \KC(\VC)$ (resp. $\KC^-, \KC^b$) denote the full subcategories consisting of complexes which are bounded below (resp. bounded above, bounded).  We write $A \simeq B$ when two complexes are homotopy equivalent (i.e.~isomorphic in $\KC(\VC)$).

\begin{definition} \label{def:homshift} Let $A[1]$ denote the complex obtained from $A$ by shifting down in homological degree and negating the differential. Thus if $A$ is concentrated in homological degree zero, then $A[1]$ is in homological degree $-1$. For a map $f \in \Homb(A,B)$, let $f[1]$ denote the same linear map, viewed as a map in $\Homb(A[1],B[1])$. We often refer to $f[1]$ simply as $f$ when the source and target are understood.
\end{definition} 

A degree $k$ cycle will be called a degree $k$ chain map. The fact that $[1]$ negates the differential implies that a degree zero cycle $A\rightarrow B[k]$ is equivalent to a degree $k$ cycle $A\rightarrow B$, which is equivalent to a degree 0 cycle $A[-k]\rightarrow B$ .  Hence there is no harm in identifying
\[
\Homb(A,B[k])  =\Homb(A,B)[k] = \Homb(A[-k],B).
\]

Given a chain map $f:X\to Y$, the mapping cone is defined as the complex $\Cone(f)=X[1]\oplus Y$ with differential $\smMatrix{-d_X &0\\ f & d_Y}$. The inclusion map $\iota \co Y \to
\Cone(f)$ and projection map $\pi \co \Cone(f) \to X[1]$ are both chain maps.

The homotopy category $\KC(\VC)$ is the prototypical example of a triangulated category (see \S\ref{subsec:triang} for a review of triangulated categories). By definition, a triangle in $\KC(\VC)$ is distinguished if and only if it is isomorphic to a triangle of the form
\[
X\buildrel f\over \longrightarrow Y\buildrel \iota\over\longrightarrow \Cone(f)\buildrel \pi\over \longrightarrow X[1],
\]
where the maps $\iota$ and $\pi$ are the canonical inclusion and projection maps. The full subcategories $\KC^{+,-,b}(\VC)\subset \KC(\VC)$ inherit triangulated structures from $\KC(\VC)$.

Let $(\BC,\otimes,\one)$ be an additive monoidal category.  Then $\KC^\pm(\BC)$ inherits the structure of a triangulated monoidal category.  The tensor product of complexes is defined by the usual rule:
\begin{equation} \label{eq:tensorproductrule}
(X\otimes Y)_k = \bigoplus_{i+j=k} X_i\otimes Y_j \ \ \ \ \ \ \  \ \text{with differential} \ \ \ \ \ \ \ \ \  d_{X_i\otimes Y_j} = d_{X_i}\otimes \Id_{Y_j} +(-1)^i \Id_{X_i}\otimes d_{Y_j}.
\end{equation}
The direct sum above is finite, given the restriction to semi-infinite complexes. 

Suppose $\BC$ is an additive monoidal category and $\VC$ is an $\BC$-module category.  Then for each  $\circ\in \{+,-,b\}$ we define an action of $\KC^\circ(\BC)$ on $\KC^\circ(\VC)$ via the usual formula for tensor product of complexes: if $C\in \K^\circ(\BC)$ and $M\in \KC^\circ(\VC)$, then define a complex $C\otimes M\in \K^\circ(\VC)$ by
 \[
 (C\otimes M)_k = \bigoplus_{i+j=k}C_i\otimes M_j
 \]
 with differential which is the sum of $d_{C_i}\otimes \Id_{M_j} + (-1)^{i} \Id_{C_i}\otimes d_{M_j}$ over all $i,j\in \Z$.  The above direct sum is finite given the boundedness conditions on $C,M$. 

The above action is compatible with mapping cones in the following sense: given any chain map $f:A\rightarrow B$ in $\KC^\circ(\BC)$ and any complex $M\in \KC^\circ(\VC)$, we have $\Cone(f\otimes \Id_M)\cong \Cone(f)\otimes M$.  Similarly, given any chain map $g:M\rightarrow N$ in $\KC^\circ(\VC)$ and any complex $A\in \KC^\circ(\BC)$, we have $\Cone(\Id_A\otimes g)\cong A\otimes \Cone(g)$.  In other words $A\otimes(-)$ and $(-)\otimes M$ are exact functors.

\subsection{Twisted complexes and convolutions}
\label{subsec-convolutions}

In this section we recall the basic construction known as convolution, inside a triangulated category $\TC$. A convolution may be thought of as an iterated mapping cone, or as a kind of
filtered complex. For now we define it under the assumption that $\TC \subset \KC(\VC)$ is a full triangulated subcategory of a homotopy category. In \S\ref{subsec-nothomotopy} we will discuss what can be done in general.

\begin{defn}\label{def:twistedCx}
Let $(I,\leq)$ be a partially ordered set, which for sake of simplicity we always assume is \emph{interval finite}: for each $j \le i$ in $I$, there are a finite number of $k \in I$ with $j \le k \le i$.  Let $\VC$ be an additive category, and $\TC \subset \KC(\VC)$.  An \emph{$I$-indexed twisted complex} in $\TC$ is a collection $\{A_i,d_{ij}\}$ where $A_i\in \TC$ are complexes indexed by $i \in I$, and $d_{ij}\in\Hom^1(A_j,A_i)$ are degree 1 linear maps with the following properties:
\begin{enumerate}
\item $d_{ij} = 0$ unless $j \le i$.
\item $d_{ii} = d_{A_i}$.
\item For each pair $j \le i$, we have $\sum_{j\le k\le i} d_{ik}d_{kj} = 0$.
\end{enumerate}
We think of $d_{ij}$ as a lower-triangular $I \times I$ matrix. We refer to each complex $A_i$ as a \emph{layer} of the twisted complex. \end{defn}

For example, when $I = \ZM$ with its usual total order, and $d_{ij}=0$ unless $i = j+1$, one obtains the definition of an (anti-commuting) bicomplex. The generalization of the total complex of a bicomplex is a convolution.

\begin{defn}\label{def:conv} Consider an $I$-indexed twisted complex $\{A_i,d_{ij}\}$. Suppose that, for each $j \in I$, we have $d_{ij}=0$ for all but finitely many $i$, so that each column of $d_{ij}$ has finitely many non-zero entries. Then (after possibly enlarging $\TC$) there is a well-defined complex
\[
\Tot^{\oplus}\{A_i,d_{ij}\} = \bigoplus_{i\in I} A_i  \ \ \ \ \text{with differential $d = \sum d_{ij}$}.
\]
Alternatively, suppose that for each $i \in I$ we have $d_{ij}=0$ for all but finitely many $j$, so that each row of $d_{ij}$ has finitely many non-zero entries. Then there is a well-defined complex
\[
\Tot^{\Pi}\{A_i,d_{ij}\} = \prod_{i\in I} A_i  \ \ \ \ \text{with differential $d = \sum d_{ij}$}.
\]
We refer to these total complexes as \emph{convolutions} of the twisted complexes $\{A_i,d_{ij}\}$.
\end{defn}
These total complexes may only be well-defined if $\TC$ contains certain countable direct sums or products (if not, we think of it as being defined over an enlargement of $\TC$). In the exposition below, we may assume tacitly that convolutions exist. We discuss the existence of convolutions further in \S\ref{subsec-locallyfinite}.

\begin{remark} The total complex $\Tot^{\oplus}$ of a twisted complex over an abelian category is the same data as a filtered complex which, upon forgetting differentials, splits in each
homological degree.\end{remark}

\begin{remark} The direct sum $\bigoplus_{i \in I} A_i$ of complexes can always be given a boring direct sum differential $\sum d_{A_i}$. A convolution endows $\bigoplus_{i \in I} A_i$
with a new, interesting differential, which is called ``twisted" because of how it differs from the boring one. It should be noted that the literature (c.f. \cite{BondalKap90}) has
primarily focused on the situation where the index set $I$ is finite, but does not assume that the differential respects a partial order on $I$. We are interested in infinite indexing
sets, but our partial order allows the convolution to remain well-behaved in many circumstances. See \S\ref{subsec-nothomotopy} for more. \end{remark}

\subsection{Local finiteness}
\label{subsec-locallyfinite}

In general, the formation of $\Tot^\oplus\{A_i, d_{ij}\}$ may require the existence of infinite direct sums. In practice, it is often the case that only finite direct sums are required, even when $I$ is infinite.

\begin{definition} \label{defn:locfinite} An $I$-indexed twisted complex $\{A_i,d_{ij}\}$ is $\TC$-\emph{locally finite} if $\Tot^\oplus\{A_i,d_{ij}\}$ and $\Tot^\Pi\{A_i,d_{ij}\}$ exist in $\TC$ and are isomorphic. \end{definition}
Note that an $I$-indexed twisted complex $\{A_i, d_{ij}\}_{i\in I}$ is $\TC$-locally finite iff $\{A_i, 0 \}_{i\in I}$ is $\TC$-locally finite.  Thus, local finiteness is a property of the complexes $\{A_i\}$ and not of the twisted differential $d_{ij}$.

\begin{example}
An example of local finiteness is the definition of the tensor product of two infinite but bounded above complexes.  Meanwhile, the ``tensor product" of a two infinite complexes, one bounded above and the other bounded below, is an example of a non-locally-finite twisted complex, and for this reason one must distinguish between the direct sum tensor product $\bigotimes^{\oplus}$ and the direct product (or \emph{completed}) tensor product $\bigotimes^\Pi$.
\end{example}

\begin{remark} \label{rmk:noGrothGpPlz1} Let us pause to state our priorities clearly. There are three things one may be concerned with when choosing a subcategory $\TC$ in $\KC(\VC)$, and taking convolutions therein:
\begin{enumerate}
\item whether $\Tot^\oplus$ and $\Tot^\Pi$ even exist for a given twisted complex, i.e. whether certain infinite direct sums or products exist in $\TC$;
\item whether $\Tot^\oplus$ and $\Tot^\Pi$ are actually isomorphic in $\TC$; and
\item whether $\TC$ has an interesting, non-zero Grothendieck group. One might also ask whether infinite complexes in $\TC$ categorify the corresponding infinite sums in the Grothendieck group (which should be convergent series).
\end{enumerate}
Properties (1) and (2) will be essential for various constructions and theorems in this paper. Property (3) is a separate issue entirely, and one that we largely ignore. There are various examples where any category large enough for certain convolutions to exist will have zero Grothendieck group, but this does not bother us. We are concerned with making valid categorical arguments, not with whether these arguments precisely categorify a statement in linear algebra.

\begin{remark}\label{rmk:noGrothGpPlz2} If $\TC$ is a triangulated category, the \emph{triangulated Grothendieck group} of $\TC$ is the abelian group formally spanned by isomorphism classes of objects of $\TC$, modulo the defining relation that $[X]-[Y]+[Z]=0$ for every distinguished triangle $X\rightarrow Y\rightarrow Z\rightarrow X[1]$ in $\TC$.  There is always a distinguished triangle of the form
\[
X\rightarrow 0 \rightarrow X[1]\rightarrow X[1],
\]
which forces the relation $[X]=-[X[1]]$ for all $X\in \TC$. If $\VC$ is an idempotent complete, additive category and $\TC=\KC^b(\VC)$, then $[\TC]$ is isomorphic to the split Grothendieck group $[\VC]$.  The isomorphism sends the class of a complex $[X]\in [\TC]$ to its Euler characteristic $\sum_{i\in \Z}(-1)^i[X_i]\in [\VC]$. On the other hand, the Grothendieck groups of categories of unbounded complexes are more difficult to analyze, and are often zero by a variant on the Eilenberg swindle. See \cite{Miyachi} for further results and discussion in this direction. \end{remark}

For this reason, we have made our definition of local finiteness context-dependent. Those readers who are concerned with the Grothendieck group will have to find appropriate triangulated
categories $\TC$ where the desired convolutions are locally finite, and appropriate versions of the Grothendieck group which behave well for $\TC$. For examples of papers which do
carefully discuss Grothendieck groups of infinite complexes, see \S 2.7.1 of \cite{CK12a}, or see the \emph{topological Grothendieck group} in \cite{AchStr}. \end{remark}

One way to guarantee local finiteness (regardless of the choice of $\TC \subset \KC(\VC)$) is via the following proposition, whose proof is obvious.

\begin{proposition}\label{prop:locfinite}
Let $A_i\in \KC(\VC)$ be given ($i\in \Z$), and denote the $j$-th chain group of $A_i$ by $A_i^j$.  Suppose that for each $j\in \Z$, we have
\[
\bigoplus_{i\in \Z} A_i^j \cong \prod_{i\in \Z} A_i^j
\]
in $\VC$.  Then $\bigoplus_i A_i $ and $\prod_i A_i$ exist in $\KC(\VC)$ and are isomorphic.  In particular $\{A_i\}_{i\in \Z}$ is $\TC$-locally finite for all $\TC\subset \KC(\VC)$.
\end{proposition}

\begin{example}[Homological local finiteness]\label{ex:homLocalFinite}
Let $\{A_i, d_{ij}\}$ be an $I$-indexed twisted complex.  Denote the $j$-th chain group of $A_i$ by $A_i^j$, and suppose that for each $j\in \Z$, $A_i^j=0$ for all but finitely many $i\in I$.  Then the direct sum $\bigoplus_{i\in I} A_i$ is finite in each degree, hence this direct sum exists and is isomorphic to the direct product $\prod_{i\in I} A_i$.  Thus, the twisted complex $\{A_i, d_{ij}\}$ is locally finite.  We call twisted complexes of this sort \emph{homologically locally finite}.
\end{example}

Most of our examples (e.g. tensor products of bounded above complexes) will be homologically locally finite.

\subsection{Diagrams and reassociation}
\label{subsec-convolutions2}

We will visualize twisted complexes and their total complexes with diagrams of arrows.  We illustrate this with examples.

\begin{example}
If $f:B\rightarrow C$ is a degree zero chain map, then its mapping cone can be written as
\[
\Cone(f) = \Tot\Big( B[1] \buildrel f\over \longrightarrow C\Big ).
\]
In this example, $I = \{1,2\}$ with $1 < 2$, and $A_1$ is the complex $B[1]$, while $A_2$ is the complex $C$. The map $d_{21}$ is $f$, viewed as a degree one map $B[1] \to C$, also known as a degree zero map $B \to C$. All of these things (the poset $I$, the complexes $A_i$, the differentials $d_{ij}$) are implicit in the diagram above. Even the fact that $f$ is a chain map is implicit, because this is exactly the condition (3) from Definition \ref{def:twistedCx}. 
We often refer to this condition (3) via its implied equality $d^2 = 0$.
\end{example}

\begin{remark} Our convention for taking total complexes differs from some others used in the literature. For example, the ``usual'' convention for total complexes would write a mapping cone as
\[
\Cone(f) = \Tot^{\text{usual}}\Big( B \buildrel f\over \longrightarrow \underline{C} \Big).
\]
In the complex $\Cone(f)$, it is $B[1]$ which appears, not $B$. That is, the usual convention for total complexes introduces extra homological shifts which are not visible in the notation (the underline $\un{C}$ is there to indicate that this complex appears with no shift). This reasonable convention becomes quite obnoxious when dealing with iterated cones.  For instance, the precise degree shifts involved in the formation of
\[
\Tot^{\text{usual}}(A \rightarrow \Tot^{\text{usual}}(\un{B}\rightarrow C) \rightarrow \un{D})
\]
can be determined but this requires bookkeeping. We prefer to use the simpler convention in which all degree shifts are already accounted for in the diagrams.

The price we pay is that all arrows in our diagrams are degree $+1$ maps, instead of degree zero maps.  We will occasionally remind the reader of our convention by writing, for instance,
\[
\Tot(A_0\buildrel[1]\over \rightarrow A_1\buildrel [1]\over\rightarrow A_2),
\]
to indicate that the connecting differentials are actually linear maps $A_i\rightarrow A_{i+1}[1]$.

The only time we use the usual convention is when we define a complex (not a convolution of complexes). So for example \[ \Big( M_1 \longrightarrow M_2 \longrightarrow \un{M}_3 \Big) \] is a three term complex where $M_1$ appears in homological degree $-2$.

\end{remark}

\begin{example}
Consider the following convolution.
\[
Z = \Tot\left( \begin{diagram}
A[2] &&\\
\dTo^{f} & \rdTo^h  &\\
B[1] &\rTo^{g} & C
\end{diagram} \ \ 
\right)
\]
Implicit here is the poset $I = \{1,2,3\}$, the complexes $A_1 = A[2]$, $A_2 = B[1]$, and $A_3 = C$, and the differential $d_{ij}$. The equation $d_Z^2=0$ implies that $f$ and $g$ correspond to degree zero chain maps $A\rightarrow B$ and $B\rightarrow C$, respectively, and $h$ is a degree $-1$ map $A \to C$ which satisfies $d\circ h+h\circ d=-g\circ f$.

Observe that $Z$ can be written as an iterated mapping cone in two ways: we can group the terms on the bottom together, obtaining
\[
Z \cong \Tot\Big(A[2] \rightarrow \Cone(g) \Big)
\]
We can also group the terms on the left together, obtaining
\[
Z \cong \Tot\Big(\Cone(-f)[1]\rightarrow C\Big)
 \]
The sign on $-f$ is to counteract the sign involved in the degree shift $[1]$.
We have just seen our first example of \emph{reassociation} of convolutions.  We now explain this in more detail.
\end{example}

\begin{defn}\label{def:convex}
Let $I$ be a poset.  A subset $J\subset I$ is called \emph{convex} if whenever $j_1\leq i\leq j_2$ and $j_1,j_2\in J$, then $i\in J$.
\end{defn}

\begin{defn} \label{def:contrib}
If $C = \Tot^{\oplus}\{C_i, d_{ij}\}_I$ is a convolution and $J\subset I$ is convex, then the \emph{contribution of $J$ to $C$} defines a  complex $C_J$.  Precisely, $C_J:=\bigoplus_{j\in J} C_j$ with differential $d_J = \sum_{j,j' \in J} d_{j,j'}$.
\end{defn}

Since $J$ is convex, one can confirm easily that $d_J^2 = 0$.  One thinks of $C_J$ as a subquotient complex of $C$, even when there is no abelian structure in the underlying category (hence no given notion of quotients, images, etc.).

If $I = J_1\sqcup \cdots \sqcup J_r$ is a partition of $I$ into disjoint subsets, one can place a preorder on the set $\{J_1,\ldots,J_r\}$, given by $J_a\leq J_b$ if there exist $i\in J_a$ and $j\in J_b$ such that $i\leq j$. This preorder may or may not be a partial order, but if it is a partial order, then the subsets $J_i$ are convex.   Partitioning $I$ in partially ordered fashion into disjoint, convex sets gives a way of reassociating the terms of $C$.

\begin{definition}\label{def-reassociation} (Reassociation) Let $C = \Tot^{\oplus}\{C_i, d_{ij}\}_I$ be a convolution over $I$.
Let $\PC$ be some partially ordered set and let $I = \bigsqcup_{a \in \PC} J_a$ be a partition of $I$ into disjoint convex subsets $J_a$, respecting the partial order on $\PC$.  Set $C_a := C_{J_a}$ as in Definition \ref{def:contrib}.  We can define a degree 1 map $d^{\prime}_{ab} \co C_b \to C_a$ as the sum of all $d_{ij}$ for $i \in J_a$ and $j \in J_b$.  This produces a new twisted complex $\{C_a, d^{\prime}_{ab}\}_{\PC}$. It is easy to see that $\Tot^{\oplus}\{C_{a},d^{\prime}_{ab}\}_{\PC}$ exists, and
\[
C  = \Tot^\oplus\{C_i, d_{ij}\}_I \cong \Tot^{\oplus}\{C_{a}, d'_{ab}\}_{\PC}.
\]
We call the twisted complex $\{C_{a}, d'_{ab}\}_{\PC}$ or its convolution a \emph{reassociation of $C$}.
\end{definition}

The following is a special example of reassociation.

\begin{lemma}[Linearity of convolutions] Let $C = \Tot^{\oplus}\{C_i, d_{i,i'}\}_I$ and $D=\Tot^\oplus\{D_j,e_{j,j'}\}_J$ be convolutions in a monoidal category.  Then the (direct sum) tensor product $C \ot D$ can also be described as a convolution $C \ot D = \Tot^{\oplus}\{C_i \ot D_j\}_{I\times J}$.  Reassociating appropriately, this can also be described as $C\otimes D \cong \Tot\{C_i\otimes D\}_I$ or $C\otimes D\cong \Tot\{C\otimes D_j\}_J$.  \end{lemma}

\subsection{Simultaneous simplifications}
\label{subsec-convolutions3}

The following crucial proposition states that convolutions are ``invariant" under homotopy equivalence.

\begin{defn} A poset is \emph{upper-finite} if it satisfies the ascending chain condition (ACC): any chain $i_1 < i_2 < \ldots$ is finite. A poset is \emph{lower-finite} if it satisfies the descending chain condition (DCC). \end{defn}

\begin{proposition}\label{prop:simplifications} (Simultaneous simplifications) Let $(I,\leq)$ be a poset and $\{C_i,d_{ij}\}_I$ a twisted complex. Suppose we are given complexes $D_i$ which are homotopy equivalent to $C_i$, for each $i \in I$. If $I$ is upper-finite then there is some twisted complex $\{D_i,d'_{ij}\}_I$ such that
$\Tot^\oplus\{C_i,d_{ij}\}_I\simeq \Tot^\oplus\{D_i,d'_{ij}\}_I$. If instead $I$ is lower-finite then the same result holds with direct sum replaced by direct product.
\end{proposition}
\begin{proof}
Follows from homotopy invariance of mapping cones and a straightforward limiting argument.
\end{proof}

That is, given a convolution, one can ``replace'' a layer of the convolution with a homotopy equivalent layer, without affecting the homotopy equivalence class of the convolution. This
does change the twisted differentials in a mysterious way which in theory can be computed, though in practice this is difficult.

In particular, suppose that each layer $C_i$ of a (direct sum) convolution $C$ is null-homotopic. If $I$ is upper-finite then $C$ is also null-homotopic. If $I$ is lower-finite but the
twisted complex is $\TC$-locally finite, then $C$ is null-homotopic because it agrees with the direct product convolution.

The following pair of examples illustrate the danger of misusing Proposition \ref{prop:simplifications}. 
\begin{example}
Consider the chain complex
\[
C \ \ : = \ \ \Z[x]\buildrel 1-x\over \longrightarrow \underline{\Z[x]},
\]
thought of as a complex of abelian groups with the underlined term in degree zero.  This complex can be thought of as the total complex of a first quadrant bicomplex
\[
\begin{diagram}[small]
&&&&\cdots\\
&&&&\dTo^{1}\\
&&\Z & \rTo^{-1} & \Z&&\\
&& \dTo^{1}&&\\
\Z & \rTo^{-1} & \Z &&\\
\dTo^{1}&&\\
\Z
\end{diagram}
\]
Since $\Z[x]$ consists of polynomials in $x$, rather than power series, $C$ is isomorphic to $\Tot^\oplus$ of this bicomplex, rather than $\Tot^\Pi$.  This total complex may be reassociated in two natural ways:
\begin{subequations}
\begin{equation}\label{auxeq:misleading1}
C = \Tot^{\oplus}(\Cone(\Id_\Z)\buildrel[1] \over \rightarrow \Cone(\Id_\Z)\buildrel[1] \over \rightarrow \cdots )
\end{equation}
\begin{equation}\label{auxeq:misleading2}
C = \Tot^{\oplus}(\cdots \buildrel[1] \over \rightarrow \Cone(-\Id_\Z)\buildrel[1] \over \rightarrow \Cone(-\Id_\Z) \buildrel[1] \over \rightarrow \Z)
\end{equation}
\end{subequations}
Since $\Cone(\Id_{\Z})$ is contractible, one may be tempted to cancel all the terms in \eqref{auxeq:misleading1}, but to do so would be a misuse of Proposition \ref{prop:simplifications}.  Since we are taking $\Tot^\oplus$ we must instead work with the upper-finite convolution \eqref{auxeq:misleading2} rather than the lower-finite convolution \eqref{auxeq:misleading1}.  Contracting the contractible terms yields
\[
\Cone(\Z[x]\buildrel 1-x \over \rightarrow \Z[x]) \simeq \Z.
\]
\end{example}
\begin{example}
Consider the above example with $\Z[x]$ replaced by the power series ring $\Z\llbracket x\rrbracket$.  The occurences of $\Tot^\oplus$ get replaced by $\Tot^\Pi$, and the lower-finite version of Proposition \ref{prop:simplifications} now applies, yielding
\[
\Cone(\Z\llbracket x\rrbracket] \buildrel 1-x \over \rightarrow \Z\llbracket x\rrbracket ) \simeq 0.
\]
Indeed, $1-x$ is invertible in the power series ring, hence the above is the cone on an isomorphism, which is contractible.
\end{example}


\subsection{General triangulated categories}
\label{subsec:triang}

Many of our crucial constructions in this paper are given as convolutions of twisted complexes, e.g. \eqref{eq:introCab}. When analyzing these constructions we will use reassociation (Definition \ref{def-reassociation}) and simultaneous simplifications (Proposition \ref{prop:simplifications}) to great effect. It is worthwhile asking what can be done in a general triangulated category, and we devote the remainder of this chapter to this discussion. The uninterested reader is welcome to skip ahead. We begin by recalling the basics of triangulated (monoidal) categories.

Let $\TC$ be an additive category with an invertible endofunctor $[1]:\TC\rightarrow \TC$, called the \emph{suspension}.  A \emph{triangle} in $\TC$ is a triple $(f,g,h)$ of composable morphisms
\begin{equation}\label{eq:triangle}
\begin{diagram}[small]
X &\rTo^{f} & Y & \rTo^g & Z & \rTo^g & X[1].
\end{diagram}
\end{equation}
A \emph{triangulated category} is an additive category $\TC$ together with a suspension functor $[1]$ and a collection of distinguished triangles satisfying a number of axioms (we adopt the axioms in \cite{MayTrace}).


If $f:X\rightarrow Y$ is a morphism in a triangulated category, then there is an object $Z$ which fits into a distinguished triangle with $f$ as in \eqref{eq:triangle}.  The object $Z$ is unique up to non-unique isomorphism; we call $Z$ the \emph{mapping cone on $f$}, and we write $Z\cong \Cone(f)$.

In order to categorify rings, we will want to consider triangulated categories that have a compatible notion of a tensor product. There seems to be no accepted definition of triangulated monoidal categories in general, though in the symmetric monoidal case there are good definitions in \cite{MayTrace,Balmer10}. We are interested in the non-symmetric monoidal case, and for this purpose it will suffice to adopt Definition 3.8 in \cite{Hog17a}.

Thus we define a \emph{triangulated  monoidal category} to be an additive category $\TC$ with the structure of a triangulated category $(\TC,[1],\{\text{triangles}\})$ and a monoidal category $(\TC,\otimes,\one)$\footnote{We will regard the associators and unitors as implicit, and omit them from the notation.}, together with natural isomorphisms
\[
\a_X: X\otimes \one[1] \cong X[1] \ \ \ \ \ \ \ \ \ \ \ \ \b_X: \one[1]\otimes  X\cong X\otimes \one[1]
\]
satisfying a number of coherence relations.  In addition to the coherence relations, we will also require that $\otimes$ be compatible with mapping cones, in the sense that
\[
X\otimes \Cone(f) \cong \Cone(\Id_X\otimes f) \ \ \ \ \ \text{and} \ \ \ \ \ \ \ \Cone(f)\otimes X \cong \Cone(f\otimes\Id_X).
\]
This is summarized by saying that $X\otimes(-)$ and $(-)\otimes X$ are exact functors $\TC\rightarrow \TC$ for all $X\in \TC$. The data above gives rise to natural isomorphisms
\[
X[1]\otimes Y \cong (X\otimes Y)[1]\cong X\otimes Y[1]
\]
for any objects $X$ and $Y$.

Among the coherence relations, it will be assumed that $\b_{\one[1]}$ is the negative (!) of the identity map of $\one[1]\otimes \one[1]$; this sign will be the ultimate reason for various signs later.

\begin{example}
Let $(\BC,\otimes,\one)$ be an additive monoidal category, and consider the triangulated monoidal category $\KC^\pm(\BC)$. Using the usual sign rule for tensor products \eqref{eq:tensorproductrule}, one sees that for any complex $X\in \KC^\pm(\BC)$, $X\otimes \one[1]$ can be naturally identifed (after a shift in homological grading) with $(X,d_X)$, while $\one[1]\otimes X$ can be identitified with $(X,-d_X)$.  Thus, the isomorphism in commuting $\one[1]$ past $X$ must involve a sign.  Naturality of this isomorphism forces the sign in commuting $\one[1]$ past $\one[1]$, as mentioned above.
\end{example}

We wish to consider actions of triangulated monoidal categories on triangulated categories.  If $(\AC,\otimes, \one)$ is an additive monoidal category, then an \emph{$\AC$-module category} is an additive category together with an additive monoidal functor $\AC\rightarrow \End(\MC)$ into the category of endofunctors of $\MC$ \cite{EGNO15}.  The action of $\AC$ on $\MC$ will be denoted $(A,M)\mapsto A\otimes M$ for $A\in \AC$, $M\in \MC$.

If $\AC$ has the structure of a triangulated monoidal category and $\MC$ has the structure of a triangulated category, then a \emph{triangulated $\AC$-module category} is an $\AC$-module category $\MC$ together with a natural isomorphism $\one[1]\otimes M\cong M[1]$ and $\otimes M[1]$ satisfying certain coherence conditions, such that $A\otimes (-)$ and $(-)\otimes M$ are exact functors for all $A\in \AC$, $M\in \MC$.  Any action of a triangulated monoidal category on a triangulated category will be assumed to be triangulated in this sense.

\subsection{Convolutions and enhanced triangulated categories}
\label{subsec-nothomotopy}

Suppose $\TC$ is a triangulated category, and $Y_0\rightarrow Y_1\rightarrow \cdots \rightarrow Y_r$ is a chain complex over $\TC$ (i.e.~an object of $\KC^b(\TC)$, in which the differential is composed of maps $\d_i:Y_i\rightarrow Y_{i+1}[1]$.  A \emph{left Postnikov system} attached to $Y_\bullet$ is a sequence of objects $Y_0=X_0,\ldots,X_r$ with maps $\sigma_i:X_i\rightarrow X_{i-1}$ and distinguished triangles
\[
\begin{diagram} Y_i  & \rTo^{\xi_i} &  X_i &\rTo^{\sigma_i}&  X_{i-1} &\rTo^{\pi_i}& Y_i[1]\end{diagram}
\]
such that $\d_i = \pi_{i+1}\circ \xi_i$.  In this case $X_r$ is called a \emph{convolution} of $Y_\bullet$.  The distinguished triangles above imply that $X_r$ is an iterated mapping cone:
\[
X_r\cong  \Cone\bigg( \cdots \Cone\Big(\Cone(Y_0[-1]\rightarrow Y_1)[-1]\rightarrow Y_2\Big)[-1]\rightarrow  \cdots \rightarrow Y_r\bigg)
\]
In particular $[X_r]=\sum_i [Y_i]$ in the Grothendieck group.  Just as the total complex of a bicomplex realizes a ``flattening'' of a bicomplex into a complex, one thinks of convolutions as realizing a ``flattening'' of a complex of objects of $\TC$ into an object of $\TC$.  In contrast with the total complex of bicomplexes, however, convolutions need not exist\footnote{Higher Massey products involving the $\d_i$ give obstructions to the existence of a convolution of $Y_\bullet$.} and they need not be unique.

In this paper we must consider generalized Postnikov systems, where the $Y_i$ are indexed by partially ordered
set $I$. Then a convolution is a family of objects $X_J$ indexed by ideals $J\subset I$, together with some additional data (restriction maps $X_{J_2}\rightarrow X_{J_1}$ if $J_1\subset
J_2$, together with a family of distinguished triangles), from which $Y_\bullet$ can be recovered. However, there does not seem to be a good theory of such objects, unless $\TC$
possesses some sort of \emph{enhancement}. An enhancement is, roughly, a choice of some sort of category $\UC$ with more structure, yielding $\TC$ upon forgetting that structure. For
instance, realizing $\TC$ as a homotopy category of complexes is a special case of a \emph{dg enhancement} of $\TC$.

The more general notion of dg enhancement is due to Bondal-Kapranov in \cite{BondalKap90}, in which the concept of a \emph{pretriangulated category} is introduced.  A pre-triangulated category $\UC$ is a certain kind of dg category whose homotopy category $H_0(\UC)$ is triangulated, and $\UC$ is said to be a dg enhancement of $H_0(\UC)$.

In a pretriangulated category, any twisted complex with finite index set (and no partially ordered assumption) can be ``flattened'' into its total complex. In this setting, generalized
Postnikov systems and convolutions correspond to the twisted complexes indexed by finite posets. It is relatively straightforward to generalize the notion of a convex set to a given
twisted complex, and to show that the reassociation operation of Definition \ref{def-reassociation}, which is crucial in this paper, will make sense in that context. Taking (co)limits of
this construction, one should be able to discuss twisted complexes with infinite indexing sets which have partial orders. We have not found this written explicitly in the literature,
however.  Simultaneous simplifications is also missing from the literature, but should work for pretriangulated categories. We expect the results of this paper to
extend to the pretriangulated categories of \cite{BondalKap90} (or other suitably enhanced triangulated categories), and it is even possible that a weaker structure will also suffice.

We often just say that $\TC$ is a \emph{homotopy category} to mean that it is a full triangulated subcategory of $\KC(\VC)$ for some additive category $\VC$. Our approach in this paper
is to work within arbitrary triangulated categories when possible, but when we require the full power of reassociation and simultaneous simplifications, to restrict to homotopy
categories. Again, we expect that we are playing it too safe, and that pre-triangulated categories are acceptable replacements for homotopy categories.

\section{Decompositions of identity}
\label{sec-decompositions}

We introduce idempotent decompositions of identity in \S\ref{subsec-idempotents}, and discuss their mutual refinements in \S\ref{subsec-idempotentreassoc}. In \S\ref{subsec:catidemp} we
discuss some abstract properties (for instance, uniqueness) of categorical idempotents. In \S\ref{subsec-idempotents} and \S\ref{subsec-idempotentreassoc} we work in the context of a
homotopy categories, because we use twisted complexes and convolutions. In \S\ref{subsec:generalDecomp} we discuss an equivalent definition which works for arbitrary triangulated
monoidal categories.

One can find a chapter on unital and counital idempotents in \cite{BD14} (see also the preprint \cite{BoyDrin-idemp}, where they are called open and closed idempotents). Many
properties of categorical idempotents are established in \cite{Hog17a}, and additional background information can be found there.

\subsection{Unital and counital idempotents}

Let $\AC$ be a monoidal category.  A \emph{counital idempotent} is an object $\PB\in \AC$ and a map $\e:\PB\rightarrow \one$ such that $\e\otimes \Id_\PB$ and $\Id_\PB\otimes \e$ are isomorphisms in $\AC$.  Similarly, a \emph{unital idempotent} is an object $\QB\in \AC$ with a map $\eta:\one\rightarrow \QB$ such that $\Id_\QB\otimes \eta$ and $\eta\otimes \QB$ are isomorphisms.

\begin{definition}
Let $(\AC,\otimes,\one)$ be a triangulated monoidal category.  We say $(\PB,\QB)$ is a \emph{pair of complementary idempotents} if $\PB\otimes \QB\cong 0 \cong \QB\otimes\PB$ and there exists a distinguished triangle of the form
\begin{equation}\label{eq:idempTriang}
\begin{diagram}[small] \PB &\rTo^{\e} & \one & \rTo^{\eta} & \QB & \rTo^{\d} & \PB[1]
\end{diagram}.
\end{equation}
In this case we also write $\QB = \PB^c$ and $\PB=\QB^c$.
\end{definition}
It is an easy exercise to show that for any choice of distinguished triangle \eqref{eq:idempTriang}, $(\PB,\e)$ and $(\QB,\eta)$ are a counital and a unital idempotent, respectively.  Conversely, any unital (resp.~counital) idempotent has a counital (resp.~unital) complement.  For example, if $(\PB,\e)$ is a counital idempotent, then $\QB:=\Cone(\e)$ has the structure of a unital idempotent, while if $(\QB,\eta)$ is a unital idempotent, then $\Cone(\eta)[-1]$ has the structure of a counital idempotent.

It turns out if $\PB\otimes \QB\cong 0\cong \QB\otimes \PB$, then any triangle of the form \eqref{eq:idempTriang} is determined by $\PB$ or $\QB$ up to unique isomorphism of triangles.  In particular, $(\PB,\e)$ determines the complement $(\QB,\eta)$ canonically, and vice versa (see \S 4.5 in \cite{Hog17a}).  Given the uniqueness of the triangle \eqref{eq:idempTriang}, we will usually  omit the structure maps $\e$ and $\eta$ from the notation, referring to the objects $\PB$ and $\QB$ as (co)unital idempotents.

\begin{example} Recall the complexes $\firstProj$ and $\secondProj$ from \S\ref{subsec:interpolation_example}. Then $\firstProj$ is a counital idempotent, and $\secondProj$ is its complementary unital idempotent. \end{example}

\begin{lemma} \label{lem:HMCexample} Let $(\BC,\otimes,\one_\BC)$ be an additive monoidal category, and let $\AC'\subset \KC^\circ(\BC)$ be a full triangulated monoidal subcategory, with
$\circ \in \{b,+,-\}$. Let $\idemp \in \AC'$ be a unital (resp. counital) idempotent. Let $\AC = \idemp \AC'\idemp\subset \AC'$ be the full subcategory consisting of complexes $C$ such
that $\idemp\otimes C\simeq C$ and $C\otimes \idemp \simeq C$. Then $\AC$ is triangulated, and has the structure of a monoidal category with monoidal identity $\one_\AC:=\idemp$.
\end{lemma}

This lemma is relatively straightforward, and applies more generally to the categorical idempotents of Definition \ref{def:categoricalIdemp}. In particular, it gives a host of new
examples of monoidal homotopy categories (c.f. Definition \ref{def:homotopyMonoidalCat}).

\subsection{Idempotent decompositions of identity}
\label{subsec-idempotents}


Let $(\AC,\otimes, \one)$ be a monoidal homotopy category, unless otherwise indicated.

\begin{definition}\label{def-resOfId} A \emph{decomposition of identity} in $\AC$ is an equivalence $\one\simeq \Tot\{\idemp_i,d_{ij}\}$ where $\{\idemp_i,d_{ij}\}$ is some twisted
complex indexed by a finite poset $(I,\leq)$. We call this an \emph{idempotent decomposition} \emph{of identity} if, in addition, $\idemp_i$ is nonzero for each $i$, $\idemp_i\otimes
\idemp_j\simeq 0$ for $i\neq j$, and $\idemp_i\otimes \idemp_i\simeq \idemp_i$. \end{definition}

Because the poset $I$ is finite, $\Tot^{\oplus}$ and $\Tot^{\Pi}$ automatically agree, so we just write $\Tot$.

The condition that $\idemp_i \ot \idemp_i \simeq \idemp_i$ is redundant. Indeed, for each $i\in I$, we have
\begin{equation} \label{eq:pipiredundant}
\idemp_i \simeq \one\otimes \idemp_i \simeq \Tot\{\idemp_j\}_{j\in I}\otimes \idemp_i \simeq \Tot\{\idemp_j\otimes \idemp_i \}_{j\in I}\simeq \idemp_i\otimes \idemp_i.
\end{equation}
In the last step, we used the fact that $\idemp_j \ot \idemp_i \simeq 0$ unless $j=i$, and we used simultaneous simplifications (Proposition \ref{prop:simplifications}) to replace the total complex of $\{\idemp_j \ot \idemp_i\}$ with its only non-null-homotopic constituent $\idemp_i \ot \idemp_i$.
 
\begin{proposition}\label{prop:shortResolutions}
In any short (i.e.~two term) idempotent decomposition of identity \begin{equation} \label{eq:short1} \one\simeq \Tot(\idemp_1\buildrel [1]\over \rightarrow \idemp_2),\end{equation} $\idemp_1$ has the structure of a  unital idempotent and $\idemp_2=\idemp_1^c$ is the complementary counital idempotent.  
\end{proposition}
\begin{proof}
An equivalence $\one\simeq \Tot(\idemp_1\buildrel [1]\over \rightarrow \idemp_2)$ corresponds to a distinguished triangle
\[
\idemp_2\rightarrow \one\rightarrow \idemp_1\rightarrow \idemp_2[1].
\]
The assumption that $\idemp_1\otimes \idemp_2\simeq 0 \simeq \idemp_2\otimes \idemp_1$ implies that $(\idemp_2,\idemp_1)$ are a pair of complementary idempotents by definition.
\end{proof}

\subsection{Reassociation and refinement of idempotent decompositions}
\label{subsec-idempotentreassoc}

When working with idempotents in an algebra, one is often concerned with finding indecomposable idempotents. A related question is whether one complete set of idempotents is courser or
finer than another. To discuss the analogous notion for idempotent decompositions of identity we need the following reassociation result, whose proof is straightforward from Proposition
\ref{prop:simplifications}.

Let $(I,\leq)$ be a finite poset.  An \emph{ideal} in $I$ is a subset $J\subset I$ such that $j\in J$ and $i\leq j$ implies $i\in J$.  A \emph{coideal} in $I$ is a subset $J\subset I$ such that $I\setminus J$ is an ideal.  Note that the set of ideals (resp.~coideals) is closed under unions and intersections.  A set $K\in I$ is convex if and only if $K = J_2\setminus J_1$ for some ideals $J_1\subset J_2\subset \PC$ (not necessarily unique).  

\begin{definition}\label{def:subquotientIdempotents}
Let $\one\simeq \Tot\{\idemp_i,d_{ij}\}$ be an $I$-indexed idempotent decomposition in $\AC$.  For every convex subset $K\subset I$, let $\idemp_K$ denote the contribution of $K$ to $\Tot\{\idemp_i,d_{ij}\}$ (Definition \ref{def:contrib}).  We refer to the objects $\idemp_K$ as \emph{subquotient idempotents}. Note that in the case of singletons we have $\idemp_i = \idemp_{\{i\}}$.
\end{definition}

The proof that $\idemp_K$ is idempotent is analogous to \eqref{eq:pipiredundant}.

\begin{proposition}\label{prop:resOfIdReass}
Let $\one\simeq \Tot\{\idemp_i,d_{ij}\}$ be an $I$-indexed idempotent decomposition in $\AC$.  The subquotient idempotents $\idemp_K$ satisfy 
\begin{enumerate}
\item $\idemp_{\emptyset}=0$ and $\idemp_I \simeq \one$.
\item $\idemp_J\otimes \idemp_K\simeq \idemp_{J\cap K}$.
\item $J\subset I$ is an ideal, then $\idemp_J$ has the structure of a unital idempotent with complementary counital idempotent $\idemp_{I\setminus J}$.
\item Suppose $I = \bigsqcup_{a \in \PC} J_a$ is a partition of $I$ into convex subsets, for some indexing set $\PC$.  Suppose that $\PC$ inherits a partial order from $I$.  Then reassociating gives an idempotent decomposition of identity $\one\simeq \Tot\{\idemp_{J_a}\}_{\PC}$.
\end{enumerate}
\end{proposition}
\begin{proof}
This is a straightforward application of reassociation and simultaneous simplifications.
\end{proof}

\begin{defn} Comparing two idempotent decompositions of identity, we say that $\one \simeq \Tot\{\idemp_i\}_I$ is \emph{finer} than $\one \simeq \Tot\{\idemp_{J_a}\}_\PC$ if the latter can be obtained from the former using a partition of $I$ into convex subsets, as in property (4) of Proposition \ref{prop:resOfIdReass}. \end{defn}

The following proposition states that commuting idempotent decompositions of identity have a mutual refinement.  Its proof is also straightforward.

\begin{proposition}\label{prop:tensoringResolutions}
Suppose $\one\simeq \Tot\{\idemp_i\}_I$ and $\one\simeq \Tot\{\otherIdemp_j\}_J$ are decompositions of $\one\in \AC$, such that $\idemp_i \ot \otherIdemp_j \cong \otherIdemp_j \ot \idemp_i$ for all $i \in I$ and $j \in J$. Then tensoring gives a decomposition of identity $\one\simeq \Tot\{\idemp_i\otimes \otherIdemp_j\}_{I\times J}$.  This tensor product decomposition of the identity is finer than both $\one\simeq \Tot\{\idemp_i\}_I$ and $\one\simeq \Tot\{\otherIdemp_j\}_J$. \end{proposition}

One key application of Proposition \ref{prop:tensoringResolutions} will be for simultaneous diagonalization. For example, suppose that $p_i$ (resp. $q_j$) are the idempotents which
project to the eigenspaces of a linear operator $f$ (resp. $g$). If $f$ and $g$ commute, then so do $p_i$ and $q_j$ for all $i, j$. The product $p_i q_j$ will project to a joint
eigenspace. Proposition \ref{prop:tensoringResolutions} will play a key role in categorifying this construction.

\begin{remark} Note that many of the tensor products $\idemp_i \ot \otherIdemp_j$ may be null-homotopic in the above proposition. Thus the minimal indexing set of this mutual refinement
is a subposet of $I \times J$, which can be more combinatorially interesting than $I \times J$ itself.\end{remark}

\begin{example}\label{ex-SBim3resolution} 
To give the flavor of reassociation and refinement, let us give an example from categorical Hecke theory (using results stated in \S\ref{subsec:HeckeConj}). We work in the homotopy
category of $\SBim_n$, the category of Soergel bimodules for $S_n$. There is an idempotent decomposition of identity $\one \cong \Tot\{\idemp_T\}$ indexed by standard tableaux $T$ with
$n$ boxes, with the dominance partial order. There is another idempotent decomposition of identity $\one \cong \Tot\{\cidemp_\l\}$ indexed by partitions of $n$. Because there is an
inclusion $\SBim_m \into \SBim_n$ for $m < n$, we also have idempotent decompositions of identity coming from partitions or tableaux of $m$.

The decomposition of identity by tableaux is obtained by refining that by partitions. More precisely, it is obtained by mutually refining the decomposition by partitions for all $m \le
n$. Conversely, the decomposition by partitions can be obtained by reassociating (i.e. coarsening) that by tableaux.
	
In case $n=3$, the decomposition of identity by tableaux will look like:
\[
\one \ \ \simeq \ \ \Tot\left(
\begin{minipage}{3in}
\begin{tikzpicture}
\node(a) at (0,0) {$\idemp_{\:\scriptscriptstyle \Yvcentermath1 \young(123)}$};
\node(b) at (3,1.5) {$\idemp_{\:\scriptscriptstyle \Yvcentermath1 \young(12,3)}$};
\node(c) at (3,-1.5) {$\idemp_{\:\scriptscriptstyle \Yvcentermath1 \young(13,2)}$};
\node(d) at (5.7,0) {};
\node at (6,-.2) {$\idemp_{\:\scriptscriptstyle \Yvcentermath1 \young(1,2,3)}$};
\path[->,>=stealth',shorten >=1pt,auto,node distance=1.8cm,
  thick]
(a) edge node {} (b)
(a) edge node {} (c)
(b) edge node {} (c)
(b) edge node {} (d)
(c) edge node {} (d);
\end{tikzpicture}
\end{minipage}
\right).
\]
Associating the middle two terms together defines a coarser decomposition of identity by partitions:
\[
\one \ \ \simeq \ \ \Tot\left(
\begin{minipage}{3in}
\begin{tikzpicture}
\node(a) at (0,0) {$\cidemp_{\:\scriptscriptstyle \Yvcentermath1 \yng(3)}$};
\node(b) at (3,0) {$\cidemp_{\:\scriptscriptstyle \Yvcentermath1 \yng(2,1)}$};
\node(d) at (5.7,0) {};
\node at (6,-.2) {$\cidemp_{\:\scriptscriptstyle \Yvcentermath1 \yng(1,1,1)}$};
\path[->,>=stealth',shorten >=1pt,auto,node distance=1.8cm,
  thick]
(a) edge node {} (b)
(b) edge node {} (d);
\end{tikzpicture}
\end{minipage}
\right).
\]
One also has
\[
\idemp_{\:\scriptscriptstyle \Yvcentermath1 \young(13,2)} \simeq \cidemp_{\:\scriptscriptstyle \Yvcentermath1 \yng(1)} \ot \cidemp_{\:\scriptscriptstyle \Yvcentermath1 \yng(1,1)} \ot \cidemp_{\:\scriptscriptstyle \Yvcentermath1 \yng(2,1)}\;,
\]
with similar homotopy equivalences for other tableaux. A tensor product of idempotents corresponding to partitions which does not lead to a standard tableau produces the zero complex:
\[0 \simeq \cidemp_{\:\scriptscriptstyle \Yvcentermath1 \yng(1)} \ot \cidemp_{\:\scriptscriptstyle \Yvcentermath1 \yng(1,1)} \ot \cidemp_{\:\scriptscriptstyle \Yvcentermath1 \yng(3)}\;. \]
\end{example}

\subsection{Categorical idempotents and uniqueness}
\label{subsec:catidemp}

The following is the fundamental theorem of unital idempotents. It actually holds in any triangulated monoidal category.

\begin{theorem}\label{thm:idemptOrder}
Let $(\PB_i,\e_i)$ ($i=1,2$) be unital idempotents in a triangulated monoidal category $\AC$.  The following are equivalent:
\begin{enumerate}
\item $\PB_1\otimes \PB_2\cong \PB_1$.
\item $\PB_2\otimes \PB_1\cong \PB_1$.
\item there exists a map $\nu:\PB_2\rightarrow \PB_1$ such that $\nu\circ \eta_2=\eta_1$.
\end{enumerate}
If either of these equivalent conditions is satisfied, then $\nu$ is unique, and it makes $\PB_1$ into a unital idempotent {relative to $\PB_2$} in the sense that $\nu\otimes \Id : \PB_2\otimes \PB_1\rightarrow \PB_1\otimes \PB_1$ and $\Id\otimes \nu : \PB_1\otimes \PB_2\rightarrow \PB_1\otimes \PB_1$ are isomorphisms. In particular, $\PB_1$ is a unital idempotent in the monoidal homotopy category $\PB_2 \AC \PB_2$. A similar result holds for counital idempotents.
\end{theorem}

For a proof see \cite{Hog17a}, \S 4.5. For a proof which does not utilize the triangulated (or additive) structure of $\AC$, see \cite{BoyDrin-idemp}.

The equivalent conditions of Theorem \ref{thm:idemptOrder} will be indicated by writing $\PB_1\leq \PB_2$. This is a well-defined partial order on the set of isomorphism classes of
unital idempotents in $\AC$. It is not difficult to show, using Theorem \ref{thm:idemptOrder}, that $\PB_1\leq \PB_2$ if and only if $\PB_1\AC\subset \PB_2\AC$ (equivalently
$\AC\PB_1\subset \AC\PB_2$). This partial order on idempotents is intimately related to the partial order on eigenmaps appearing in categorical diagonalizations (see
\S\ref{subsec:spectrum}).

If $\PB_1\leq \PB_2$ are unital idempotents, we refer to $\nu:\PB_2\rightarrow PB_1$ as the \emph{canonical map} between comparable unital idempotents. Of course, if $\AC$ is a homotopy category, then $\nu$ is unique only up to homotopy.
 
\begin{defn} If $\PB_1\leq \PB_2$, then we define an object $\DB(\PB_2,\PB_1)$ by choosing a distinguished triangle
\[
\DB(\PB_2,\PB_1) \rightarrow \PB_2 \rightarrow \PB_1 \rightarrow \DB(\PB_2,\PB_1)[1]
\]
in which the second map is the canonical map.\end{defn}

The object $\DB(\PB_2,\PB_1)$ should be thought of as the ``difference'' of $\PB_2$ and $\PB_1$.  it is easily seen that $\DB(\PB_2,\PB_1)$ is an idempotent which is fixed by $\PB_2$ and annihilated by $\PB_1$ on the left and right. 
\begin{definition}\label{def:categoricalIdemp}
A \emph{locally unital idempotent} or \emph{categorical idempotent} in $\AC$ is any object isomorphic to $\DB(\PB_2,\PB_1)$ for some unital idempotents $(\PB_i,\e_i)$.
\end{definition}
Locally unital idempotents were called \emph{locally closed idempotents} in \cite{BD14}.

\subsection{A convolution-free description}
\label{subsec:generalDecomp}

We have defined idempotent decompositions of identity for homotopy categories using twisted complexes and convolutions. Using reassociation we have already discussed how an idempotent decomposition of identity indexed by a poset $\PC$ gives a subquotient idempotent $\idemp_K$ for each convex subset $K$ of $\PC$, and a unital idempotent $\idemp_I$ for each ideal $I$ in $\PC$. The goal of this subsection is to provide an equivalent definition of an idempotent decomposition of identity, which begins with unital idempotents attached to ideals, and works in more generality (i.e. it does not require twisted complexes or convolutions to define). Now let $\AC$ be an arbitrary triangulated monoidal category.

\begin{defn} \label{defn:alternatedefnidempres} Let $(\PC,\leq)$ be a finite poset.  An (ideal-indexed) idempotent decomposition of $\one\in \AC$ is a collection of unital idempotents $\{\idemp_I,\eta_I\}$ in $\AC$, indexed by the ideals $I\subset \PC$, satisfying the conditions:
\begin{enumerate}
\item[(P1)] $\idemp_\emptyset = 0 $ and $\idemp_\PC \cong \one$.
\item[(P2)] $\idemp_{I\cap J} \cong \idemp_I\otimes \idemp_J$.
\item[(P3)] $\idemp_{I\cup J} = \idemp_I\vee \idemp_J$.
\end{enumerate}
The notation $\vee$ is to be interpreted in the usual manner for posets: $\idemp_I\vee \idemp_J$ is the minimum among all unital idempotents $\QB$ with $\idemp_1,\idemp_2\leq \QB$. This minimum exists by Proposition 5.1 in \cite{Hog17a}, which applies since $\idemp_I$ and $\idemp_J$ commute by (P2). \end{defn}

Condition (P2) implies that $\idemp_I\leq \idemp_J$ whenever  $I\subset J$. Condition (P3) implies that the following excision property holds: $\DB(\idemp_{I\cup J}, \idemp_J) \cong \DB(\idemp_I,\idemp_{I\cap J})$ (see Proposition 5.8 in \cite{Hog17a}).

\begin{lemma} Suppose one has an ideal-indexed idempotent decomposition of $\one$ as above. When $I \subset J$, the set $K = J \setminus I$ is convex, and we define $\idempotent_K$ to be $\DB(\idemp_J,\idemp_I)$. Then $\idempotent_K$ depends only on the convex set $K$, not on the choice of $I \subset J$ with $K = J \setminus I$, up to unique isomorphism. \end{lemma}

\begin{lemma} Suppose that $\AC$ is a monoidal homotopy category, and one has an idempotent decomposition of identity as in Definition \ref{def-resOfId}, indexed by a poset $\PC$. Then the subquotient idempotents $\idemp_K$ of Definition \ref{def:subquotientIdempotents} give rise to an ideal-indexed idempotent decomposition of identity. Conversely, given an ideal-indexed idempotent decomposition, the idempotents $\idemp_K$ corresponding to singletons $K$ give rise to an idempotent decomposition of identity. \end{lemma}
	
Because of this lemma, our two definitions of idempotent decompositions agree where both make sense.

\begin{remark}
We may regard $\PC$ as a topological space in which the open sets are the poset ideals $I\subset \PC$.  In this language, an idempotent decomposition of identity may be regarded as sheaf of idempotents on $\PC$.  The closed subsets $I\subset \PC$ correspond to unital idempotents $\PB_I$, the open subsets $\PC\setminus I$ correspond to counital idempotents $\PB_{\PC\setminus I}$, and more generally the locally closed subsets $K=I\setminus J$ correspond to locally unital idempotents $\PB_K = D(\PB_J,\PB_I)$.  The usual sheaf condition corresponds to the existence of a weak pullback square relating the idempotents appearing in (P2) and (P3).  See \S 5.5 in \cite{Hog17a} for more details.
\end{remark}

\begin{remark} \label{rmk:semiorthogonal}
Let $\MC$ be a triangulated category on which $\AC$ acts by exact functors.  Let $\{\idemp_I,\e_I\}$ be an idempotent decomposition of identity.  For each poset ideal $I\subset \PC$ , denote the image of $\PB_I$ in $\MC$ by $\MC_I$.  If $I\subset J$, then $\MC_I\subset \MC_J$.  In fact the inclusion $\MC_I\rightarrow \MC_J$ has a right adjoint given by tensoring with $\PB_I$.  Thus, idempotent decompositions of identity are a natural extension of the theory of semi-orthogonal decompositions \cite{BonOrl-semi}.
\end{remark}

%
%

\section{Pre-diagonalizability and diagonalizability}
\label{sec:prediag}

\subsection{Scalar objects}
\label{subsec:scalars}

Recall the following definition.

\begin{defn} Let $\CC$ be a monoidal category. The \emph{Drinfeld center} $\ZC(\CC)$ is the category whose objects are pairs $(Z,\phi)$, where $Z$ is an object of $\CC$ and $\phi$, the \emph{braiding morphism}, is a natural isomorphism of endofunctors from $Z \ot (-)$ to $(-) \ot Z$.  Morphisms in the Drinfeld center are morphisms $Z \to Z'$ in $\CC$ whose induced maps $Z \ot (-) \to Z' \ot (-)$ and $(-) \ot Z \to (-) \ot Z'$ commute with the braiding morphisms. \end{defn}

Thus the Drinfeld center has a forgetful map to $\CC$. The Drinfeld center is naturally braided, though it is typically not symmetric.

\begin{defn} \label{def:scalarcat} Let $\CC$ be an monoidal category which is either additive or triangulated. A \emph{category of scalars in $\CC$} is an additive symmetric monoidal category $\KS$, equipped with a braided monoidal functor $\rho \co \KS \to \ZC(\CC)$ from $\KS$ to the Drinfeld center of $\CC$. \end{defn}

\begin{remark}\label{rmk:braiding}
Let $\tau_{\l,C}$ denote the braiding morphism $\l\otimes C\rightarrow C\otimes \l$ associated to a scalar object $\l$.  The assumption that $\KS$ is symmetric implies that if $\l,\mu$ are scalar objects then $\tau_{\l,\mu}$ and $\tau_{\mu,\l}\inv$ equal (or homotopic in the case of homotopy categories of complexes) as isomorphisms $\l\otimes \mu\rightarrow \mu\otimes \l$.  This assumption will not play a role for the most part, except in \S \ref{subsec:structure}.
\end{remark}

We usually abuse notation and write $\rho(\lambda)$ simply by $\lambda$. For each $\lambda\in \KS$, we have the functor $B\mapsto \lambda\otimes B$, which is naturally isomorphic to $B\mapsto B\otimes \lambda$.

\begin{definition}\label{def:scalar} Objects of $\KS$ will be called \emph{scalar objects}, and the functors $\lambda\otimes(-)$ with $\lambda\in \KS$ will be called \emph{scalar functors}.  A \emph{scalar object of $\CC$} is any object isomorphic to $\l\otimes \one$, for $\l\in \KS$. We sometimes abuse notation and write $\l\otimes \one$ simply by $\l$, and omit the tensor product symbol when acting with $\l$, e.g. $B \mapsto \l B$.

We will be particularly interested in \emph{invertible scalar objects/functors} $\l$, for which there is a scalar object $\l^{-1}$ such that $\l \ot \l^{-1} \cong \one \in \KS$.
\end{definition}

\begin{remark} One reason we view $\KS$ as a separate additive category, rather than thinking of it as a subcategory of the (possibly triangulated) category $\ZC(\CC)$, is that it leads to a larger Grothendieck group. \end{remark} 

\begin{ex}\label{ex:gradedScalars}
When $R$ is a ring, and $\AC=\KC^b(R\bimod)$, it is common to let $\KS$ denote the category of graded abelian groups which are finitely generated and free as $\Z$-modules.  The action of $\KS$ on $\AC$ is determined by the following:
\begin{enumerate}
\item If $\Z_i$ denotes a copy of $\Z$ sitting in degree $i$, then $\Z_i\ot B = \Z_i \ot_\Z B \cong B[-i]$.
\item $(\l\oplus\mu)\ot B\cong \l\ot B \oplus \mu\ot B$.
\end{enumerate}
In this setting, the invertible scalar objects in $\AC$ are of the form $R[i]$, and they categorify $(-1)^i$. The Grothendieck group of the category $\KS$ itself is larger: $[\KS] = \Z[T,T\inv]$, and the action of $[\KS]$ on $[\AC]$ factors through $T = -1$.
\end{ex}

\begin{ex}\label{ex:bigradedScalars}
If $R$ is a graded ring, and $\AC=\KC^b(R\gbimod)$, it is common to let $\KS$ denote the category of bigraded abelian groups which are finitely generated and free as ungraded $\Z$-modules.  The action of $\KS$ on $\AC$ is determined by the following:
\begin{enumerate}
\item If $\Z_{i,j}$ denotes a copy of $\Z$ sitting in bidegree $i,j$, then $\Z_{i,j}\otimes B\cong B[-i](-j)$. 
\item $(\l\oplus\mu)\otimes_\Z B\cong (\l\otimes B) \oplus (\mu\otimes B)$.
\end{enumerate}
In this setting, the invertible scalar objects in $\AC$ are of the form $R[i](j)$, and they categorify powers of $v$ with signs. The Grothendieck group of $\KS$ itself is $\Z[v,v\inv,T,T\inv]$.
\end{ex}


In the above examples, the invertible scalar functors are just called \emph{shifts}. We may also call them \emph{monomial scalar objects}, because they categorify monomials.

\begin{ex}\label{ex:sheavesScalars}
When $\AC$ is a category of sheaves on an algebraic variety $X$, it is typical to let $\KS$ denote the category of vector bundles on $X$, which acts on $\AC$ by the usual tensor product.  In this setting, the invertible scalar objects are the line bundles.\end{ex}

\begin{definition}\label{def:gradedHoms}
Let $\CC$ be a monoidal category with category of scalars $\KS$.  Let $\KS^\times\subset \KS$ denote the full subcategory consisting of the invertible objects.  Let $A,B\in \CC$ be given.  For any $\l\in \KS^\times$, \emph{a morphism $A\rightarrow B$ of degree $\l$} will mean a morphism $f:\l\otimes A\rightarrow B$.  This can be assembled into the \emph{enriched hom space}
\[
\Homg_\CC(A,B) =\bigoplus_{\l\in \KS^\times} \Hom_{\CC}(\l\otimes A,B).
\]
The composition of a degree $\mu$ map $g:\mu  A\rightarrow B$ and a degree $\l$ map $f:\l  B\rightarrow C$ is by definition $f\circ(\l g)$.  In this way one can form the enriched category $\CC'$ with the same objects as $\CC$, but with morphism spaces $\Homg_\CC(A,B)$.  The tensor product of maps $f:\l A\rightarrow B$ and $g:\mu C\rightarrow D$ is by definition the composition
\[
\l \mu (A\otimes C)  \cong (\l A) \otimes (\mu C) \buildrel f\otimes g\over \longrightarrow B\otimes D.
\]
This makes $\CC'$ into a sort of $\KS^\times$-graded monoidal category, but we will not spell this out explicitly here.
\end{definition}

\subsection{Commutativity and eigencones}
\label{subsec:commutativitystuff}

In \S\ref{subsec:intro_catDiag2} we discussed one of the major differences between linear algebra and categorical linear algebra: operators which are guaranteed to commute in linear algebra need not commute categorically. For example, any two polynomials in an operator $f$ commute with each other, while two eigencones need not commute.  In this section we discuss various obstruction theoretic conditions which have important consequences for eigenmaps and the commutativity of their eigencones, as well as other technical statements. 



\begin{definition}\label{def:aCommutesWithF}
Let $\l$ be a scalar object and $\a:\l\rightarrow F$ a map.  We say $\a$ \emph{commutes with $F$} if there is an element $h_\a\in \Hom^{-1}(\l\otimes  F, F\otimes F)$ such that
\[
d\circ h + h\circ d = \a  F - F\a.
\]
\end{definition}
Note that this definition only makes sense when $\AC$ is a monoidal homotopy category.  In this equation, notation such as $\a F$ is shorthand for $\a\otimes \Id_F$, and we have omitted an occurence of the isomorphism $\l F\rightarrow F \l$ from the second term on right-hand side.  We make similar such abbreviations throughout this section.  Compare with Definition \ref{def:gCommutesWithF}.

\begin{proposition}\label{prop:commutavityProps}
Suppose $\a:\l\rightarrow F$ and $\b:\mu\rightarrow F$ commute with $F$, and let $h_\a,h_\b$ denote the corresponding homotopies.  Then
\begin{enumerate}\setlength{\itemsep}{2pt}
\item $\Cone(\a)\otimes F\simeq F\otimes \Cone(\a)$.
\item the element $z_{\b,\a}:= h_{\b}\circ \a + h_\a\circ \b$ is a cycle in $\Hom^{-1}(\mu\l,FF)$.  If $z_{\b,\a}\simeq 0$, then
\[
\Cone(\a)\otimes \Cone(\b)\simeq \Cone(\b)\otimes \Cone(\a).
\]
\item the element  $w_\a=h_\a\circ \a$ is a cycle in $\Hom^{-1}(\l\l,FF)$.  If $w_\a\simeq 0$, then
\[
\Cone(\a)^{\otimes 2} \simeq (F\oplus \l[1])\otimes \Cone(\a).
\]
\end{enumerate}
\end{proposition}
\begin{proof}
Statements (1) and (2) are proven in \S \ref{subsec:commutativity}, and statement (3) is proven in \S \ref{subsec:selfObstruction}.
\end{proof}

\begin{definition}\label{def:roughObstructionFree}
If $\SC$ is a collection of morphisms $\a_i:\l_i\rightarrow F$, then we say $\SC$ is a \emph{weakly commuting family} if each $\a_i$ commutes with $F$.  We say $\SC$ is a \emph{strongly commuting family if} the homotopies $h_i\in \Hom^{-1}(\l_i F,FF)$ can all be chosen so that $z_{\a,\b}\simeq 0$.  If, in addition, the $h_i$ can be chosen so that the obstructions $w_\a$ vanish, then we say that $\SC$ is \emph{obstruction-free}.  
\end{definition}

Various consequences of these definitions are explored in Appendix \ref{sec:appendix}. For now we only mention that the existence of $h_\a$ is in practice relatively easy to check, but
the conditions that the cycles $w_\a$ and $z_{\a,\b}$ be null-homotopic represents a serious technical difficulty to be overcome. In all examples we know of, the only way to prove that
the obstructions $w_\a$ and $z_{\a,\b}$ vanish is to argue that they must be zero for degree reasons.

These commutativity hypotheses are also used to prove that the periodicity maps for eigenprojections are related to the eigenmaps as discussed in \S\ref{subsec:intro_period}, see Proposition \ref{prop:eigenaction}.

\subsection{Pre-diagonalizability}
\label{subsec:donotcommute}

In this section we discuss prediagonalizability of $F\in \AC$, which is a necessary condition for diagonalizability.  The definition is quite technical, but simplifies tremendously under the commutativity conditions discussed in the previous section. 



\begin{defn}\label{def:prediag}
We say $F\in \AC$ is \emph{categorically prediagonalizable} if it is equipped with a finite set of maps $\{\a_i:\l_i\rightarrow F\}_{i\in I}$ (where $\l_i$ are scalar objects), such that
\begin{enumerate}
\item[(PD1)] $\bigotimes_{i\in I} \Cone(\a_i)\simeq 0$ for all orderings of the factors.
\item[(PD2)] the indexing set $I$ is minimal with respect to (PD1).
\item[(PD3)] for any surjective map of sets $\pi \co \{1, \ldots, s\} \to I$ we have
\begin{equation}\label{eq:bigTensorVanishes}
\bigotimes_{a = 1}^s \Cone(\a_{\pi(a)})\simeq 0.
\end{equation}
\end{enumerate}
In this case the collection $\{\a_i\}_{i \in I}$ is called a \emph{saturated collection of eigenmaps} or a \emph{prespectrum} for $F$. \end{defn}

Clearly (PD3) implies (PD1).  If there is an isomorphism $\Cone(\a_i) \ot \Cone(\a_j) \cong \Cone(\a_j) \ot \Cone(\a_i)$ for all $i, j \in I$, then (PD1) implies (PD3), and can be checked on a single ordering of the tensor factors. Thus Definition \ref{def:prediag} and Definition \ref{def:prediagintro} are equivalent in this circumstance.  More generally we have the following.

\begin{lemma} \label{lem:pd1vspd3} If $\{\a_i\}_{i\in I}$ is a weakly (resp.~strongly) commuting family then (PD3) is equivalent to the statement that
\begin{equation} \label{eq:lessBigTensorVanishes}
\bigotimes_{i \in I} \Cone(\a_i) \cong 0
\end{equation}
for all possible orderings (resp.~at least one ordering) of the tensor factors. \end{lemma}

\begin{proof}  If $\{\a_i\}$ is a strongly commuting family then the cones $\{\Cone(\a_i) \}$ commute up to homotopy, and the lemma is obvious.  So assume that $\{\a_i\}$ is a weakly commuting family.

One direction is clear. Now suppose \eqref{eq:lessBigTensorVanishes}. Let $\pi$ be surjective and consider \[Z = \bigotimes_{a = 1}^s \Cone(\a_{\pi(a)}).\] One can reassociate this
tensor product, replacing any given cone with two terms, one copy of $\l_i$ and one copy of $F$. Both of these commute past the remaining cones in the tensor product up to isomorphism,
$\l_i$ because it lives in the Drinfeld center, and $F$ by Proposition \ref{prop:commutavityProps}. Consequently, by Proposition \ref{prop:simplifications}, if some ordered subset of
the cones $\Cone(\a_{\pi(a)})$ will tensor to zero, then the entire complex $Z$ is nulhomotopic. \end{proof}

\begin{remark}In \S\ref{sec:extendedintro} we simplified the exposition by assuming that eigencones commute, but in fact everything holds true verbatim in the general case so long as one uses Definition \ref{def:prediag} for prediagonalizability.
\end{remark}



\subsection{Basics of categorical diagonalization}
\label{subsec-diagonalizable}

We recall the definition from the introduction.

\begin{definition}\label{def:catdiag}
Let $F\in \AC$ be an object of a homotopy monoidal category.  Let $I$ be a finite poset, and suppose we are given scalar objects $\l_i\in \AC$, maps $\a_i:\l_i\rightarrow F$, and nonzero objects $\PB_i\in \AC$, indexed by $i\in I$.  We say that $\{\PB_i,\a_i\}_{i\in I}$ is a \emph{diagonalization of $F$} if
\begin{itemize}\setlength{\itemsep}{2pt}
\item there is an idempotent decomposition of identity $\one \simeq\Big(\bigoplus_{i\in I}\PB_i, d\Big)$.
\item $\Cone(\a_i)\otimes \PB_i\simeq 0 \simeq \PB_i\otimes \Cone(\a_i)$.
\end{itemize} 
We say that $F$ is \emph{diagonalizable} if there exists a diagonalization of $F$. We say that the diagonalization is \emph{tight} if, whenever $\AC$ acts on a triangulated category $\VC$ and $M$ is an object of $\VC$, one has $\Cone(\a_i) \ot M \simeq 0$ if and only if $\PB_i \ot M \simeq M$.
\end{definition}

The next proposition states that if $F$ is tightly diagonalized then it is prediagonalized.

\begin{prop}\label{prop:bigtensorIsZero}
If $\{\PB_i,\a_i\}_{i\in I}$ is a diagonalization of $F$, then (PD3) holds for $\{\a_i\}_{i\in I}$. If the diagonalization is tight, then (PD2) holds as well, so that $\{\a_i\}_{i\in I}$ is a saturated collection of eigenmaps for $F$, and $F$ is prediagonalizable. \end{prop}

If the diagonalization is not tight, then nothing prevents the $\a_i$ from all being equal, as in \S\ref{subsec:semisimple}. Consequently, (PD2) can fail.

The proof uses two lemmas.

\begin{lemma} Suppose that $\{\PB_i,\a_i\}_{i\in I}$ is a diagonalization of $F$.  If $Z\in \AC$ satisfies $Z\otimes \idemp_i\simeq 0$ for all $i\in I$, then $Z\simeq 0$.
\end{lemma}

\begin{proof}
Tensoring $Z$ with $\one\simeq \Tot\{\idemp_i\}$ yields $Z\simeq \Tot\{Z\otimes \idemp_i\}$.  Each term is contractible by hypothesis, hence the total complex is contractible by Proposition \ref{prop:simplifications}.
\end{proof}

\begin{lemma} \label{lem:PConeCommute}
If $\{\PB_i,\a_i\}_{i\in I}$ is a diagonalization of $F$, then $\idemp_i \ot \Cone(\a_j) \cong \Cone(\a_j) \ot \idemp_i$, for all $i$ and $j$. \end{lemma}

\begin{proof} The full subcategory $\ZC_{\idemp_i}\subset \AC$ consisting of complexes $C$ which commute with a categorical idempotent $\idemp_i$ is triangulated (see \S 4.4 in \cite{Hog17a}).  Now, $\l_j$ commutes with $\idemp_i$ since $\l_j$ commutes with everything, and $F$ commutes with $\idemp_i$ by definition. Thus $\Cone(\a_j)$ commutes with $\idemp_i$ for all $i,j\in I$. \end{proof}

\begin{proof}[Proof of Proposition \ref{prop:bigtensorIsZero}]
Let $Z$ be the tensor product in \eqref{eq:bigTensorVanishes}. By definition, we have $\Cone(\a_i) \otimes \idemp_i \simeq 0$. Thanks to Lemma \ref{lem:PConeCommute}, we deduce that $Z \ot \idemp_i \simeq 0$ for all $i$. By the lemma preceding that, we conclude that $Z\simeq 0$. Thus (PD3) holds.

Suppose that the diagonalization is tight. If (PD2) fails and there are $r+1$ projectors, then some tensor product of $r$ distinct cones (in some order) is zero. The partial order on
$I$ will not be relevant for this part of the proof, so let us rename our eigenmaps and suppose that \[ Y_1 = \Cone(\a_1) \ot \cdots \ot \Cone(\a_r) \simeq 0,\] where $\Cone(\a_0)$ is
absent. We will deduce that $\idemp_0 = 0$, a contradiction.

Let $Y_k = \Cone(\a_k) \ot \cdots \ot \Cone(\a_r)$ for $1 \le k \le r$. Because $\Cone(\a_1) \ot Y_2 \simeq 0$, we have $\idemp_1 Y_2 \simeq Y_2$.   This implies that $\idemp_i \ot Y_2 \simeq \idemp_i\ot \idemp_1 \ot  Y_2\simeq 0$ for all $i \neq 1$.


Since $\idemp_0 \ot Y_2 \simeq 0$ then (using Lemma \ref{lem:PConeCommute} again) we see that $\Cone(\a_2)$ kills $\idemp_0 \ot Y_3$. Thus $\idemp_2$ fixes $\idemp_0 \ot Y_3$. But $\idemp_2 \ot \idemp_0 \simeq 0$, meaning that $\idemp_0 \ot Y_3 \simeq \idemp_2 \ot \idemp_0 \ot Y_3 \simeq 0$.

Since $\idemp_0 \ot Y_3 \simeq 0$ then $\Cone(\a_3)$ kills $\idemp_0 \ot Y_4$. Thus $\idemp_3$ fixes $\idemp_0 \ot Y_4$, so that $\idemp_0 \ot Y_4 \simeq 0$.

Repeating this argument, at the last step we obtain that $\idemp_0 \simeq \idemp_r \ot \idemp_0 \simeq 0$. \end{proof}

\subsection{Categorified quasi-idempotents}
\label{subsec:quasiIdemp1}

We do not know of a complete characterization of which pre-diagonalizable objects are diagonalizable in general.  However, for prediagonalizable objects there is a closely family of objects $\KB_i$ which, while not idempotent, are interesting in their own right.

\begin{definition}\label{def:Ki}
If $\{\a_i:\l_i\rightarrow F\}_{i\in I}$ is a saturated collection of eigenmaps, then set $\KB_i=\bigotimes_{j\neq i}\Cone(\a_j)$ (choose an ordering).  
\end{definition}

\begin{proposition} \label{prop:KBiprops}
Assume that $\{\a_i:\l_i\rightarrow F\}$ is a saturated collection of eigenmaps which is obstruction-free.  Then
\begin{enumerate}\setlength{\itemsep}{2pt}
\item $\Cone(\a_i)\otimes \KB_i\simeq 0 \simeq \KB_i\otimes \Cone(\a_i)$.
\item $\KB_i$ does not depend on the ordering of the factors $\Cone(\a_j)$ up to equivalence.
\item $\KB_i^{\otimes 2} \simeq \bigotimes_{j\neq i}(\l_i\oplus \l_j[1]) \otimes \KB_i$.
\end{enumerate}
\end{proposition}
\begin{proof}
Statements (1) and (2) are clear, while statement (3) is an immediate consequence of Proposition \ref{prop:commutavityProps}.
\end{proof}

Thus, the $\KB_i$ are $\a_i$-eigenobjects, and they are quasi-idempotent up to homotopy.  Recall that in linear algebra an element $k\in A$ of a $\KM$-algebra is quasi-idempotent if there is a scalar $c$ such that $k^2= c k$.  If $c$ is invertible, then $p:=c\inv k$ is an idempotent.  Thus, the problem of categorical diagonalization comes down to: if $\{\a_i:\l_i\rightarrow F\}$ is an obstruction-free saturated collection of eigenmaps, how can we construct an idempotent $\idemp_i$ from $\KB_i$?  What is the categorical analogue of division by the scalar $\bigotimes_{j\neq i}(\l_i\oplus \l_j[1])$?

To explore this question, we assume that $F$ is diagonalizable with eigenmaps $\a_i:\l_i\rightarrow F$ and eigenprojections $\idemp_i$, and examine the relationship between $\KB_i$ and $\idemp_i$.

\begin{definition} \label{def:rhoi}
For each scalar object $\mu$ and each map $\b:\mu\rightarrow F$, let $\rho(\b):\mu\idemp_i\rightarrow \l_i\idemp_i$ denote the composition
\[
\begin{diagram}
\mu\otimes \idemp_i &\rTo^{\b\otimes \Id_{\idemp_i}} & F\otimes \idemp_i & \rTo^{\simeq } & \l_i\otimes \idemp_i,
\end{diagram}
\]
where the second map is the inverse of $\a_i \ot \idemp_i$.
\end{definition}

Observe that $\rho_i(\a_i)\simeq \Id_{\l_i\idemp_i}$.

\begin{proposition}\label{prop:KasKoszulComplex}
If $F$ is tightly diagonalized with eigenmaps $\{\a_i:\l_i\rightarrow F\}_{i\in I}$ and eigenprojections $\idemp_i$, then
\begin{equation}\label{eq:KasKoszulComplex}
\bigotimes_{j\neq i}\Cone\Big(\begin{diagram}[small]\l_j\idemp_i & \rTo^{\rho_i(\a_j)} & \l_i\idemp_i\end{diagram}\Big)  \ \ \simeq \ \ \KB_i
\end{equation}
\end{proposition}

\begin{proof}
Recall that  $\KB_i$ is an $\a_i$-eigenobject on the left and right, hence (by tightness) $\PB_i\otimes \KB_i\simeq \KB_i\otimes \PB_i$.  On the other hand,
\[
\Cone(\rho_i(\a_j)) \simeq \Cone(\a_j\otimes \Id_{\idemp_i}) \cong \Cone(\a_j)\otimes \idemp_i \cong \idemp_i \ot \Cone(\a_j).
\]
The first isomorphism used that two ``isomorphic'' maps have isomorphic cones; the second used the compatibility of cones and tensor products. Thus,
\[
\KB_i\simeq \KB_i\otimes \idemp_i =\bigotimes_{j\neq i} \Cone(\a_j)\otimes \PB_i\simeq \bigotimes_{j\neq i} \Cone(\rho_i(\a_j))\otimes \idemp_i,
\]
as claimed. Here, we used Lemma \ref{lem:PConeCommute} several times, to duplicate $\idemp_i$ and bring it past the various eigencones $\Cone(\a_j)$.
\end{proof}

This proposition is the categorical analogue of the identity $\prod_{j\neq i}(\k_i-\k_j)p_i = \prod_{j\neq i}(f-\k_j)$.  It states that $\KB_i$ is homotopic to a sort of Koszul complex associated to $\idemp_i$ together with the endomorphisms $\rho_i(\a_j)$.  Based on this relationship, one wishes to construct $\PB_i$ by a categorical analogue of dividing by $\prod_{j\neq i}(\k_i-\k_j)$.  Making sense of this is a very delicate issue, one which we do not know how to solve.  However, if each $\l_i$ is invertible then this can be accomplished essentially by mimicking the classical Koszul duality.  This is the subject of the next chapter.

\section{Interpolation complexes for invertible eigenvalues}
\label{sec:interpolation}

In this chapter we state and prove our main diagonalization theorem, which gives a sufficient condition for a complex to be categorically diagonalizable. We restrict ourselves to the
case of invertible complexes. Also, because we use convolutions and reassociations in our proofs, we restrict to homotopy monoidal categories (see Definition
\ref{def:homotopyMonoidalCat}). Although the main constructions work in more generality, our proofs do not.

We have endeavored to prove our main theorem, Theorem \ref{thm-firstEigenthm}, without assuming any of the simplifying hypotheses discussed in \S\ref{subsec:commutativitystuff}. The
consequence is that several proofs are complicated by the lack of commutativity of eigencones. The first-time reader should try to ignore such technical issues.

\subsection{Completed scalars and smallness}
\label{subsec-completed}

We wish to categorify the Lagrange interpolation formula $p_i=\prod_{j\neq i}\frac{f-\k_j}{\k_i-\k_j}$ for the eigenprojections of a diagonalizable linear operator $f$.  Roughly speaking, we will categorify denominators of the form $\frac{1}{1-\k}$ by expanding into a geometric series and working with infinite complexes.  

We first extend the category of scalar objects to include certain ``formal series'' of scalar objects.

\begin{definition}\label{def:completedScalars}
Let $\CC$ be an (additive or triangulated) monoidal category with category of scalar objects $\KS\subset \ZC(\CC)$.  Let $\hat{\KS}\subset \CC$ denote the full subcategory of objects $\Gamma$ such that
\begin{enumerate}
\item There exist scalar objects $\l_0,\l_1,\ldots,\in \KS$ such that $\Gamma\cong \bigoplus_{k\geq 0}\l_k\cong \prod_{k\geq 0}\l_k$, and
\item $\bigoplus_{k\geq 0} \l_k\otimes B\cong \prod_{k\geq 0}\l_k\otimes B$ and both exist, for all $B\in \CC$.
\end{enumerate}
Then we call $\hat{\KS}$ the \emph{completed category of scalars}.
\end{definition}

\begin{example}\label{ex:homologicalCompletion}
Let $\BC$ be an additive monoidal category.  Let $\AC=\KC^-(\BC)$.  For the category of scalars $\KS\subset \ZC(\AC)$, let us choose the finite direct sums of shifted copies of $\one$ (as in Example \ref{ex:bigradedScalars}).  In this context, $\hat{\KS}\subset \AC$ contains all complexes with zero differential whose terms are scalar objects in $\BC$.  That is to say, $\l_0\oplus \l_1[1]\oplus \l_2[2]\oplus \cdots $ is in $\hat{\KS}$ for any sequence of scalar objects $\l_i \in \BC$, as is any shift thereof.
\end{example}

\begin{definition}\label{def:polyObjects}
Let $\CC$ be an (additive or triangulated) monoidal category with scalar objects $\KS\subset \CC$.  Let $u$ be a formal indeterminate with formal degree $\deg(u)=\l\in \KS$.  Then $\one[u]$ will denote the object\footnote{The reader should think of $\one[u]$ as being like a polynomial ring, rather than viewing the brackets as indicating a shift.}
\[
\one[u] = \one\oplus (\l\otimes \one) \oplus (\l^2\otimes \one) \oplus \cdots .
\]
This infinite direct sum may or may not be well-defined in $\CC$.  We say that $\l$ is \emph{small} if $\one[u]$ is defined, and is an object of $\hat\KS$.  In this case, we set $B[u]:=\one[u]\otimes B$ for all $B\in \CC$.
\end{definition}



The most common instances of smallness come from the following example (compare with the notion of homological local finiteness from Example \ref{ex:homLocalFinite}).

\begin{ex}[Homological smallness]\label{ex:smallnessOfShifts}
Suppose $\BC$ is an additive monoidal category and $\AC=\KC^-(\BC)$.  Then $\l:=\one[i]$ is small in $\AC$ if and only if $i>0$. 

Similarly, if $i<0$, then $\one[i]$ is small as an object of $\KC^+(\BC)$.
\end{ex}

\subsection{Categorified interpolation complexes}
\label{subsec-cateInterpoly}

We recall a number of definitions from \S\ref{sec:extendedintro}. Fix an object $F$ in a monoidal homotopy category $\AC$, with a category of scalars $\KS$.
 
\begin{definition}\label{def:Cij}
Suppose $\l,\mu\in \KS$ are scalar objects such that $\l\mu\inv$ is small, and $\a,\b\in \Homg(\one,F)$ are maps of degree $\l$ and $\mu$, respectively.  Then define the following complexes $\Cij{\b}{\a}(F)$ and $\Cij{\a}{\b}(F)$ by:
\begin{equation}\label{eq:Cab}
\Cij{\b}{\a}(F) \ \  := \ \  \left(
\begin{minipage}{1.7in}
\begin{tikzpicture}
\node (aa) at (0,8){$\one$ };
\node (ab) at (0,6){ $\displaystyle\frac{\l}{\mu} $};
\node (ac) at (0,4){ $\displaystyle\frac{\l^2}{\mu^2}$ };
\node (ad) at (0,2){ $\displaystyle\frac{\l^3}{\mu^3}$ };
\node (ae) at (0,0){$\cdots$};
\node (ba) at (3,8){};
\node (bb) at (3,6){$\displaystyle \frac{1}{\mu} F[-1]$};
\node (bc) at (3,4){$\displaystyle \frac{\l}{\mu^2} F[-1]$};
\node (bd) at (3,2){$\displaystyle \frac{\l^2}{\mu^3} F[-1]$};
\node (be) at (3,0){$\cdots$};
\path[->,>=stealth',shorten >=1pt,auto,node distance=1.8cm,font=\small]
(aa) edge node[auto] {$\b$} (bb)
(ab) edge node {$\a$} (bb)
(ab) edge node {$-\b$} (bc)
(ac) edge node {$\a$} (bc)
(ac) edge node {$-\b$} (bd)
(ad) edge node {$\a$} (bd)
(ad) edge node {$-\b$} (be);
\end{tikzpicture}
\end{minipage}
 \right)
 \ \ \ \ \ \ \ \ \ \ \ \ \ \ 
\Cij{\a}{\b}(F) \ \  := \ \  \left(
\begin{minipage}{1.7in}
\begin{tikzpicture}
\node (ab) at (0,6){ $\displaystyle\frac{\l}{\mu}[1] $};
\node (ac) at (0,4){ $\displaystyle\frac{\l^2}{\mu^2}[1]$ };
\node (ad) at (0,2){ $\displaystyle\frac{\l^3}{\mu^3}[1]$ };
\node (ae) at (0,0){$\cdots$};
\node (bb) at (3,6){$\displaystyle \frac{1}{\mu} F$};
\node (bc) at (3,4){$\displaystyle \frac{\l}{\mu^2} F$};
\node (bd) at (3,2){$\displaystyle \frac{\l^2}{\mu^3} F$};
\node (be) at (3,0){$\cdots$};
\path[->,>=stealth',shorten >=1pt,auto,node distance=1.8cm,font=\small]
(ab) edge node {$\a$} (bb)
(ab) edge node {$-\b$} (bc)
(ac) edge node {$\a$} (bc)
(ac) edge node {$-\b$} (bd)
(ad) edge node {$\a$} (bd)
(ad) edge node {$-\b$} (be);
\end{tikzpicture}
\end{minipage}
 \right)
\end{equation}
Because of the smallness condition, these complexes are well-defined objects in $\AC$.
\end{definition}

The following is clear from the definitions.
\begin{lemma} \label{lem:CijTriangle}
There exists a distinguished triangle of the form
\[
\Cij{\b}{\a}(F)\rightarrow \one\rightarrow \Cij{\a}{\b}(F)\rightarrow \Cij{\b}{\a}(F)[1],
\]
so that $\one\simeq (\Cij{\a}{\b}(F)\buildrel [1]\over \rightarrow \Cij{\b}{\a}(F))$.\qed
\end{lemma}

The periodicity in the construction of $\Cij{\a}{\b}$ is represented by a distinguished endomorphism, which we refer to as the \emph{periodicity map}.

\begin{definition}\label{def:periodicityMapCab}
Let $\mu$ and $\l$ be invertible scalar objects.  If $\mu\l\inv$ is small, then the following diagram defines a chain map $f_{\a,\b}:\frac{\l}{\mu} \Cij{\a}{\b}(F)\rightarrow \Cij{\a}{\b}(F)$:
\[
\begin{tikzpicture}
\node(aa) at (0,2.5) {$0$};
\node(ba) at (2.5,2.5) {$0$};
\node(ca) at (5,2.5) {$\displaystyle \frac{\l}{\mu^2} F$};
\node(da) at (7.5,2.5) {$\displaystyle \frac{\l^2}{\mu^2} $};
\node(ea) at (10,2.5) {$\displaystyle \frac{\l^2}{\mu^3} F$};
\node(fa) at (12.5,2.5) {$\cdots$};
\node(ab) at (0,0) {$\displaystyle \frac{1}{\mu} F$};
\node(bb) at (2.5,0) {$\displaystyle \frac{\l}{\mu} $};
\node(cb) at (5,0) {$\displaystyle \frac{\l}{\mu^2} F$};
\node(db) at (7.5,0) {$\displaystyle \frac{\l^2}{\mu^2} $};
\node(eb) at (10,0) {$\displaystyle \frac{\l^2}{\mu^3} F$};
\node(fb) at (12.5,0) {$\cdots$};
\tikzstyle{every node}=[font=\small]
\path[->,>=stealth',shorten >=1pt,auto,node distance=1.8cm]
(ba) edge node[above] {} (aa)
(ba) edge node[above] {} (ca)
(da) edge node[above] {$ \frac{\l}{\mu^2}\a$} (ca)
(da) edge node[above] {$ -\frac{\l^2}{\mu^3}\b$} (ea)
(fa) edge node[above] {$ \frac{\l^2}{\mu^3}\a$} (ea)
(bb) edge node[above] {$ \frac{1}{\mu} \a$} (ab)
(bb) edge node[above] {$ -\frac{\l}{\mu^2} \b $} (cb)
(db) edge node[above] {$ \frac{\l}{\mu^2}\a$} (cb)
(db) edge node[above] {$ -\frac{\l^2}{\mu^3}\b$} (eb)
(fb) edge node[above] {$ \frac{\l^2}{\mu^3}\a$} (eb)
(aa) edge node[left] {} (ab)
(ba) edge node[left] {} (bb)
(ca) edge node[left] {$\Id$} (cb)
(da) edge node[left] {$\Id$} (db)
(ea) edge node[left] {$\Id$} (eb);
\end{tikzpicture}.
\]
Similarly, if $\mu\l\inv$ is small, the following diagram defines a chain map $f_{\a,\b}: \frac{\l}{\mu}\Cij{\a}{\b}(F)\rightarrow \Cij{\a}{\b}(F)$:
\[
\begin{tikzpicture}
\node(xa) at (-7.5,2.5) {$\frac{\l}{\mu}$};
\node(ya) at (-5,2.5) {$\frac{1}{\mu}F$};
\node(za) at (-2.5,2.5) {$\one$};
\node(aa) at (0,2.5) {$\displaystyle \frac{1}{\l} F$};
\node(ba) at (2.5,2.5) {$\displaystyle \frac{\mu}{\l} $};
\node(ca) at (5,2.5) {$\cdots$};
\node(xb) at (-7.5,0) {0};
\node(yb) at (-5,0) {0};
\node(zb) at (-2.5,0) {$\one$};
\node(ab) at (0,0) {$\displaystyle \frac{1}{\l} F$};
\node(bb) at (2.5,0) {$\displaystyle \frac{\mu}{\l} $};
\node(cb) at (5,0) {$\cdots$};
\tikzstyle{every node}=[font=\small]
\path[->,>=stealth',shorten >=1pt,auto,node distance=1.8cm]
(xa) edge node[above] {$-\frac{\l}{\mu}\a$} (ya)
(za)edge node[above] {$\frac{1}{\mu}\b $}  (ya)
(za) edge node[above] {$-\frac{1}{\l}\a$} (aa)
(ba) edge node[above] {$ \frac{1}{\l} \b$} (aa)
(ba) edge node[above] {$ -\frac{\mu}{\l^2} \a $} (ca)
(xb) edge node[above] {} (yb)
(zb)edge node[above] {}  (yb)
(zb) edge node[above] {$-\a$} (ab)
(bb) edge node[above] {$ \frac{1}{\l} \b$} (ab)
(bb) edge node[above] {$ -\frac{\mu}{\l^2} \a $} (cb)
(xa) edge node[left] {} (xb)
(ya) edge node[left] {} (yb)
(za) edge node[left] {$\Id$} (zb)
(aa) edge node[left] {$\Id$} (ab)
(ba) edge node[left] {$\Id$} (bb);
\end{tikzpicture}.
\]
We regard $f_{\a,\b}$ as an endomorphism of $\Cij{\a}{\b}$ of degree $\frac{\l}{\mu}$.  We refer to $f_{\a,\b}$ as the \emph{periodicity map} or \emph{periodicity endomorphism} of $\Cij{\a}{\b}$. 
\end{definition}

The following is clear from the definitions.
\begin{lemma}\label{lemma:coneFromCone}
We have $\Cone(f_{\a,\b})\simeq \mu\inv \Cone(\a)$.\qed
\end{lemma}

\begin{defn}\label{def:interpolationCx}
Suppose $\l_i\in \KS$ is a collection of invertible scalar objects indexed by a totally ordered finite set $I$ so that $\l_i\l_j\inv$ is small whenever $i > j$, and let $\a_i\in \Homg(\one,F)$ be maps of degree $\l_i$ ($1\leq i\leq r$).  Set $\Cij{i}{j}(F):=\Cij{\a_i}{\a_j}(F)$ for all $i\neq j$.  We define the \emph{interpolation complex} $\PB_i(F)$ by
\[
\PB_i(F)=\bigotimes_{j\neq i} \Cij{j}{i}(F).
\]
When the complex $F$ is understood, we omit it from the notation, writing $\PB_i=\PB_i(F)$. We may identify $I$ with the ordered set $\{0, 1, \ldots, r\}$ below.
\end{defn}

\begin{remark}\label{rmk:orderSmallness}
In case $\AC=\KC^\pm(\BC)$, a sufficient (but not necessary) condition for smallness is that the $\l_i$ are ordered by homological degree, c.f. Example \ref{ex:smallnessOfShifts}.
\end{remark}

Note that these definitions can all be made without assuming that $\{\a_i\}$ is a saturated collection of eigenmaps.

This definition of $\PB_i$ relies on a choice of order for the tensor product $\bigotimes_{j \ne i} \Cij{j}{i}$. However, when $\{\a_i\}$ is a saturated collection of eigenmaps, we will prove in Theorem \ref{thm-firstEigenthm} that $\PB_i$ is the projection to the $\a_i$-eigenspace (for any choice of order). Using this and its consequences, we prove in Lemma \ref{lem:orderdoesntmatter} that $\PB_i$ does not depend on the ordering of the tensor factors, up to homotopy equivalence.

In nice situations, the $\Cij{i}{j}$ commute with one another anyway.

\begin{proposition}
If $\{\a_i:\l_i\rightarrow F\}_{i=0}^r$ is a strongly commuting family, then the $\Cij{i}{j}$ tensor commute with one another, up to homotopy.
\end{proposition}
This is restated and proven in Proposition \ref{prop:zigzagsCommute}.

\subsection{A compact description of the interpolation complexes}
\label{subsec:compact}

In this section we develop some technology to discuss the precise relation between $\PB_i(F)$ and $\KB_i(F):=\bigotimes_{j\neq i}\Cone(\a_j)$, or between $\Cij{j}{i}$ and $\Cone(\a_j)$. The material in this section can be safely skipped, after the following notational convention.

\begin{definition}\label{def:periodicityMapPi} Retain notation as in Definition \ref{def:interpolationCx}. Let $u_i$ be formal indeterminates of degree $\l_i$. We will denote the
endomorphism $f_{ji}$ of $\Cij{j}{i}(F)$ from Definition \ref{def:periodicityMapCab} by $u_j/u_i$. By abuse, we denote the induced endomorphism of $\PB_i(F):=\bigotimes_{j\neq
i}\Cij{j}{i}(F)$ also by $u_j/u_i$. By convention $u_i/u_i$ is the identity morphism of $\Cij{j}{i}(F)$ or $\PB_i(F)$. \end{definition}

These maps make $\Cij{j}{i}$ a module over the ring $\Z[u_j/u_i]$, and $\PB_i$ a module over the polynomial ring $\Z[u_1/u_i,\ldots,u_r/u_i]$. That is, homogeneous polynomials act by
chain maps of the appropriate degree. This module structure is manifest from the construction of $\Cij{j}{i}$ as a periodic complex built from copies of $\Cone(\a_j)$ (with shifts), or
of $\PB_i$ as a multiply periodic complex built from copies of $\KB_i$. Here is a framework in which to discuss complexes which are also modules over polynomial rings.

\begin{notation} \label{not:polyMinus} Let $u$ be a formal indeterminate. Let $\Z[u]^+$ denote $\Z[u]$, viewed as a left module over itself. Let $\Z[u]^-$ denote the $\Z[u]$-module
\[
\Z[u]^- = \Z[u,u\inv]/u\Z[u].
\]
That is, $\Z[u]^-$ is isomorphic to $\Z[u^{-1}]$ as graded abelian groups, but with a locally nilpotent action of $u$.  If $u_1,\ldots,u_r$ are formal indeterminates and $\e=(\e_1,\ldots,\e_r)\in \{\pm 1\}$ is a sequence of signs, then we let $\Z[u_1,\ldots,u_r]^\e$ denote the tensor product over $\Z$ of $\Z[u_i]^{\e_i}$.

Now fix a small scalar object $\l$, and set $\deg(u) = \l$. Then we may also denote $\one[u]^\pm \ot (-)$ by $\Z[u]^\pm \ot_\Z (-)$, see Definition \ref{def:polyObjects}. The advantage of the latter notation is that it makes the $\Z[u]$-module structure on $\Z[u] \ot_\Z (-)$ clearer. \end{notation}

Fix $\l$ small. For an object $N \in \AC$, the complex
\[ \Z[u]^+ \ot_\Z N \cong N \oplus (\l N) \oplus (\l^2 N) \oplus \cdots \]
is a convolution indexed by $\Z_{\ge 0}$, with layers $\l^k N$, and the boring twisted differential $1 \ot d_N$. It also has a $\Z[u]$-module structure, meaning that $u$ commutes with this differential. However, $\Z[u] \ot_\Z N$ can be equipped with more interesting twisted differentials which commute with $u$, which are called \emph{$\Z[u]$-equivariant twisted differentials}. The total complex will therefore still have an action of $\Z[u]$. Any equivariant twisted differential has the form
\[ d = \sum_{k \ge 0} u^k \ot \pa_k,\]
for some linear maps $\pa_k$ from $N$ to $N$ of the appropriate degrees (not necessarily chain maps). We assume $\pa_0 = d_N$ so that the layers are still shifts of $N$.

Note that, if $\pa_k = 0$ for all $k \ge 2$, then
\[ d = 1 \ot d_N + u \ot \pa_1 \]
is a twisted differential if and only if $\pa_1$ is a chain map with $\pa_1^2 = 0$.

Similar statements can be made for $\Z[u]^- \ot N$.

We now give a compact description of $\Cij{\a}{\b}(F)$.

\begin{definition}\label{def:exteriorAction}
Let $\l$, $\mu$ be invertible scalar objects and let $\a:\l\rightarrow F$, $\b:\mu\rightarrow F$ be chain maps.  Define a chain map $\partial_{\a,\b}:\frac{\mu}{\l}[-1] \Cone(\a)\rightarrow \Cone(\a)$ by the following diagram:
\[
\begin{tikzpicture}[baseline=-4em]
\node at (4.9,0) {$\Big)$};
\node at (.7,0) {$\Big($};
\node at (1.3,-2) {$\Big)$};
\node at (-2.4,-2) {$\Big($};
\node at (-3.5,0) {$=$};
\node at (-3.5,-2) {$=$};
\node(f) at (-5.7,0) {$\displaystyle \frac{\mu}{\l}[-1] \Cone(\a)$};
\node(h) at (-5.7,-2) {$\Cone(\a)$};
\node(a) at (4,0) {$\displaystyle\frac{\mu}{\l} F[-1]$};
\node(b) at (1,0) {$\mu $};
\node(c) at (1,-2) {$F$};
\node(d) at (-2,-2) {$\l[1] $};
\tikzstyle{every node}=[font=\small]
\path[->,>=stealth',shorten >=1pt,auto,node distance=1.8cm]
(b) edge node[above] {$\a $} (a)
(b) edge node[right] {$-\b$} (c)
(d) edge node[above] {$\a$} (c)
(f) edge node[right] {$\partial_{\a,\b}$} (h);
\end{tikzpicture}.
\]
\end{definition}

\begin{lemma} \label{lemma:CabAsConv}
Continue the notation as above.  Let $u$ and $v$ denote formal indeterminates of degree $\l$ and $\mu$.  Assume that $\l\mu\inv$ is small.  Then
\begin{equation} \label{eq:CabAsConv}
\CB_{\a,\b}(F) = \Z\Big[\frac{u}{v}\Big]^+\otimes_\Z \mu\inv \Cone(\a)\ \ \ \ \ \text{ with differential } \ \ \ \ \ 1\otimes d + \frac{u}{v}\otimes \partial_{\a,\b}.
\end{equation}

\begin{equation} \label{eq:CbaAsConv}
\CB_{\b,\a}(F) = \Z\Big[\frac{v}{u}\Big]^-\otimes_\Z \mu\inv \Cone(\b)[-1] \ \ \ \ \ \text{ with differential } \ \ \ \ \ 1\otimes d -\frac{v}{u}\otimes \partial_{\b,\a},
\end{equation}
\end{lemma}

\begin{proof} This is obvious from the definitions. \end{proof}
	
Lemma \ref{lemma:coneFromCone} can be rephrased more generally as follows.

\begin{lemma} \label{lem:NconeM} Let $M = \Z[u]^+ \ot_\Z N$ with $\Z[u]$-equivariant twisted differential, and let $\l$ be the degree of $u$. Then $\Cone(\l M \buildrel u \over \to M) \cong N$. Similarly, when $M = \Z[u]^- \ot_\Z N$, then $\Cone(\l M \buildrel u \over \to M) \cong \l N[1]$. \qed \end{lemma}

There are analogous results for the complexes $\PB_i(F)$.

\begin{notation}
Let $u_1,\ldots,u_r$ be formal indeterminates.  For $\e=(\e_1,\ldots,\e_r)\in \{\pm 1\}^r$, let $\Z[u_1,\ldots,u_r]^\e:=\Z[u_1]^{\e_1}\otimes_\Z\cdots \otimes_\Z \Z[u_r]^{\e_r}$, thought of as a $\Z[u_1,\ldots,u_r]$-module.
\end{notation}

\begin{proposition}\label{prop:PiDescription_general}
Let $\l_0,\ldots,\l_r$ be invertible scalar objects such that $\l_i\l_j\inv$ is small when $i>j$.  Let $\a_i:\l_i\rightarrow F$ be given maps, and let $\partial_{i,j}$ denote the endomorphisms of $\Cone(\a_i)$ from Definition \ref{def:exteriorAction}.  Fix $i\in I$.  Then
\begin{equation}\label{eq:PasKoszulComplex}
\l_0\cdots\l_{i-1}\l_i^{r-i}\PB_i(F) \cong \Z\Big[\frac{u_0}{u_i},\ldots,\frac{\hat{u_i}}{u_i},\ldots,\frac{u_r}{u_i}\Big]^{\e}\otimes_{\Z} \KB_i[-i]
\end{equation}
with differential
\[
1\otimes d + \sum_{j<i}\frac{u_j}{u_i}\otimes \partial_{j,i}-\sum_{j>i}\frac{u_j}{u_i}\otimes \partial_{j,i},
\]
where $\KB_i:=\bigotimes_{j\neq i}\Cone(\a_j)$ and $\e$ is the sequence of signs
\[
\e_j=\begin{cases}
-1 & \text{ if $j>i$}\\
+1 & \text{ if $j<i$}
\end{cases}
\]
\end{proposition}

\begin{proof} This follows from Lemma \ref{lemma:CabAsConv}.\end{proof}


Moreover, there is an ``inverse'' to Proposition \ref{prop:PiDescription_general}, stating that $\KB_i$ can be reconstructed from $\PB_i$ by tensoring with a Koszul complex for the operators $u_j/u_i$. This is analogous to Lemma \ref{lemma:coneFromCone} or Lemma \ref{lem:NconeM}, which is the one variable case. We leave the details to the reader. Note that a variant of this Koszul complex reconstruction was already stated in Proposition \ref{prop:KasKoszulComplex}, though using eigenmaps instead of the periodicity maps. The comparison between these two approaches will be taken up in the next section.

%

\subsection{Compatibility of eigenmaps and periodicity}
\label{subsec:structure}

Let $\AC$ be a triangulated monoidal category with scalars $\KS\rightarrow \ZC(\AC)$, and let $F\in \AC$ be a fixed object.  Compare the following with Definition \ref{def:rhoi}.

\begin{definition}\label{def:rho}
Let $\a,\b\in \Homg(\one,F)$ be maps of degree $\l,\mu$, respectively, and suppose that $M\in \AC$ is a $\b$-eigenobject.  Then define $\rho_M(\a,\b)$ to be the composition
\[
\begin{diagram}
\l\otimes M & \rTo^{\a\otimes \Id_M} & F\otimes M & \rTo^{(\b\otimes \Id_M)\inv} & \mu\otimes M.
\end{diagram}
\]
Up to an application of $\mu\inv$, we think of $\rho_M(\a,\b)$ as an endomorphism of $M$ with degree $\l/\mu$.   We will abuse notation, denoting this morphism $\l\inv \rho_M(\a,\b)$ simply by $\displaystyle \frac{\a}{\b}$.
\end{definition}

\begin{remark}
If $\a\in \Homg(\one,F^{\otimes k})$ and $M$ is a $\b$-eigenobject, then we can define $\displaystyle \frac{\a}{\b^k}\in \Endg(M)$ in a similar way.  We will not use this construction for $k>1$ in this paper.
\end{remark}

Now suppose that $\{\a_i\}$ is a saturated collection of eigenmaps for $F$ satisfying the conditions of Definition \ref{def:interpolationCx}, and let $\PB_i=\PB_i(F)$ denote the interpolation complexes. We will prove shortly in Lemma \ref{lem:piaizero} that this implies $\PB_i \ot \Cone(\a_i) \cong \Cone(\a_i) \ot \PB_i \cong 0$, so that $\PB_i$ has an endomorphism $\frac{\a_j}{\a_i}$ for any $j \ne i$. By construction, the projectors $\PB_i$ are equipped with a commuting family of endomorphisms $u_j/u_i$ of degree $\l_j/\l_i$ (Definition \ref{def:periodicityMapPi}). It would be desirable if the periodicity maps $u_j/u_i$ were homotopic to the maps $\a_j/\a_i$.

Unfortunately, we do not know how to prove this, and believe it could be false in general.  Fortunately, it turns out that the vanishing of the same homological obstructions considered in \S \ref{subsec:commutativitystuff} will guarantee this equality (up to homotopy).  First we state a lemma which does not assume any sort of pre-diagonalizability.

\begin{lemma}\label{lemma:eigenaction}
Let $\l$ and $\mu$ be invertible scalar objects such that $\l\mu\inv$ is small, and let $\a:\l\rightarrow F$, $\\b:\mu\rightarrow F$ be chain maps.  Recall the periodicity maps from Definition \ref{def:periodicityMapCab}.  If $\{\a,\b\}$ is obstruction free, then
\begin{enumerate}
\item $\displaystyle\b \otimes f_{\a,\b}\simeq \a\otimes \Id$ as chain maps $\mu\otimes \Cij{\a}{\b}(F)\rightarrow F\otimes \Cij{\a}{\b}(F)$.\vskip8pt
\item $\displaystyle\a\otimes f_{\b,\a}\simeq \b\otimes \Id$ as chain maps $\l\otimes  \Cij{\b}{\a}(F)\rightarrow F\otimes \Cij{\b}{\a}(F)$.
\end{enumerate}
\end{lemma}
\begin{proof}
We will only prove statement (1).  Statement (2) is similar.  Let us ignore all occurences of $\l$ and $\mu$ for simplicity of notation; this abuse of notation is ultimately justified by the canonicalness of the isomorphism $\l\mu\simeq \mu\l$ (see Remark \ref{rmk:braiding}) .   Since $\{\a,\b\}$ is a strongly commuting family, we may choose homotopies $h_\a\in \Hom^{-1}(F, FF)$ so that the cycles $z_{\a,\b}$, $z_{\b,\a}$, $w_\a$, and $w_\b$ are boundaries.  In other words, there exist homotopies $k_{\a,\b}\in \Hom^{-2}(\one, FF)$ and $k_\a\in \Hom^{-2}(\one,FF)$ such that
\begin{eqnarray}
\label{eq:kHomotopy1} d\circ k_{\a,\b}- k_{\b,\a}\circ d &=& h_\b \circ \a + h_\a\circ \b\\
\label{eq:kHomotopy2} d\circ k_\a - k_\a\circ d &=& h_{\a}\circ \a.
\end{eqnarray}
and similarly with the roles of $\a$ and $\b$ reversed.  Given this, the following diagram is an explicit homotopy for $\a\otimes \Id -\b\otimes \frac{u}{v}$:
\[
\begin{tikzpicture}
\node(aa) at (0,2.5) {$ F$};
\node(ba) at (2.5,2.5) {$\one $};
\node(ca) at (5,2.5) {$ F$};
\node(da) at (7.5,2.5) {$\one $};
\node(ea) at (10,2.5) {$ F$};
\node(fa) at (12.5,2.5) {$\cdots$};
\node(ab) at (0,0) {$ FF$};
\node(bb) at (2.5,0) {$F $};
\node(cb) at (5,0) {$ FF$};
\node(db) at (7.5,0) {$F $};
\node(eb) at (10,0) {$ FF$};
\node(fb) at (12.5,0) {$\cdots$};
\node(z) at (12.5,1.25) {};
\tikzstyle{every node}=[font=\small]
\path[->,>=stealth',shorten >=1pt,auto,node distance=1.8cm]

(ba) edge node[above] {$ \a $} (aa)
(ba) edge node[above] {$ -\b$} (ca)
(da) edge node[above] {$\a$} (ca)
(da) edge node[above] {$ -\b $} (ea)
(fa) edge node[above] {$\a$} (ea)
(bb) edge node[below] {$ F\a $} (ab)
(bb) edge node[below] {$ -F\b$} (cb)
(db) edge node[below] {$F\a$} (cb)
(db) edge node[below] {$- F\b $} (eb)
(fb) edge node[below] {$F\a$} (eb)
%
(ba) edge node[xshift=-.8cm,yshift=-.4cm] {$-k_\a$} (ab)
(da) edge node[xshift=-.8cm,yshift=-.4cm] {$-k_\a$}(cb)
(fa) edge node[xshift=-.8cm,yshift=-.4cm] {$-k_\a$}(eb)
%
(aa) edge node[left] {$h_\a$}  (ab)
(ca) edge node[left] {$h_\a$}  (cb)
(ea) edge node[left] {$h_\a$}  (eb);
%
\draw[frontline,->,>=stealth',shorten >=1pt,auto,node distance=1.8cm]
(aa) to node[xshift=-1cm,yshift=0cm] {$\Id$}  (bb);
\draw[frontline,->,>=stealth',shorten >=1pt,auto,node distance=1.8cm]
(ba) to node[xshift=-.3cm] {$k_{\a\b}$}(cb);
\draw[frontline,->,>=stealth',shorten >=1pt,auto,node distance=1.8cm]
(ca) to node[xshift=-1cm,yshift=0cm]{$\Id$} (db);
\draw[frontline,->,>=stealth',shorten >=1pt,auto,node distance=1.8cm]
(da) to node[xshift=-.3cm] {$k_{\a\b}$} (eb);
\draw[frontline,->,>=stealth',shorten >=1pt,auto,node distance=1.8cm]
(ea) to node[xshift=-1cm,yshift=0cm] {$\Id$} (fb);
%
\draw[frontline,->,>=stealth',shorten >=1pt,auto,node distance=1.8cm]
(aa) to node[xshift=-.2cm,yshift=-.1cm]{$-h_\a$} (cb);
\draw[frontline,->,>=stealth',shorten >=1pt,auto,node distance=1.8cm]
(ca) to node[xshift=-.2cm,yshift=-.1cm] {$-h_\a$} (eb);
\draw[frontline]
(ea) to node[xshift=-.2cm,yshift=-.15cm] {} (z);
\end{tikzpicture}
.\]
\end{proof}

\begin{remark}
Our omission of the scalars $\l,\mu$ everywhere in the proof above obscures the fact that the morphisms $\l\otimes F\simeq F\otimes \l$ and $\l\otimes \mu\simeq \mu\otimes \l$, and even $\l\l\inv\simeq \one$ may be nontrivial homotopy equivalences. We leave it to the reader formulate a more precise argument involving all the relevant homotopies.
\end{remark}

\begin{proposition}\label{prop:eigenaction}
Suppose that $\{\a_i\}$ is a saturated collection of eigenmaps for $F$. If the collection of eigenmaps $\{\a_i\}_{i\in I}$ is obstruction-free, then
\[
\frac{u_j}{u_i} \simeq \frac{\a_j}{\a_i}
\]
as endomorphisms of $\idemp_i$, for all $j\in \{1,\ldots, \hat{i},\ldots,r\}$.  In other words, the periodicity endomorphisms of $\idemp_i$ coincide with the eigenmap endomorphisms.
\end{proposition}

\begin{proof}
Recall that $\idemp_i$ is defined to be $\bigotimes_{j\neq i}\Cij{j}{i}$, for some ordering of the factors.  The order is irrelevant up to homotopy, but let us be careful nonetheless.  Fix an ordering of the factors, and write $\idemp_i = A\otimes B \otimes C$, where $A$ and $C$ are tensor products of certain complexes of the form $\Cij{\ell}{i}$ with $\ell\neq i,k$, and $B:=\Cij{k}{i}$.  Given this, the periodicity map $\frac{u_k}{u_i}\in \Endg(\idemp_i)$ should more properly be written as $A\frac{u_k}{u_i}C$.  We have a commutative diagram
\[
\begin{tikzpicture}
\node(aa) at (0,2.5) {$ FABC$};
\node(ba) at (2.5,2.5) {$AFBC $};
\node(ca) at (5,2.5) {$AFBC$};
\node(da) at (7.5,2.5) {$FABC $};
\node(ab) at (0,0) {$\l_k ABC$};
\node(bb) at (2.5,0) {$A\l_k BC$};
\node(cb) at (5,0) {$A\l_k BC$};
\node(db) at (7.5,0) {$\l_k ABC$};
\tikzstyle{every node}=[font=\small]
\path[->,>=stealth',shorten >=1pt,auto,node distance=1.8cm]
(ab) edge node {$\a_k A B C$} (aa)
(bb) edge node {$A\a_k BC$} (ba)
(cb) edge node {$A \a_i (\frac{u_k}{u_i})C$} (ca)
(db) edge node {$\a_i A(\frac{u_k}{u_i})C $} (da)
(aa) edge node {$\simeq$} (ba)
(ba) edge node {$\Id$} (ca)
(ca) edge node {$\simeq $} (da)
(ab) edge node {$\simeq$}(bb)
(bb) edge node {$\Id$} (cb)
(cb) edge node {$\simeq $} (db);
\end{tikzpicture}
\]
The left and right squares commute by the proof of Lemma \ref{lemma:conesCommute}, particularly \eqref{eq:coneCommuteSquare}, and the middle square commutes by Lemma \ref{lemma:eigenaction}.  The compositions along the rows are homotopic to identity maps, hence commutativity of the diagram implies that 
\[
\a_k ABC \simeq (\a_i ABC) \circ \Big(A\frac{u_k}{u_i}C\Big)
\]
as maps $\l_k A B C \rightarrow F A B C$.  Postcomposing with the homotopy equivalence $F ABC = F \idemp_i \simeq \l_i \idemp_i$ yields
\[
\frac{\a_k}{\a_i} \simeq \frac{\a_i}{\a_i} \circ \Big(A\frac{u_k}{u_i}C\Big) \simeq A\frac{u_k}{u_i}C,
\]
as claimed.
\end{proof}

\begin{remark} \label{rmk:whysigns} In order for this proposition to hold, the signs introduced to the complex $\Cij{j}{i}$ were required. This, ultimately, is the reason for our sign convention.
\end{remark}

\subsection{The projector is an eigenobject}
\label{subsec-piaizero}

\begin{lemma} \label{lem:piaizero} We have $\PB_i \ot \Cone(\a_i) \cong \Cone(\a_i) \ot \PB_i \cong 0$. \end{lemma}

\begin{proof} Clearly $\KB_i \ot \Cone(\a_i) \cong \Cone(\a_i) \ot \KB_i \cong 0$, where $\KB_i$ is the tensor product of all the $\Cone(\a_j)$ for $j \ne i$. This was proven in Proposition \ref{prop:KBiprops}, and it doesn't depend on what order the tensor factors appear in $\KB_i$. Now $\PB_i$ is built as a convolution with layers $\KB_i$, so if we could use simultaneous simplifications to this convolution tensored with $\Cone(\a_i)$, we could deduce the lemma.
	
The convolution governing $\PB_i$ is locally finite, but it satisfies neither the ACC nor the DCC. It is indexed by an orthant of $\ZM^r$, those tuples of integers $(n_1, \ldots, n_r)$ where $n_k$ is either positive or negative depending on whether the $k$-th tensor factor of $\PB_i$, which is $\Cij{j}{i}$ for some $j \ne i$, satisfies $j<i$ or $j > i$. See Proposition \ref{prop:PiDescription_general} for the details.

Thankfully, this convolution for $\PB_i$ can be reassociated, as a convolution indexed by $\ZM_{\ge 0}^i$, with layers each of which is a convolution indexed by $\ZM_{\le 0}^{r-i}$, with layers each of which is $\KB_i$. This is done by lumping together the various coordinates where the orthant is positive, and the ones where it is negative. Now each individual convolution has either the ACC or the DCC, so by applying simulaneous simplifications twice we are done. \end{proof}

\section{The diagonalization theorem}
\label{sec:diag}

Let us recall our main theorem from the introduction, which will prove in this chapter.

\begin{thm}[Diagonalization Theorem]\label{thm-firstEigenthm} Let $\AC$ be a homotopy monoidal category (Definition \ref{def:homotopyMonoidalCat}) with category of scalars $\KS$.  Fix $F \in \AC$ which is categorically prediagonalizable with pre-spectrum $\{\a_i \co \l_i \to F\}_{i \in I}$. Assume that: \begin{itemize}
\item Each scalar object $\l_i$ is invertible.
\item The set $I$ is a finite totally-ordered set, which we identify with $\{0,1,\ldots,r\}$. Whenever $i > j$, the scalar object $\l_i \l_j\inv$ is small.
\end{itemize}
Under these assumptions, we have constructed complexes $\PB_i$ in \S\ref{subsec-cateInterpoly}. Then $\{\PB_i,\a_i\}_{i \in I}$ is a tight diagonalization of $F$, as in Definition \ref{def:catdiagzed}. \end{thm}

\subsection{Remarks about the proof}
\label{subsec-remarksonproof}

The fine details of the construction of $\PB_i$ are not actually essential to this proof. Here are the properties we use. \begin{itemize}
\item There is a complex $\Cij{j}{i}$ which is a locally finite convolution with layers $\Cone(\a_j)$. The poset governing this convolution satisfies either the ACC or the DCC.
\item $\PB_i$ is a tensor product of $\Cij{j}{i}$ for $j \ne i$, in some order.
\item There is a distinguished triangle (see Lemma \ref{lem:CijTriangle})
\[
\Cij{\b}{\a}(F)\rightarrow \one\rightarrow \Cij{\a}{\b}(F)\rightarrow \Cij{\b}{\a}(F)[1].
\]
\end{itemize}

We may be more precise for pedagogical reasons, but this is all that is necessary. This remark will be important if one wishes to generalize our results beyond the known cases, such as
to mixed eigenmaps (see \S\ref{sec:generalizations}).

\subsection{The relative diagonalization theorem}
\label{subsec-relative}

We will actually prove a refinement of Theorem \ref{thm-firstEigenthm} which is more technical to state, but more useful in applications.

\begin{theorem}[Relative Diagonalization Theorem]\label{thm:relDiag} Let $\AC$ be a monoidal homotopy category (Definition \ref{def:homotopyMonoidalCat}) with category of scalars $\KS$.  Fix $F \in \AC$ which is categorically prediagonalizable with pre-spectrum $\{\a_y \co \l_y \to F\}_{y \in \YC}$. Assume that: \begin{itemize}
\item Each scalar object $\l_y$ is invertible.
\item The set $\YC$ is a finite \emph{partially-ordered} set. Whenever $y_1 > y_2$ in $\YC$, the scalar object $\l_{y_1} \l_{y_2}\inv$ is small.
\end{itemize}
Now fix an idempotent decomposition of identity $\one\simeq \Tot\{\idemp_x, d\}$ indexed by a poset $\XC$, and assume that: \begin{itemize}
\item $\PB_{x}\otimes F\simeq F\otimes \PB_{x}$ for all $x\in \XC$.
\item For each $x\in \XC$ there exists a totally ordered subset $\YC_x \subset \YC$ such that
\[
\PB_{x}\otimes \bigotimes_{y\in \YC_x}\Cone(\a_y)\simeq 0.
\]
\end{itemize}
Then $F$ is diagonalizable.  More specifically, there is an $\XC\times \YC$-indexed diagonalization $\{(\PB_{(x,y)},\a_y)\}$ of $F$ such that reassociation recovers $\PB_x$:
\begin{equation}\label{eq:PxFromPxy}
\PB_{x} \simeq \Tot\{\PB_{(x,y)}\}_{y\in \YC}.
\end{equation}
\end{theorem}

\begin{remark}
It follows from \eqref{eq:PxFromPxy} $\PB_x\otimes \PB_{(x,y)}\simeq \PB_{(x,y)}\simeq \PB_{(x,y)}\otimes \PB_x$.  Further, if we define reassociated idempotents $\PB_{y} \simeq \Tot\{\PB_{(x,y)}\}_{x\in \XC}$, then $\PB_{(x,y)}\simeq \PB_x\otimes \PB_y$.  Some of these idempotents may be null-homotopic, so that the diagonalization $\{(\PB_x\otimes \PB_y,\a_y)\}$ is more properly indexed by a subset of $\XC\times \YC$. 
\end{remark}

\begin{remark}When $\XC=\{x\}$ is a singleton and $\PB_x = \one$, the partial order on $\YC$ must be a total order, and the construction of $\PB_{(x,y)}$ in the proof implies that $\PB_y$ is the
interpolation complex constructed in \S\ref{subsec-cateInterpoly}. Thus Theorem \ref{thm:relDiag} implies all but the tightness statement in Theorem \ref{thm-firstEigenthm}.  In fact, we do not expect any tightness results for the relative situation, because it is possible that eigencategories for
incomparable elements of $\YC$ may overlap.
\end{remark}

\subsection{The proof of relative diagonalization}
\label{subsec:theproof}

Assume the hypotheses of Theorem \ref{thm:relDiag}.

\begin{definition}\label{def:Pxy}
Define
\[
\PB_{(x,y)}:= \PB_x\otimes \bigotimes_{y'}\CB_{y',y},
\]
where the tensor product is over all $y'\in \YC$ which are comparable, but not equal to $y$.
\end{definition}

We first note that $\PB_x$ commutes with all the relevant complexes.

\begin{lemma}\label{lemma:PxCommuting}
We have $\PB_x\otimes \Cone(\a_y)\simeq \Cone(\a_y)\otimes \PB_x$ for all $x\in \XC$ and all $y\in \YC$.  Similarly $\PB_x\otimes \CB_{y,y'}\simeq \CB_{y,y'}\otimes \PB_x$ for all $x\in \XC$ and all comparable $y\neq y'\in \YC$.
\end{lemma}

\begin{proof}
Similar to \ref{lem:PConeCommute}.
\end{proof}

\begin{notation}
In the following, we will abbreviate by writing $\L_i:=\Cone(\a_i)$.  We will also frequently omit the tensor product symbol, writing $A\otimes B$ simply as $AB$.
\end{notation}

We begin with some lemmas. 
Recall that $\YC_x\subset \YC$ inherits a total order from $\YC$ and satisfies $\PB_x\otimes \bigotimes_{y\in \YC_x}\L_y\simeq 0$.

\begin{lemma}\label{lemma:tooManyCones}
Choose $x\in \XC$, and let $m\geq |\YC_x|$ be given.  Suppose we have maps $y_1,y_2:\{1,\ldots,m\}\rightarrow \{\YC_x\}$ so that $y_1$ is surjective.  Then $\PB_x\otimes \bigotimes_{i=1}^m \Cij{y_1(i)}{y_2(i)}\simeq 0$.
\end{lemma}

\begin{proof} Observe that $\bigotimes_{i=1}^m \Cij{y_1(i)}{y_2(i)}$ is the total complex of a convolution whose layers are shifted copies of $Z=\bigotimes_{i=1}^m \Lambda_{y_1(i)}$. If
each index $y\in \YC_x$ appears among the $y_1(i)$, then $\PB_x\otimes Z\simeq 0$ by hypothesis. Contracting all of the copies of $Z$ is legal by repeated use of Proposition
\ref{prop:simplifications}: the convolution is indexed by an orthant in $\ZM^m$, and we argue as in the proof of Lemma \ref{lem:piaizero}. \end{proof}

\begin{lemma}\label{lemma:PxyProps}
The complexes $\PB_{(x,y)}$ satisfy
\begin{enumerate}
\item $\PB_x\otimes \PB_{(x,y)}\simeq \PB_{(x,y)}\simeq \PB_{(x,y)}\otimes \PB_x$.
\item $\PB_{(x,y)} \otimes \PB_{(x',y')}$ is homotopy equivalent to $\PB_{(x,y)}$ if $(x,y)=(x',y')$, and is contractible otherwise.
\item $\Cone(\a_y)\otimes \PB_{(x,y)}\simeq 0 \simeq \PB_{(x,y)}\otimes \Cone(\a_y)$.
\item $\CB_{y,y'} \ot \PB_{(x,y)} \simeq 0 \simeq \PB_{(x,y)} \ot \CB_{y,y'}$ for all $y' \ne y$ comparable to $y$.
\item $\CB_{y',y}\otimes \PB_{(x,y)}\simeq \PB_{(x,y)}\simeq \PB_{(x,y)}\otimes \CB_{y',y}$ for all $y' \ne y$ comparable to $y$.
\end{enumerate}
\end{lemma}
\begin{proof}
Property (1) holds by construction, since $\PB_x$ is idempotent.
The proofs of (3) and (4) follow along the same lines as the proof of Lemma \ref{lemma:tooManyCones}.  Property (5) follows from (4) together with the existence of the distinguished triangle relating $\CB_{y,y'}$, $\CB_{y',y}$, and $\one$. Property (2) follows from (4) and (5), as well as the fact that $\PB_x$ and $\PB_{x'}$ are orthogonal for $x \ne x'$.
\end{proof}

It is now clear that the $\PB_{(x,y)}$ are mutually orthogonal idempotents, and statement (3) of the lemma is part of what we need to prove for diagonalization.  It remains to construct an idempotent decomposition of identity of the form
\begin{equation}\label{eq:relDecompOfOne}
\one\simeq \Tot\{\PB_{(x,y)}\}_{\XC\times \YC}.
\end{equation}

First recall (Lemma \ref{lem:CijTriangle}) that for each $y_1<y_2$ we have an equivalence
\begin{equation}\label{eq:shortDecomp}
\one\simeq (\CB_{y_2,y_1}\buildrel[1]\over \longrightarrow \CB_{y_1,y_2}).
\end{equation}

\begin{remark}\label{rmk:differentialOrder}
The differential strictly decreases the first coordinate, and strictly increases the second.
\end{remark}
Let us denote the right hand side of this equivalence by $\IB_b$, where $b=\{y_1,y_2\}$.    Now, let $B$ denote the subset of two-element subsets $\{y_1,y_2\}\subset \YC$ with $y_1\neq y_2$ comparable.   Choose once and for all an arbitrary total order on $B$, for the purposes of determining the order of tensor factors.  Tensoring together the complexes $\IB_b$ over all $b\in B$ yields
\begin{equation}\label{eq:crazyCube}
\one \simeq \bigotimes_{b\in B} \IB_b.
\end{equation}

The right hand side of the above is the total complex of a huge cube-like convolution of dimension equal to $|B|$ (i.e. having $2^{|B|}$ layers). In order to organize the layers of this convolution (i.e.~vertices of this cube) we introduce some interesting combinatorics.

Let $E$ denote the set of pairs $(b,y)$ with $b\in B$ and $y\in b$.  Let $\pi:E\rightarrow B$ denote the projection.

Let $\W$ denote the set of sections $B\rightarrow E$, that is the set of functions $\sigma:B\rightarrow E$ such that $\pi\circ \sigma = \Id_B$.  In other words, a section $\sigma$ is simply a choice of preferred element $y=\sigma(\{y,y'\})$ for each unordered two-element set of comparable elements $\{y,y'\}$.  The set of sections is partially ordered, via $\sigma\leq \sigma'$ if $\sigma(b)\leq \sigma'(b)$ for all $b\in B$.

Associated to each section $\sigma$ we have the complementary section $\sigma^\ast$ defined so that $b = \{\sigma(b),\sigma^\ast(b)\}$ for all $b\in B$.    

The right hand side of \eqref{eq:crazyCube} can be expanded into a convolution of the form
\begin{equation}\label{eq:sumOverSections}
\one\simeq \Tot\{\CB_\sigma,d\}_{\sigma\in \W}.
\end{equation}
Here,
\[
\CB_\sigma:=\bigotimes_{b\in B} \CB_{\sigma(b),\sigma^\ast(b)},
\]
in which the ordering of the tensor factors is assumed to be the same as in \eqref{eq:crazyCube}.

For each $y\in \YC$, we say that $\sigma$ is $y$-avoiding if $\sigma(b)\neq y$ for all $b\in B$.  Remark that if $y$ and $y'$ are comparable, then a section $\sigma$ cannot be both $y$-avoiding and $y'$-avoiding, since $\sigma(\{y,y'\})$ is either $y$ or $y'$.

The following is an immediate consequence of Lemma \ref{lemma:tooManyCones}.

\begin{lemma}
If $\sigma$ is not $y$-avoiding for any $y\in \YC_x$, then $\PB_x\otimes \CB_\sigma\simeq 0$.\qed
\end{lemma}

For each pair $(x,y)\in \XC\times \YC$ with $y\in \YC_x$, let $\W_{x,y}$ denote the set of $y$-avoiding sections $\sigma\in \W$.  Since $\YC_x$ inherits a total order from $\YC$, it is not possible for $\sigma$ to be $y$-avoiding and $y'$-avoiding with $y\neq y'\in \YC_x$.  Thus, $\W_{x,y}\cap \W_{x,y'}=\emptyset$ for $y\neq y'$.

\begin{lemma}
For each $x$, $\W$ can be partitioned into convex subsets
\[
\W = \bigsqcup_{y\in \YC_x}\W_{x,y}\sqcup \bigsqcup_{\sigma} \{\sigma\},
\]
where the disjoint union on the right is indexed by all sections $\sigma$ which are not $y$-avoiding for any $y\in \YC_x$.  This is a poset partition.
\end{lemma}

\begin{proof}
Fix $x\in \XC$ for the proof.  Suppose we have three sections $\sigma\leq \sigma'\leq \sigma''$.  If $\sigma,\sigma''\in \W_{x,y}$ for some $y\in \YC$.  We claim that $\sigma'\in \W_{x,y}$ as well.  For any $y'$ comparable to $y$, the above inequality implies
\[
\sigma(\{y,y'\})\leq \sigma'(\{y,y'\})\leq \sigma''(\{y,y'\}).
\]
The first and third are equal to $y'$ since $\sigma, \sigma''$ are $y$-avoiding.  Thus, the middle term equals $y'$ as well, and $\sigma'$ is $y$-avoiding.

Suppose now we have four sections $\sigma_1,\sigma_1
,\sigma_2,\sigma_2'$ such that $\sigma_i,\sigma_i'\in \W_{(x,y_i)}$ for $i=1,2$ and $\sigma_1\leq \sigma_2$, while $\sigma_2'\leq \sigma_1'$.   We must show that $y_2=y_1$.    By hypothesis, $\YC_x$ is totally ordered, hence $y_2\in \YC_x$ is comparable to $y_1$.  If $y_2\neq y_1$ then $\sigma_i(\{y_1,y_2\})$ and $\sigma_i(\{y_1,y_2\})$ are defined, and since $\sigma_i,\sigma'$ are $y_i$-avoiding, we have
\[
y_2=\sigma_1(\{y_1,y_2\})\leq \sigma_2(\{y_1,y_2\})=y_1
\]
and
\[
y_1=\sigma_2'(\{y_1,y_2\})\leq \sigma_1'(\{y_1,y_2\})=y_2.
\]
Thus $y_1=y_2$, a contradiction.  This proves that the disjoint union in the statement respects the partial order in $\W$.
\end{proof}

\begin{definition}
For each $y\in \YC_x$, let $\QB_{x,y}$ denote the contribution of the $y$-avoiding sections to the right-hand side of (\ref{eq:sumOverSections}).  That is,
\begin{equation}\label{eq:defQxy}
\QB_{x,y} = \Tot\{\CB_\sigma\}_{\sigma\in \W_{x,y}}.
\end{equation}
If $y\not\in \YC_x$, then we set $\QB_{x,y}:=0$.
\end{definition}

\begin{lemma}\label{lemma:PxyRes}
We have $\PB_x\otimes \QB_{x,y}\simeq \PB_{(x,y)}$, and there is an idempotent decomposition of identity $\one\simeq \Tot\{\PB_x\otimes \QB_{x,y},d\}_{(x,y)\in \XC\times \YC}$. 
\end{lemma}
\begin{proof}
Tensoring (\ref{eq:crazyCube}) with $\one\simeq \bigoplus_{x\in \XC}\PB_x$ (with twisted differential) yields
\[
\one\simeq \Tot\{\PB_x\otimes \CB_\sigma\}_{(x,\sigma)\in \XC\times \W}
\]
We reassociate as follows:
\[
\one\simeq \bigoplus_{x\in \XC}\left( \bigoplus_{y\in \YC_x} \PB_x\otimes \QB_{x,y}\ \oplus \ \underbrace{\bigoplus_{\sigma} \PB_x\otimes \CB_\sigma}_{\simeq 0}\right) \ \ \ \text{ with twisted differential}.
\]
where in the sum on the right, $\sigma$ runs over the sections which are not $y$-avoiding for any $y\in \YC_x$.  If $\sigma$ is such a section, then $\PB_x\otimes \CB_\sigma\simeq 0$.  Contracting these terms yields a decomposition of identity
\[
\one\simeq \Tot\{\PB_x\otimes \QB_{x,y}\}_{(x,y)\in \XC\times \YC}.
\]
The differential respects the product partial order on $\XC\times \YC$ in the sense that $\PB_x\otimes \QB_{x,y}$ maps to $\PB_{x'}\otimes \QB_{x',y'}$ only if $x\leq x'$ and $y\leq y'$ (see Remark \ref{rmk:differentialOrder}).

To prove Theorem \ref{thm:relDiag} it remains only to show that $\PB_x\otimes \QB_{x,y}\simeq \PB_{(x,y)}$.  To prove this we first factor $\QB_{x,y}$ in the following way
\[
\QB_{x,y} = \bigotimes_{\{y,y'\}\in B} \CB_{y',y} \otimes \bigotimes_{b\not\owns y} \IB_b,
\]
where $\IB_b$ is the right-hand side of (\ref{eq:shortDecomp}), and the ordering of the tensor factors is the same as in (\ref{eq:crazyCube}).  That $\QB_{x,y}$ can be written in this way follows easily from its definition (\ref{eq:defQxy}). Each $\IB_b$ is homotopy equivalent to $\one$, hence
\[
\QB_{x,y}\simeq \bigotimes_{\{y,y'\}\in B}\CB_{y',y}.
\]
This tensor product can also be regarded as the tensor product over all elements $y'\neq y\in \YC$ which are comparable to $y$.  Thus $\PB_x\otimes \QB_{x,y}\simeq \PB_{(x,y)}$, by definition of $\PB_{(x,y)}$.
\end{proof}

This concludes the proof of Theorem \ref{thm:relDiag}.

\subsection{Observations about the relative projectors}

We never proved or assumed that the complexes $\L_y$ or $\Cij{y'}{y}$ tensor commute up to homotopy. The definition of $\PB_y$ therefore depends on a choice of order for the tensor product of the various $\Cij{y'}{y}$. In fact, the projectors $\PB_{(x,y)}$ are independent of this choice of order.

\begin{lemma} \label{lem:orderdoesntmatter} In Definition \ref{def:Pxy}, the ordering of the tensor factors is irrelevant up to isomorphism.\end{lemma}

\begin{proof} Let $\PB_{(x,y)}'$ denote a tensor product with another order.  Lemma \ref{lemma:PxyProps} implies that $\PB_{(x,y)}' \ot \PB_{(x,y)} \simeq \PB_{(x,y)}$ because $\PB_{(x,y)}$ absorbs each $\Cij{y'}{y}$ in with $y'$ comparable to $y$.  A similar argument shows that $\PB_{(x,y)}' \ot \PB_{(x,y)} \simeq \PB_{(x,y)}$, hence $\PB_{(x,y)}'\simeq \PB_{(x,y)}$.\end{proof}

Let us provide a slightly smaller description of $\PB_{(x,y)}$ which is sometimes useful.
\begin{lemma}\label{lemma:smallerDescription}
We have
\[
\PB_{(x,y)}\simeq \PB_x\otimes \bigotimes_{y'\in \YC_x\setminus\{y\}}\CB_{y',y}.
\]
\end{lemma}

\begin{proof}
Let $\PB_{(x,y)}'=\PB_x\otimes \bigotimes_{y'\in \YC_x\setminus\{y\}}\CB_{y',y}$.  The proof of Lemma \ref{lemma:PxyProps} carries over nearly verbatim to $\PB_{(x,y)}'$.  Then the equivalence $\PB_{(x,y)}'\simeq \PB_{(x,y)}$ proceeds as in the proof of Lemma \ref{lem:orderdoesntmatter}.
\end{proof}

\subsection{Tightness}
\label{subsec:tightness}

Now we prove the final piece of Theorem \ref{thm-firstEigenthm}, that the diagonalization is tight. We do not have a tightness result for relative diagonalization (nor do we expect one), but we begin the discussion in this context. Assume that $\AC$ acts on $\VC$, and $M$ is an object of $\VC$.

\begin{lemma} One has $\L_y \ot M \simeq 0$ if and only if $\Cij{y}{y'} \ot M \simeq 0$ for a single $y' \ne y$ comparable, if and only if $\Cij{y}{y'} \ot M \simeq 0$ for all $y' \ne
y$ comparable. \end{lemma}

\begin{proof} Note that $\Cij{y}{y'}$ is isomorphic to a convolution with layers $\L_y$, by construction, and $\L_y$ is a convolution with layers $\Cij{y}{y'}$ (by taking the cone of
the periodicity map, Lemma \ref{lemma:coneFromCone}). Both these convolutions are locally finite and have a poset with either the ACC or the DCC, so simultaneous simplifications can
apply. From this, the lemma is obvious. \end{proof}

\begin{cor} \label{cor:spbfom} If $\L_y \ot M \simeq 0$ then $\PB_{y'} \ot M \simeq 0$ for any $y' \ne y$ comparable in $\YC$. \end{cor}

\begin{proof} If $\L_y \ot M$ then $\Cij{y}{y'} \ot M \simeq 0$ for all $y' \ne y$ comparable, which implies that $\PB_{(x,y')} \ot M \simeq 0$ as well, for all $x \in \XC$ (this uses
Lemma \ref{lem:orderdoesntmatter} to place the $\Cij{y}{y'}$ tensor factor of $\PB_{(x,y')}$ at the end). Thus $\PB_{y'} \ot M \simeq 0$ for all $y' \ne y$ comparable. \end{proof}

Applying the idempotent decomposition $\one\simeq \Tot\{\idemp_y, d\}$ to $M$, we see that $M \simeq \idemp_y \ot M \oplus \oplus (\bigoplus \idemp_{z} \ot M, d)$ where $z$ ranges over
those elements of $\YC$ which are incomparable to $y$. Because the idempotent decomposition respected the partial order on $\YC$, there are no terms in the differential between
$\idemp_y$ and the various $\idemp_z$, though there may be differentials between different $z$ incomparable to $y$ but comparable to each other.

Thankfully, in the situation of Theorem \ref{thm-firstEigenthm}, $\YC$ is equipped with a total order, and Corollary \ref{cor:spbfom} implies that $\PB_y \ot M \simeq M$. Thus the
diagonalization is tight.


\section{Generalizations and the Casimir element}
\label{sec:generalizations}

There are many categorical notions which could serve as a version of ``eigenvalue" and ``eigenobject," and one must judge their worth by the value these notions add to their study, and
by the presence of examples. The theory of forward eigenmaps pursued in this paper will be extremely useful in categorical Hecke theory. Everything done in this paper with forward
eigenmaps can be adapted to backward eigenmaps mutatis mutandis. However, these are not the only possibilities.

This short section illustrates some of the places ripe for further exploration. We begin with a generalization we hope will have applications: mixed eigenmaps.

\subsection{Mixed eigenmaps}
\label{subsec:mixed}

As usual, we assume that $\AC$ is a monoidal homotopy category, acting on itself (or perhaps on another homotopy category $\VC$) by tensor product.  We abuse notation by referring to objects of $\AC$ as functors, with the understanding that $F\in \AC$ corresponds to $F\otimes -$. 

\begin{defn}\label{def:mixedEigenmap} Fix $F \in \AC$. Suppose that $\lambda$ and $\mu$ are scalar functors, and $\a \co \lambda \to F$ and $\b \co F \to \mu$ are natural transformations satisfying $\b
\circ \a = 0$. We call $M \ne 0$ a \emph{split $(\a,\b)$-eigenobject of $F$} if both $\a_M \co \lambda(M) \to F(M)$ and $\b_M \co F(M) \to \mu(M)$ are split, giving together an isomorphism
$F(M) \cong \lambda(M) \oplus \mu(M)$.  In this case, we call $(\a,\b)$ a \emph{mixed eigenmap of $F$}, having \emph{(categorical) eigenvalue} $\lambda \oplus \mu$. \end{defn}

Forward or backward eigenmaps are examples of mixed eigenmaps, where one of the scalars $\mu$ or $\l$ is zero.  In contrast with the case of forward and backward eigenmaps, the category of split mixed eigenobjects is typically not closed under mapping cones.

\begin{example} \label{ex:badmixed} Work in the homotopy category of $\Z$-modules.  Let
\[
F = \Z \buildrel 2\over \rightarrow \underline{\Z}.
\]
There is a mixed eigenmap $(\a,\b)$ with eigenvalue $\Z \oplus \Z[1]$, where $\a$ is the inclusion of $\Z$, and $\b$ is the projection to $\Z[1]$.  If $C$ is a complex of finitely generated abelian groups, then $C$ is a split $(\a,\b)$-eigenobject if and only if $C$ is isomorphic to a direct sum of shifts of $\Z/2$ and $F$.  On the other hand, there exists a distinguished triangle
\[
F \rightarrow C \rightarrow F \rightarrow F[1]
\]
where $C:=\Cone(\Z\buildrel 4\over\rightarrow \Z)$ is not a split $(\a,\b)$-eigenobject.  This shows that the split eigenobjects are not closed under mapping cones.
\end{example}

One could attempt to fix this by extending the notion of split eigenobject as follows.

\begin{defn}\label{def:mixedEigencone} Given a mixed eigenmap $(\a,\b)$ as above, define the $(\a,\b)$-\emph{eigencone} to be the convolution \begin{equation} \L_{\a,\b} := \Big(\lambda[1]
\buildrel \a \over \longrightarrow F \buildrel \b \over \longrightarrow \mu[-1] \Big), \end{equation}
which exists since $\b \circ \a = 0$.  We say that $M$ is a \emph{generalized eigenobject} for $(\a,\b)$ if $\L_{\a,\b}\otimes M\simeq 0$.  The full subcategory of $\VC$ consisting of the generalized $(\a,\b)$-eigenobjects will be called the \emph{$(\a,\b)$-eigencategory}.
\end{defn}

Note that $\b$ induces a map $\hat{\b}:\Cone(\a)\rightarrow \mu$ and $\a$ induces a map $\hat{\a}:\l\rightarrow \Cone(F\rightarrow \mu)$.
\begin{remark}
The eigencategories for $\hat{\a}$, $\hat{\b}$, and $(\a,\b)$ coincide.  Each is triangulated.
\end{remark}

Observe that if $\L_{\a,\b}\otimes M\simeq 0$, then
\begin{equation}\label{eq:nonsplit}
F\otimes M \simeq (\mu M\buildrel [1]\over\rightarrow \l M).  
\end{equation}
Occasionally, one is lucky and the connecting morphism $\d : \mu M \rightarrow \l M[1]$ in the above convolution is null-homotopic.  In this case, $M$ is a split $(\a,\b)$-eigenobject in the sense of Definition \ref{defn:mixeddiag}.
Of course, one is not always so lucky.

\begin{ex}\label{ex:badmixed2}
In the setting of Example \ref{ex:badmixed}, $\L_{\a,\b}\simeq 0$, hence every object of $\VC$ is a generalized $(\a,\b)$-eigenobject.
\end{ex}

\begin{remark} The fact that $\l$ and $\mu$ differ in homological degree by exactly one in Example \ref{ex:badmixed} is what permits the possiblility of a nonsplit extension between $\l
M$ and $\mu M$ for an object $M$ of the additive category.  
\end{remark}

\begin{remark} In important examples, the space of chain maps $\mu M\rightarrow \l M[1]$ vanishes for degree reasons. For instance, categorical quantum groups and categorical Hecke
algebras possess interesting triangulated \emph{parity} subcategories. For a complex in this parity subcategory, the Hom complex of endomorphisms is concentrated in purely even
homological degree. Therefore, if $\l$ and $\mu$ have the same parity of homological shift, there can be no degree 1 maps from $\mu M$ to $\l M$ for any parity object $M$, and the split
eigencategory inside the parity subcategory is triangulated. Because of this feature, we believe mixed eigenmaps are best suited to these settings. \end{remark}

In any case, the following is one way to exclude pathological situations such as Example \ref{ex:badmixed2}.

\begin{definition}
Let $\a:\l\rightarrow F$, $\b:F\rightarrow \mu$ be given.  Form the convolution $\L_{\a,\b}$ as in Definition \ref{def:mixedEigencone}.  We say that $(\a,\b)$ is a \emph{(split generated) mixed eigenmap} if the $(\a,\b)$-eigenobjects is generated by the split $(\a,\b)$-eigenobjects with respect to sums, summands, shifts, mapping cones, and locally finite convolutions.
\end{definition}

\begin{example}
Suppose $F$ is an invertible complex and let $\a:\l\rightarrow F$ be a forward eigenmap.  Let $\a^\ast:F\inv\rightarrow \l\inv$ be the dual map.  Let $G:=F\oplus F\inv$ with $\b_1:=(\a,0):\l\rightarrow G$ and $\b_2:=(0,\a^\ast):G\rightarrow \l\inv$.  Then
\[
\L_{\b_1,\b_2}= \Cone(\a)\oplus \Cone(\a^\ast) \simeq \Cone(\a)\oplus \l\inv F\inv \Cone(\a)
\]
It follows that $\L_{\b_1,\b_2}\otimes M\simeq 0$ iff $\Cone(\a)\otimes M\simeq 0$.  Thus  every generalized $(\b_1,\b_2)$-eigenobject is a $(\b_1,\b_2)$-eigenobject, and $(\b_1,\b_2)$ is a split-generated mixed eigenmap.
\end{example}

This trivial looking example is a special case in a conjectural family of mixed eigenmaps for Casimir complexes.

\begin{ex} \label{ex:casimirInSbim} Recall the setup of \S\ref{subsec:typeA1}. Let $C = \FT_2(-1)[1] \oplus \FT_2\inv(1)[-1]$. Then any eigenobject for $\FT_2$ is also a split eigenobject
for an appropriate mixed eigenmap of $C$. For instance $B_s$ is a split eigenobject of $C$ with eigenvalue $\one(-3)[1]\oplus \one(3)[-1]$. Under categorical Schur-Weyl duality, the object $C$ above corresponds to the Casimir complex of quantum $\sl_2$ acting on the zero weight space in an $\sl_2$ categorification of
$(\C^2)^{\otimes 2}$. We will discuss this briefly in \S\ref{subsec:casimir}. \end{ex}

\begin{defn} \label{defn:mixeddiag} We say that $F$ is \emph{categorically prediagonalizable} with \emph{prespectrum} consisting of mixed eigenmaps $\{ (\a_i, \b_i) \}$ if the
obvious analogs of (PD1), (PD2), and (PD3) from Definition \ref{def:prediag} hold, replacing the usual eigencones with $\L_{\a_i,\b_i}$. In particular, \[ \bigotimes
\L_{\a_i, \b_i} \simeq 0 \] for any ordering of the tensor factors. The definition of what it means for $F$ to be \emph{categorically diagonalized} is again the obvious analog of Definition \ref{def:catdiag}. \end{defn}

\begin{remark} It is possible to write down homological obstructions for the commutativity of the mixed eigencones $\L_{\a_i,\b_i}$, similar to what is done in the appendix.
\end{remark}

Given a prediagonalizable functor $F$ with mixed eigenmaps $\a_i \co \l_i \to F$, $\b_i \co F \to \mu_i$, one might hope for an analog of Theorem \ref{thm-firstEigenthm}, at least when
the scalar functors $\l_i$ and $\mu_i$ are all invertible.  We do not yet know how to accomplish this, because we do not know how to construct the appropriate complexes $\Cij{j}{i}$.
Essentially, we do not know how to interpret the denominator of the Lagrange interpolation formula as a complex; we do not know the appropriate analog of Koszul duality which applies to this context.

With this future application in mind, we have written the proof of Theorem \ref{thm-firstEigenthm} in a way that will adapt to mixed eigenvalues, once the appropriate complexes $\Cij{j}{i}$ are constructed. That is, suppose that one can find a total order on the mixed eigenmaps and can construct objects $\Cij{j}{i}$ satisfying \begin{itemize}
\item $\Cij{j}{i}$ is homotopy equivalent to a convolution (which is locally finite, indexed by a bounded above poset) whose layers are $\L_{\a_j,\b_j}$.
\item Conversely, $\L_{\a_j,\b_j}$ is homotopy equivalent to a convolution whose layers are $\Cij{j}{i}$.
\item If $i < j$ then there is a resolution of identity $\1 \simeq \left( \Cij{i}{j} \rightarrow \Cij{j}{i} \right)$.
\end{itemize}
Then the entire argument of \S\ref{subsec:theproof} holds almost verbatim, and one can deduce that $\PB_i = \bigotimes_{j \ne i} \Cij{j}{i}$ will be a categorical eigenprojection.

\subsection{The Casimir operator}
\label{subsec:casimir}

We now give a conjectural example of mixed eigenmaps, in the context of $\sl_2$ categorification. It takes a great deal of background and setup to present this here completely, and since
we have no diagonalization theory yet to apply, it is not worth it. We merely wish to record the idea for posterity, because we imagine there are many interested readers. We assume the
reader is familiar with Lauda's 2-category $\UC(\sl_2)$ from \cite{LauSL2}, and only outline the proofs of various assertions. Note that $\UC(\sl_2)$ is a 2-category with objects $n \in
\Z$, so that the endomorphism category of $n$ (or rather, its homotopy category) is a monoidal category where we can discuss categorical diagonalization.

In work of Beliakova-Khovanov-Lauda \cite{BKL}, a complex $C_{n}$ was constructed for each $n \in \Z$, which categorifies the action of the Casimir operator in the quantum group of $\sl_2$
as it acts on the $n$-th weight space.\footnote{They denote this complex $C' \1_{n}$.} Here is their complex $C_{n}$ (see \emph{loc.~cit.~} for the precise differentials).

\[
\begin{tikzpicture}
\node(a) at (0,2) {$FE\one_n(-2)$};
\node(b) at (0,0) {$\one_n(-n-1)$};
\node(c) at (3,2) {$FE\one_n$};
\node(d) at (3,0) {$FE\one_n$};
\node(e) at (6,2) {$FE\one_n(2)$};
\node(f) at (6,0) {$\one_n(n+1)$};
\node at (0,1) {$\oplus$};
\node at (3,1) {$\oplus$};
\node at (6,1) {$\oplus$};
\tikzstyle{every node}=[font=\small]
\path[->,>=stealth',shorten >=1pt,auto,node distance=1.8cm]
(a) edge node[above] {} (c)
(a) edge node[above] {} (d)
(b) edge node[above] {} (d)
(c) edge node[above] {} (e)
(c) edge node[above] {} (f)
(d) edge node[above] {} (f);
\draw[frontline,->,>=stealth',shorten >=1pt,auto,node distance=1.8cm]
(b) to node[xshift=-1cm,yshift=0cm] {}  (c);
\draw[frontline,->,>=stealth',shorten >=1pt,auto,node distance=1.8cm]
(d) to node[xshift=-.3cm] {}(e);
\end{tikzpicture}
\]

When the Casimir operator acts on a highest weight representation with highest weight $\l \in \Z_{\ge 0}$, it acts by multiplication by the scalar $-(q^{\l+1} + q^{-\l-1})$. One should
expect this scalar to be categorified by the functor $\1(\l+1)[-1] \oplus \1(-\l-1)[1]$. After all, let $M$ be an object in a 2-representation of highest weight $\l$, living in the highest weight space: we call $M$ a \emph{highest weight object}. Then $E M \cong 0$, so that only two terms in $C_{\l}$ act nontrivially on $M$, and \[C_{\l} M \cong
M(\l+1)[-1] \oplus M(-\l-1)[1].\]

Fix $\l \ge 0$. One might hope then, for any $n \in \Z$, to construct chain maps \[ \a_{\l,n}:\one_{n}(-\l-1)[1]\rightarrow C_{n}, \quad \b_{\l,n}:C_{n}\rightarrow \one(\l+1)[-1],\] such that $(\a_{\l,n},\b_{\l,n})$ is a mixed eigenmap for any object in weight $n$ inside a categorical representation of highest weight $\l$.
\[
\begin{tikzpicture}
\node(a) at (0,2) {$FE\one_n(-2)$};
\node(b) at (0,0) {$\one_{n}(-n-1)$};
\node(c) at (3,2) {$FE\one_{n}$};
\node(d) at (3,0) {$FE\one_{n}$};
\node(e) at (6,2) {$FE\one_{n}(2)$};
\node(f) at (6,0) {$\one_{n}(n+1)$};
\node(g) at (-2,-2) {$\one_{n}(-\l-1)[1]$};
\node(h) at (8,4) {$\one_{n}(\l+1)[-1]$};
\node at (0,1) {$\oplus$};
\node at (3,1) {$\oplus$};
\node at (6,1) {$\oplus$};
\tikzstyle{every node}=[font=\small]
\path[->,>=stealth',shorten >=1pt,auto,node distance=1.8cm]
(g) edge node[above] {} (a)
(g) edge node[above] {} (b)
(e) edge node[above] {} (h)
(f) edge node[above] {} (h)
(a) edge node[above] {} (c)
(a) edge node[above] {} (d)
(b) edge node[above] {} (d)
(c) edge node[above] {} (e)
(c) edge node[above] {} (f)
(d) edge node[above] {} (f);
\draw[frontline,->,>=stealth',shorten >=1pt,auto,node distance=1.8cm]
(b) to node[xshift=-1cm,yshift=0cm] {}  (c);
\draw[frontline,->,>=stealth',shorten >=1pt,auto,node distance=1.8cm]
(d) to node[xshift=-.3cm] {}(e);
\end{tikzpicture}
\]
Unfortunately, one can compute that no such chain maps exist.

Fix $\l \ge 0$, and consider instead $\UC(\sl_2)_\l$, the \emph{$\l$-th cyclotomic quotient} of $\UC(\sl_2)$. This is the quotient 2-category where one kills any diagram possessing
either \begin{itemize} \item Any weight $\one_\mu$ for $\mu > \l$. \item Any bubble (real or fake) of positive degree inside a region labeled $\l$. \end{itemize} Preliminary
investigations lead us to conjecture that, in this quotient, one can construct the chain maps $(\a_{\l,n},\b_{\l,n})$ for all $n$, in such a way that they commute with the actions of $E$
and $F$. \footnote{Beliakova-Khovanov-Lauda provide chain maps for the commutation isomorphisms $C_{n+2} \ot E \cong E \ot C_n$ and $C_{n-2} \ot F \cong F \ot C_n$, proving that the
family of maps $C_n$ is in the Drinfeld center. Our conjectural maps $\a_{\l,n}$ and $\b_{\l,n}$ are morphisms in the Drinfeld center, and hence so is the cone $\Lambda_{\a,\b}$.} The
maps $\a_{\l,\l}$ and $\b_{\l,\l}$ would just be the inclusions of the obvious chain factors.

Given these conjectural mixed eigenmaps, it is not difficult to prove that any object in an isotypic highest weight $\l$ 2-representation of $\UC(\sl_2)_\l$ is actually a split
eigenobject for $(\a_{\l,n},\b_{\l,n})$. We have already discussed how any highest weight object would be a split eigenobject. Using Drinfeld centrality, we deduce that $F^k M$ is a
split eigenobject, for any highest weight object $M$. One can use results of Chuang and Rouquier \cite{ChuRou} to prove that (in an abelian 2-representation of $\UC(\sl_2)_\l$, or its
derived category) any weight category is generated by the objects $F^k M$ for highest weight objects $M$.

It still does not seem possible to construct $(\a_{\mu,n},\b_{\mu,n})$ for $0 \le \mu < \l$. The obstruction to the existence of such a chain map is easy to measure, and perhaps more
sophisticated generalizations of eigenmaps will allow one to solve this conundrum. Alternatively, it may be possible to reconstruct the cyclotomic quotient $\UC(\sl_2)_{\l - 2}$ if one
manages to construct the projection away from the $(\a_{\l,n},\b_{\l,n})$-eigencategory inside $\UC(\sl_2)_{\l}$.

In some special cases, one can construct enough mixed eigenmaps to pre-diagonalize the Casimir operator.

\begin{ex} A categorical version of Schur-Weyl duality is constructed by Mackaay-Stosic-Vaz in \cite{MSV09}. The \emph{Schur quotient} $\SC(2,2)$ of $\UC(\sl_2)$ is obtained by killing
all diagrams with weights not appearing in the representation $(\C^2)^{\ot 2}$ of $\sl_2$. This quotient factors through the cyclotomic quotient $\UC(\sl_2)_2$. They give a functor from
the Soergel category to the endomorphism category of the $0$-weight space in $\SC(2,2)$, sending $B_s$ to $FE \one_0$. Under this functor, the complex $C$ from Example
\ref{ex:casimirInSbim} is sent to the Casimir complex $C_0$, and the mixed eigenmap for $B_s$ is exactly our mixed eigenmap $(\a_{2,0},\b_{2,0})$. Also note that $B_s$ is in the image
under $F$ of a highest weight object.

One can also construct a mixed eigenmap for the isotypic component of highest weight $0$ (using the other eigenmap of $\FT_2$). The arrows point in the opposite direction to those of
$(\a_{2,0},\b_{2,0})$, projecting away from the $FE \one_0$ terms. We leave it as an exercise to prove that $C_0$ is categorically pre-diagonalizable. \end{ex}

\subsection{Eigenconvolutions}
\label{subsec:eigenconvs}

Let $\langle F\rangle $ denote the smallest full triangulated subcategory of $\AC$ It is sensible to suggest that the true generalization of an eigenmap, or rather of an eigencone, is an object in the triangulated hull convolution built out of the functor $F$ (appearing once) and various scalar
functors. Even if $F \ot M$ is a non-split extension of various scalar functors applied to $M$, this is not the most unreasonable definition of a categorical eigenobject. Because we have
no interesting examples, we say no more.

\subsection{Generalized eigenobjects}
\label{subsec:geneigen}

In linear algebra, a \emph{generalized eigenvector} $m$ with \emph{generalized eigenvalue} $\lambda$ for an operator $f$ satisfies $(f - \lambda)^{N} m = 0$ for some $N$. Clearly one
could define a \emph{generalized eigenmap/eigenobject} analogously, by the condition $\Lambda_{\a,\b}^{\ot N} \ot M = 0$ for a mixed map $(\a,\b)$. Again, we have no interesting examples, so
remains to be seen whether this definition is a fruitful one.

\subsection{Powers}
\label{subsec:powers}

If $f$ is diagonalizable then so is $f^2$, and the eigenvalues of $f^2$ are the squares of the eigenvalues of $f$.

Given a forward eigenmap $\a$ of $F$ with eigenobject $M$, one has an eigenmap $\a \ot \a$ of $F \ot F$ with eigenobject $M$. It seems possible that the eigencategory of $\a \ot \a$ may
be strictly larger than the eigencategory of $\a$. It is certainly possible that $F \ot F$ has more eigenmaps and eigencategories than those of the form $\a \ot \a$. After all, this
happens for $G$ and $F = G \ot G$ from \S\ref{subsec:nonexamples}.

Inspired by the example of $\OC(1)$ in \S\ref{subsec:cohp1}, one may wish to study the eigenmaps of the entire family $F^{\ot n}$ for $n \ge 0$. We do not pursue this idea any further,
but it clearly has interest.

\subsection{Eigenroofs}
\label{subsec:roofs}

We shift tacks to another potentially useful generalization of an eigenmap. Here is a motivational example from the categorical representation theory of Hecke algebras.

Let $\FT_n$ denote the full twist on $n$ strands, or its corresponding
Rouquier complex. Then the Jucys-Murphy element can be defined as $J_n = \FT_n \FT_{n-1}^{-1}$.  As discussed in \S\ref{subsec:HeckeConj}, the complex $\FT_n$ admits many forward eigenmaps, and $\FT_{n-1}\inv$ admits many backward eigenmaps, but not vice versa. Typically $J_n$ admits neither forward nor backward eigenmaps!\footnote{Nor does it admit mixed eigenmaps. The eigenvalues are invertible, and hence indecomposable, so any mixed eigenmap would be either a forward or a backward eigenmap.}

For example, let $M\in \K(\SBim_3)$ denote the complex
\[
\un{B_s}(-3)\rightarrow B_{sts}(-1) \rightarrow B_{sts}(1)\rightarrow B_s(3)
\]
where the first and last maps are canonical (they live in one-dimensional morphism spaces), and the third map is multiplication by $x_3\otimes 1 - 1\otimes x_3$ in the ground ring.  This complex is an eigencomplex for $\FT_2\inv$ with eigenmap $\a:\FT_2\inv\rightarrow \one(2)$, and is an eigencomplex for $\FT_3$ with eigenmap $\b:\one[-2]\rightarrow \FT_3$.  Thus, $J_3 M\simeq M(2)[-2]$.  However, there is no forward or backward eigenmap $\one(2)[-2]\leftrightarrow J_3$ which induces this equivalence!  In fact, any chain map of this degree is null-homotopic.

Instead, one can use the eigenmaps for $\FT_3$ and $\FT_2^{-1}$ to create an \emph{eigenroof}: a pair of maps $\one(2)[-2] \longleftarrow \FT_2\inv[-2] \longrightarrow J_3$ where
each map becomes a homotopy equivalence once applied to $M$. Perhaps a theory of eigenroofs can be developed, similar to the theory of localizing classes in categories.

An alternative approach is the observation that $J_3$ does admit a genuine forward eigenmap when restricted to the eigencategory for $\b$, which is a monoidal category with monoidal
identity given by $\PB_\b$, the projection onto the $\b$-eigencategory. Similarly, it admits a backward eigenmap when restricted to the eigencategory for $\a$. One could attempt to study
relative eigenmaps in this fashion.

\appendix


\section{Commuting properties}
\label{sec:appendix}


\begin{notation}
We often abbreviate $F\otimes G$ simply by $FG$.
\end{notation}

Let $A$ be a ring and $f, g, g' \in A$, such that $f$ commutes with $g$ and with $g'$. Then $f$ commutes with $g - g'$, of course.

Let $\AC$ be a triangulated monoidal category and $F, G, G' \in \AC$ be objects, such that $FG \cong GF$ and $FG' \cong G' F$. If $\a \co G \to G'$, then it is
entirely possible that $F \ot \Cone(\a) \ncong \Cone(\a) \ot F$.

As just illustrated, many basic aspects of ring theory relating to commutativity will fail when rings are replaced by (triangulated) monoidal categories. This is a significant technical
issue in this paper. Diagonalization of a linear operator $f$ is essentially the in-depth study of the commutative ring generated by $f$. Unfortunately, the subcategory generated by an
object in a triangulated monoidal category need not be tensor commutative. Here is an example.

\begin{example} Let $\AC=\KC^b(\SBim_2)$, and let $F$ denote the Rouquier complex associated to a crossing. For any polynomial $g\in \Q[x_1,x_2]$, left or right multiplication by $g$
gives rise to bimodule endomorphisms $gB$ and $Bg$ for any bimodule $B$ (or complex of bimodules). For $F$ we have $x_1F\simeq Fx_2$ and $x_2F\simeq F x_1$. Then $\Cone(x_1F)$ does not
commute with $F$, since $F\Cone(x_1F)\simeq \Cone(x_2FF)$, while $\Cone(x_1F)F\cong \Cone(x_1FF)$. It is an exercise to show that these complexes are not isomorphic in $\AC$.
\end{example}

\begin{remark} This is one example which is not more easily done with $A$-modules, essentially because the monoidal identity of $\SBim_2$ has nontrivial endomorphisms $x_1$ and $x_2$, while the monoidal identity of $A$-modules has no nontrivial endomorphisms.
	
Let $F$ be the two-term complex $(A \to \Z)$. Then multiplication by $x$ does give an endomorphism of this complex (which acts by $1$ on the $\Z$ term), as does multiplication by any
scalar in $\Z$. One should think of left multiplication by $x_1$ (resp. $x_2$) in $\SBim_2$ as analogous to multiplication by $x$ (resp. $-x$) on $A$-modules, and of right
multiplication by $x_1$ (resp. $x_2$) as the action of $x$ (resp. $-x$) on the monoidal identity $\Z$, then tensored with $F$. One still does have $x\cdot F \simeq (-1)\cdot F$ and
$(-x)\cdot F \simeq F\cdot (+1)$. Unlike the example above, $\Cone(x\cdot F) \cong 0$ since multiplication by $x$ is an isomorphism, so we conclude that $\Cone(x \cdot F)$ does commute
with $F$. \end{remark}

All our constructions (eigencones, the idempotents $\idemp_i$, the complexes $\Cij{i}{j}$) are built as iterated cones of $F$ and scalar functors, so they are analogous to polynomials
in $F$, but it takes a lot of work to show (under certain hypotheses) that they commute with each other up to homotopy equivalence. Proving certain results amounts to chasing down
homotopies, homotopies of homotopies, etcetera.

In our main theorem (Theorem \ref{thm-firstEigenthm}), we made no assumptions about whether eigencones or the complexes $\Cij{i}{j}$ commute, although we did assume the
overly complex condition (PD3) in our definition of pre-diagonalizability (Definition \ref{def:prediag}). We have no examples where we can prove (PD3) without proving that the eigencones
commute up to isomorphism. Moreover, other compabilities in the theory, such as the ability to express the periodicity maps for projectors in terms of eigenmaps, relied on additional commutativity assumptions. This appendix carefully studies the homological obstructions relevant for these commutativity properties.



For the appendix, let $\BC$ be an additive monoidal category, and let $\AC\subset \KC(\BC)$ be a monoidal homotopy category.  Note that the tensor product and associator in $\AC$ are induced from those in $\BC$, but the monoidal identity and, hence unitor may be different.  Thus, let us denote the monoidal identity in $\AC$ by $\IB$. equal to those in $\KC(\BC)$   The categories $\BC$ and $\AC$ will be fixed throughout.  We will also fix a category of scalar objects $\KS \subset \ZC(\AC)$.  For each scalar object $\l$, let $\tau_\l$ denote the natural isomorphism $\tau_\l: \l\otimes (-) \rightarrow (-) \otimes \l$.  We will assume that:
\begin{itemize}
\item For every pair of scalar objects we have $\tau_{\l,\mu}\simeq \tau_{\mu,\l}\inv$.
\end{itemize}
We leave it as an exercise to show that this property are satisfied in our main examples, where scalar objects are just shifted copies of $\IB$.  The importance of this assumption will become clear momentarily. 

Recall the $\KS^\times$-graded hom space $\Homg_\AC(C,D)= \bigoplus_\l \Hom_\AC(\l C,D)$, where $\l$ runs over (a chosen set of representatives for) the isomorphism classes of invertible scalar objects $\l\in \KS^\times$.

A related but different object is the $\Z$-graded internal Hom complex $\Homb_\AC(C,D)$.  The differential on $\Homb_\AC(C,D)$ maps $f \mapsto [d,f]:=d_C\circ f - (-1)^{|f|} f\circ d_A$.  Here $|f|\in \Z$ denotes the homological degree of $f$.  The homology of $\Homb_\AC(A,B)$ is $\bigoplus_k\Hom(C[-k],D)$, which is a homogeneous subspace of $\Homg_\AC(C,D)$.

The following notation will help us to abbreviate many commutative diagrams in this section.
\begin{notation}
For any object $C$, we also denote by $C$ the operation of taking the tensor product (of a morphism) with $\Id_C$. Thus, $Cf=\Id_C\otimes f$ and $fC=f\otimes \Id_C$.
\end{notation}

\subsection{Obstructions to commutativity}
\label{subsec:commutativity}


\begin{lemma}\label{lemma:FcommutesWithCone}
Suppose we are given complexes $F,G_0,G_1\in \AC$, homotopy equivalences $\phi_i:G_i\otimes F\rightarrow F\otimes G_i$ ($i=0,1$), and a chain map $g:G_0\rightarrow G_1$.  If $\phi_1\circ (g\otimes \Id_F) - (\Id_F\otimes g)\circ \phi_0$ is null-homotopic, then the equivalences $\phi_i$ induce a homotopy equivalence $\Cone(g)\otimes F\rightarrow F\otimes \Cone(g)$.
\end{lemma}
\begin{proof}
Consider the following diagram:
\[
\begin{diagram}
G_0\otimes F & \rTo^{g\otimes \Id_F} & G_1\otimes F & \rTo & \Cone(g)\otimes F & \rTo & G_0\otimes F[1]\\
\dTo^{\phi_0} && \dTo^{\phi_1} && \dDashto^{\Phi} && \dTo^{\phi_0[1]}\\ 
F\otimes G_0 & \rTo^{\Id_F\otimes g} & F\otimes G_1 & \rTo & F\otimes \Cone(g) & \rTo & F\otimes G_0[1]
\end{diagram}
\]
in which the rows are distinguished triangles and the left square commutes up to homotopy by hypothesis.  There exists a map $\Phi$ making this diagram commute up to homotopy (this is one of the axioms of triangulated categories).  This $\Phi$ is a homotopy equivalence by the 5-lemma for triangulated categories.
\end{proof}

The conditions of Lemma \ref{lemma:FcommutesWithCone} occur often enough that they deserve a name.

\begin{definition}\label{def:gCommutesWithF}
Let $\phi_i:G_i\otimes F\rightarrow F\otimes G_i$ be given homotopy equivalences ($i=1,2$).  We say $g:G_1\to G_2$ \emph{commutes with $F$} if $\phi_1\circ (g\otimes \Id_F) - (\Id_F\otimes g)\circ \phi_0=[d,h_g]$ for some $h_g\in \Hom^{-1}(F\otimes G_1, G_2\otimes F)$.
\end{definition}

We now wish to discuss a homological obstruction which may prevent two cones $\Cone(F_0\buildrel f\over \rightarrow F_1)$ and $\Cone(G_0\buildrel g\over \rightarrow G_1)$ from commuting with one another, even when their terms commute ($G_i\otimes F_j\simeq F_j\otimes G_i$).

Let $F_i,G_i\in \AC$ be complexes ($i=0,1$), and let $\phi_{i,j}:G_i\otimes F_j\rightarrow F_j\otimes G_i$ be homotopy equivalences.   Suppose that the chain maps $\phi_{1,j}\circ(g\otimes \Id_{F_j})-(\Id_{F_j}\otimes g)\circ \phi_{0,j}$ and $\phi_{i,1}\circ (\Id_{G_i}\otimes f) - (f\otimes \Id_{G_i})\circ \phi_{i,0}$ are null-homotopic.  Choose homotopies $h_i\in \Hom^{-1}(G_i\otimes F_0,F_1\otimes G_i)$ and $k_i\in \Hom^{-1}(G_0\otimes F_i,F_i\otimes G_1)$ such that
\begin{eqnarray}
\label{eq:h} d\circ h_i + h_i\circ d &=& \phi_{i,1}\circ (\Id_{G_i} \otimes f) - (f\otimes \Id_{G_i})\circ \phi_{i,0}\\
\label{eq:k} d\circ k_i + k_i\circ d &=& \phi_{1,i}\circ (g \otimes \Id_{F_i}) - (\Id_{F_i}\otimes g)\circ \phi_{0,i}
\end{eqnarray}
\begin{proposition}\label{prop:zCycle}
The element $z_{f,g}\in \Hom^{-1}(G_0F_0,F_1G_1)$ defined by
\begin{equation}\label{eq:zCycle}
z_{f,g} \ \ := \ \ h_1\circ (gF_1) - (F_0g)\circ h_0-k_1\circ(G_0f)+(fG_0)\circ k_0
\end{equation}
is a cycle.\qed
\end{proposition}

\begin{definition}\label{def:secondaryObstruction}
Let $F_i,G_i\in \AC$ be complexes ($i=0,1$), and let $f:F_0\rightarrow F_1$, $g:G_0\rightarrow G_1$ be chain maps.  We refer to the cycle $z_{f,g}\in \Hom^{-1}(G_0F_0,F_1G_1)$ as the \emph{obstruction to commuting $\Cone(f)$ past $\Cone(g)$}. We also refer to $z_{f,g}$ as a \emph{secondary obstruction}, the primary obstructions being the chain maps in the statement of Lemma \ref{lemma:FcommutesWithCone}.
\end{definition}
\begin{remark}The definition of $z_{f,g}$ requires several choices: we must choose homotopy equivalences $\phi_{i,j}:G_i\otimes F_j\rightarrow F_j\otimes G_i$ ($i,j\in\{0,1\}$), as well homotopies $h_i,\in \Hom^{-1}(G_iF_0,F_1G_i)$, $k_i\in \Hom^{-1}(G_0F_i,F_iG_1)$ ($i=1,2$) satisfying equations (\ref{eq:h}) and (\ref{eq:k}).  We usually regard the $\phi_{i,j}$ as being fixed, while the homotopies $h_i,k_i$ are chosen arbitrarily.  Whenever we say that $z_{f,g}$ is null-homotopic, we mean that $z_{f,g}$ is null-homotopic for \emph{some} choice of homotopies $h_i,k_i$.
\end{remark}

\begin{lemma}\label{lemma:conesCommute}
Let $F_i,G_i\in \AC$ be chain complexes ($i=0,1$), and let $\phi_{i,j}:G_i F_j\rightarrow F_jG_i$ be given homotopy equivalences, and let $f:F_0\rightarrow F_1$, $g:G_0\rightarrow G_1$ be chain maps such that the secondary obstruction $z_{f,g}$ is defined.  If $z_{f,g}$ is null-homotopic, then $\Cone(f)\otimes \Cone(g)\simeq \Cone(g)\otimes \Cone(f)$.
\end{lemma}
\begin{proof}
From the hypothesis that $z_{f,g}$ is defined, we know that each of the following squares commutes up to homotopy:
\[
\begin{diagram}
G_0F_j & \rTo^{g\Id} & G_1 F_j \\
\dTo^{\phi_{0,j}} && \dTo^{\phi_{1,j}}\\
F_j G_0 &\rTo^{\Id g} & F_j G_1
\end{diagram}
\ \ \ \ \ \ \ \ \ \ \ \ \ \ \ \ 
\begin{diagram}
G_iF_0 & \rTo^{\Id f} & G_i F_1 \\
\dTo^{\phi_{1,j}} && \dTo^{\phi_{i,0}}\\
F_0 G_i &\rTo^{f \Id} & F_1 G_i
\end{diagram}
\]
Let us choose homotopies $h_i$ and $k_i$ which realize the comutativity of these squares up to homotopy; these elements satisfy (\ref{eq:h}) and (\ref{eq:k}).  By hypothesis, we may choose $h_i$ and $k_i$ such that the obstruction $z_{f,g}\simeq 0$.  Thus, we also have a homotopy $\ell\in\Hom^{-2}(G_0F_0,F_1G_1)$ such that $d\circ \ell- \ell\circ d = z_{f,g}$.  Below, we will often omit the subscripts from $h_i$ and $k_i$, for readability.

 Commutativity of the left squares implies that we have equivalences $\Cone(g)\otimes F_j\rightarrow F_j\otimes \Cone(g)$.  We intend to show that these equivalences fit into a homotopy commutative square
\begin{equation}\label{eq:coneCommuteSquare}
\begin{diagram}
\Cone(g)\otimes F_0 & \rTo^{\Id f} & \Cone(g)\otimes F_1\\
\dTo^{\simeq}&& \dTo^{\simeq}\\
F_0\otimes \Cone(g) & \rTo^{f\Id} & F_1\otimes \Cone(g)
\end{diagram}
\end{equation}
This will imply that $\Cone(g)$ commutes with $\Cone(f)$ up to homotopy by Lemma \ref{lemma:FcommutesWithCone}.  Now, the top and right arrows of diagram (\ref{eq:coneCommuteSquare}) can be expressed by the diagrams:
\[
\begin{diagram}
G_0F_0 & \rTo^{g\Id} & G_1 F_0\\
\dTo^{\Id f} && \dTo^{\Id f}\\
G_0F_1 & \rTo^{g\Id}& G_1F_1
\end{diagram}
\ \ \ \ \ \ \ \ \ \ \ \ \ \ \ \ 
\begin{diagram}
G_0F_1 & \rTo^{g\Id}& G_1F_1\\
\dTo^{\phi_{1,1}} &\rdTo^{k}& \dTo^{\phi_{0,1}}\\
F_1 G_0& \rTo^{\Id g}& F_1G_1
\end{diagram}
\]
Let $\Psi_1:\Cone(g)\otimes F_0\rightarrow F_1\Cone(g)$ denote the composition of these maps.  The bottom and left arrows of diagram (\ref{eq:coneCommuteSquare}) can be expressed by the diagrams
\[
\begin{diagram}
G_0F_0 & \rTo^{g\Id} & G_1 F_0\\
\dTo^{\phi_{1,0}} &\rdTo^{k}& \dTo^{\phi_{0,0}}\\
F_0G_0 & \rTo^{\Id g}& F_0 G_1
\end{diagram}
\ \ \ \ \ \ \ \ \ \ \ \ \ \ \ \ 
\begin{diagram}
F_0G_0 & \rTo^{\Id g}& F_0 G_1\\
\dTo^{f \Id} && \dTo^{f\Id}\\
F_1 G_0& \rTo^{\Id g}& F_1G_1
\end{diagram}
\]
Let $\Psi_0:\Cone(g)\otimes F_0\rightarrow F_1\Cone(g)$ denote the composition of these maps.  Then $\Psi_1-\Psi_0$ can be expressed diagrammatically by:
\[
\begin{diagram}
G_0F_0&\rTo^{g\Id}& G_1F_0\\
\dTo^{\phi\circ (\Id f) - (f\Id)\circ \phi} & \rdTo^{k\circ (\Id f) - (f\Id)\circ k} & \dTo_{\phi\circ (\Id f)- (f\Id)\circ\phi }\\
F_1 G_0 &\rTo^{\Id g}& F_1G_1
\end{diagram}
\]
This map is null-homotopic, with homotopy given by the diagram
\begin{equation}\label{eq:coneCommuteHomotopy}
\begin{diagram}
G_0F_0&\rTo^{g\Id}& G_1F_0\\
\dTo^{-h} & \rdTo^{\ell} & \dTo_{h}\\
F_1 G_0 &\rTo^{\Id g}& F_1G_1
\end{diagram}.
\end{equation}
Here, the homotopies $h$ satisfy
\[
d\circ h +h\circ d = \phi\circ (\Id f)- (f\Id)\circ\phi,
\]
and the homotopy $\ell$ satisfies
\[
d\circ \ell - \ell\circ d = h\circ (g \Id) - (\Id g)\circ h - k\circ(\Id f)+(f\Id)\circ k.
\]
The sign on the left arrow in the homotopy (\ref{eq:coneCommuteHomotopy}) is there because in forming the mapping cone on the top (or bottom) arrows, the left most terms aquire a shift of $[1]$, which negates the differential.  The sign on $-h$ counteracts this.

This proves that $\Psi_1$ is homotopic to $\Psi_0$, hence the square \eqref{eq:coneCommuteSquare} commutes up to homotopy.  This completes the proof.
\end{proof}

\subsection{Zigzag complexes}
\label{subsec:zigzags}

In this section we introduce a construction which generalizes that of mapping cone, and we show that the commutativity results of the previous section extend to this more general setting.

\begin{definition}\label{def:zigzagCx}
Let $\CC$ be an additive category.  Suppose we are given a collection of chain complexes $X_0,X_1,\ldots,\in \KC(\CC)$ and chain maps as in the following diagram:
\[
\begin{tikzpicture}
\node(aa) at (0,0) {$X_0$};
\node(ba) at (2.5,1) {$X_1$};
\node(ca) at (5,0) {$X_2$};
\node(da) at (7.5,1) {$X_3$};
\node(ea) at (10,0) {$\cdots$};
\tikzstyle{every node}=[font=\small]
\path[->,>=stealth',shorten >=1pt,auto,node distance=1.8cm]
(ba) edge node[above] {$ g_0$} (aa)
(ba) edge node[above] {$ g_1 $} (ca)
(da) edge node[above] {$ g_2$} (ca)
(da) edge node[above] {$ g_3$} (ea);
\end{tikzpicture}.
\]
Note the direction of these arrows; the even numbered maps point to the left, while the odd numbered maps point to the right.  Let us abbreviate this diagram by writing $\{X_k,g_k\}$, or simply $X_\bullet$.  The \emph{total complex} of $X_\bullet$ is the complex obtained by first applying the shift $[1]$ to each odd-numbered term $X_{2i+1}$, then taking total complex in the usual manner for twisted complexes (see \S \ref{subsec-convolutions}).  In other words:
\[
\Tot(X_\bullet)  \ \ = \ \ \Tot\Big(\begin{diagram} X_0 & \lTo & X_1[1] & \rTo & X_2 & \lTo & X_3[1] & \rTo & \cdots \end{diagram} \Big).
\]
We do not necessarily assume any local finiteness properties, so there are two potentially distinct total complexes, $\Tot^\oplus(X_\bullet)$ and $\Tot^\Pi(X_\bullet)$.  These complexes live in $\KC(\BC^\oplus)$ and $\KC(\BC^\Pi)$, where $\BC^\oplus$ and $\BC^\Pi$ are the categories obtained from $\BC$ by adjoining countable direct sums and products, respectively.  We refer to any such complex obtained in this way as a \emph{zigzag complex}.  We sometimes let $e\in \{\oplus,\Pi\}$, and we write $\Tot^e$ to mean either $\Tot^\oplus$ or $\Tot^\Pi$.
\end{definition}
\begin{remark}
If $g_0:X_1\rightarrow X_0$ is a chain map, then $\Cone(g_0)$ is a zigzag complex with $0=X_2=X_3=\cdots $.  In general, zigzag complexes fit into a distinguished triangle of the form:
\[
\bigoplus_{i\geq 0} X_{2i+1}\ \  \longrightarrow \ \  \bigoplus_{i\geq 0} X_{2i} \ \  \longrightarrow \ \ \Tot^\oplus(X_\bullet) \ \  \longrightarrow \ \ \bigoplus_{i\geq 0} X_{2i+1}[1],
\]
and similarly for $\Tot^\Pi(X_\bullet)$.  These distinguished triangles characterize the zigzag complexes up to isomorphism, and can be used to define the zigzag construction any triangulated category which contains countable sums, respectively products.  We will not need this fact.
\end{remark}
\begin{remark}\label{rmk:indexingSet}
The indexing set for the convolution $\Tot^e(X_\bullet)$ is $\Z_{\geq 0}$ with a non-standard partial order: $i<j$ iff $i$ is odd and $j=i\pm 1$.  This poset is both upper finite and lower finite (it satisfies the ascending and descending chain conditions). 
\end{remark}

\begin{remark}\label{rmk:Cab}
Let $F\in \KC(\BC)$ be given, let $\a,\b\in \Homg(\one,F)$ be maps of degrees $\l,\mu$, respectively, and assume that $\l\mu\inv$ is small.  The complexes $\Cij{\a}{\b}(F)$ from Definition \ref{def:interpolationCx} are special cases (in fact for us they are the main examples) of zigzag complexes.   Indeed, $\Cij{\a}{\b}(F)$ is the total complex of the following zigzag diagram:
\begin{equation}\label{eq:appCab}
\begin{tikzpicture}[baseline=1em]
\node(aa) at (0,0) {$\displaystyle \mu\inv F$};
\node(ba) at (3,1) {$\displaystyle \l\mu\inv $};
\node(ca) at (6,0) {$\displaystyle \l\mu^{-2} F$};
\node(da) at (9,1) {$\displaystyle\l^2\mu^{-2} $};
\node(ea) at (12,0) {$\cdots$};
\tikzstyle{every node}=[font=\small]
\path[->,>=stealth',shorten >=1pt,auto,node distance=1.8cm]
(ba) edge node[above] {$\mu\inv\a$} (aa)
(ba) edge node[above] {$ -\l\mu\inv \b $} (ca)
(da) edge node[above] {$ \l\mu^{-2}\a$} (ca)
(da) edge node[above] {$- \l^2\mu^{-2}\b$} (ea);
\end{tikzpicture}.
\end{equation}
Similarly, $ \Cij{\b}{\a}(F)[1]$ (note the homological shift) is the total complex of the following zigzag diagram:
\begin{equation}\label{eq:Cba}
\begin{tikzpicture}[baseline=1em]
\node(za) at (-3,1) {$\one $};
\node(aa) at (0,0) {$\mu\inv F$};
\node(ba) at (3,1) {$\l\mu\inv $};
\node(ca) at (6,0) {$\l\mu^{-2} F$};
\node(da) at (9,1) {$\cdots$};
\tikzstyle{every node}=[font=\small]
\path[->,>=stealth',shorten >=1pt,auto,node distance=1.8cm]
(za) edge node[above] {$-\mu\inv\b$} (aa)
(ba) edge node[above] {$ \mu\inv \a$} (aa)
(ba) edge node[above] {$ -\l\mu^{-2}\b $} (ca)
(da) edge node[above] {$ \l\mu^{-2}\a$} (ca);
\end{tikzpicture}.
\end{equation}
\end{remark}

The following is a sort of homotopy invariance for zigzag complexes.
\begin{lemma}\label{lemma:zigzagSimp}
Consider a diagram of the form:
\[
\begin{tikzpicture}[baseline=0em]
\node(aa) at (0,2) {$X_0$};
\node(ba) at (2.5,2) {$X_1$};
\node(ca) at (5,2) {$X_2$};
\node(da) at (7.5,2) {$X_3$};
\node(ea) at (10,2) {$\cdots$};
\node(ab) at (0,0) {$Y_0$};
\node(bb) at (2.5,0) {$Y_1$};
\node(cb) at (5,0) {$Y_2$};
\node(db) at (7.5,0) {$Y_3$};
\node(eb) at (10,0) {$\cdots$};
\tikzstyle{every node}=[font=\small]
\path[->,>=stealth',shorten >=1pt,auto,node distance=1.8cm]
(aa) edge node[above] {$f_0$} (ba)
(ca) edge node[above] {$f_1$} (ba)
(ca) edge node[above] {$f_2$} (da)
(ea) edge node[above] {$f_3$} (da)
(ab) edge node[above] {$g_0$} (bb)
(cb) edge node[above] {$g_1$} (bb)
(cb) edge node[above] {$g_2$} (db)
(eb) edge node[above] {$g_3$} (db)
(aa) edge node {$\phi_0$} (ab)
(ba) edge node {$\phi_1$} (bb)
(ca) edge node {$\phi_2$} (cb)
(da) edge node {$\phi_3$} (db);
\end{tikzpicture}
\]
where the $X_i$ and $Y_i$ are arbitrary chain complexes and the arrows are chain maps.  Let $e\in \{\oplus,\Pi\}$.  If the squares commute up to homotopy, then the above diagram induces a chain map of zigzag complexes $\Phi:\Tot^e(X_\bullet)\rightarrow \Tot^e(Y_\bullet)$.  If the vertical arrows are homotopy equivalences, then $\Phi$ is a homotopy equivalence.
\end{lemma}
\begin{proof}
Note that $f_i$ maps $X_{i+1}\rightarrow X_i$ if $i$ is even, and maps $X_i\rightarrow X_{i+1}$ if $i$ is odd, and similarly for $g_i$.  For each $k\in \Z_{\geq0}$ let $h_{2k}\in \Hom^{-1}(X_{2k+1},X_{2k})$ and $h_{2k+1}\in \Hom^{-1}(X_{2k+1},X_{2k+2})$ denote the homotopies such that
\[
d\circ h + h \circ d = \phi\circ f - g\circ \phi \ \ \ \ \ \ (\text{indices omitted}).
\]
It is straightforward to construct the chain map $\Phi:\Tot^e(X_\bullet)\rightarrow \Tot^e(Y_\bullet)$ using the $\phi_i$ and the homotopies $h_i$.  For the reader's convenience, we include a diagram which represents the mapping cone of $\Phi$:
\[
\Cone(\Phi) \ \cong \ \Tot^e\left(
\begin{tikzpicture}[baseline=2.3em]
\node(aa) at (0,2) {$X_0[1]$};
\node(ba) at (2.5,2) {$X_1[2]$};
\node(ca) at (5,2) {$X_2[1]$};
\node(da) at (7.5,2) {$X_3[2]$};
\node(ea) at (10,2) {$\cdots$};
\node(ab) at (0,0) {$Y_0$};
\node(bb) at (2.5,0) {$Y_1[1]$};
\node(cb) at (5,0) {$Y_2$};
\node(db) at (7.5,0) {$Y_3[1]$};
\node(eb) at (10,0) {$\cdots$};
\tikzstyle{every node}=[font=\small]
\path[->,>=stealth',shorten >=1pt,auto,node distance=1.8cm]
(aa) edge node[above] {$-f_0$} (ba)
(ca) edge node[above] {$-f_1$} (ba)
(ca) edge node[above] {$-f_2$} (da)
(ea) edge node[above] {$-f_3$} (da)
(ab) edge node[above] {$g_0$} (bb)
(cb) edge node[above] {$g_1$} (bb)
(cb) edge node[above] {$g_2$} (db)
(eb) edge node[above] {$g_3$} (db)
(aa) edge node[left] {$\phi_0$} (ab)
(ba) edge node {$\phi_1$} (bb)
(ca) edge node[left] {$\phi_2$} (cb)
(da) edge node {$\phi_3$} (db)
(ba) edge node[xshift=-.6cm, yshift=.5cm]   {$h_0$} (ab)
(da) edge node[xshift=-.6cm, yshift=.5cm]  {$h_2$} (cb)
(ba) edge node[xshift=-.2cm, yshift=.1cm]  {$h_1$} (cb)
(da) edge node[xshift=-.2cm, yshift=.1cm]  {$h_3$} (eb);
\end{tikzpicture}
\right)
\]
Checking that $\Phi$ is a chain map amounts to checking that $d^2=0$ for the above chain complex.  This is straightforward, and left to the reader.  Now, suppose the $\phi_i$ are homotopy equivalences, and set $Z_i:=\Cone(\phi_i)$.  The mapping cone on $\Phi$ can be reassociated  as follows:
\[
\Cone(\Phi)\ \ = \ \ \Tot^e\Big(\begin{diagram} Z_0 & \lTo & Z_1[1] & \rTo & Z_2 & \lTo & Z_3[1] & \rTo & \cdots \end{diagram} \Big).
\]
The indexing set for this convolution is both upper finite and lower finite (Remark \ref{rmk:indexingSet}), hence Proposition \ref{prop:simplifications} applies: each $Z_i$ is contractible, so the convolution $\Cone(\Phi)\cong \Tot^e(Z_\bullet)$ is contractible as well.  This shows that $\Phi$ is a homotopy equivalence, and completes the proof.
\end{proof}

A zigzag complex $\Tot(\{G_i,g_i\})$ is ``built" from the mapping cones $\Cone(g_i)$.  Roughly speaking, commuting a complex $F$ past $\Tot(\{G_i,g_i\})$ essentially ammounts to commuting $F$ past each $\Cone(g_i)$.  The following makes this precise.

\begin{lemma}\label{lemma:FcommutesWithZigzag}
Let $F\in \AC$ be given, and suppose we have a zigzag diagram $\{G_k,g_k\}_{k\in \Z_{\geq 0}}$.  Suppose we are given homotopy equivalences $\phi_i:G_i\otimes F \rightarrow F\otimes G_i$ such that, for each odd index $i=2k+1$ the chain maps $\phi_{i-1}\circ (g_{i-1}F)-(Fg_{i-1})\circ \phi_i$ and $\phi_{i+1}\circ (g_i F) - (F g_i) \circ \phi_i$ are null-homotopic.
%
%
Then the $\phi_i$ induce a homotopy equivalence $\Phi:\Tot^e(G_\bullet)\otimes F\rightarrow F\otimes \Tot^e(G_\bullet)$.
\end{lemma}
\begin{proof}
Follows immediately from Lemma \ref{lemma:zigzagSimp}.
\end{proof}

Now we have our main commutativity result.

\begin{proposition}\label{prop:zigzagsCommute}
Suppose we have zigzag diagrams
\[
\begin{tikzpicture}
\node(aa) at (0,0) {$F_0$};
\node(ba) at (2.5,0) {$F_1$};
\node(ca) at (5,0) {$F_2$};
\node(da) at (7.5,0) {$F_3$};
\node(ea) at (10,0) {$\cdots$};
\tikzstyle{every node}=[font=\small]
\path[->,>=stealth',shorten >=1pt,auto,node distance=1.8cm]
(ba) edge node[above] {$ f_0$} (aa)
(ba) edge node[above] {$ f_1 $} (ca)
(da) edge node[above] {$ f_2$} (ca)
(da) edge node[above] {$ f_3$} (ea);
\end{tikzpicture}
\]
and
\[
\begin{tikzpicture}
\node(aa) at (0,0) {$G_0$};
\node(ba) at (2.5,0) {$G_1$};
\node(ca) at (5,0) {$G_2$};
\node(da) at (7.5,0) {$G_3$};
\node(ea) at (10,0) {$\cdots$};
\tikzstyle{every node}=[font=\small]
\path[->,>=stealth',shorten >=1pt,auto,node distance=1.8cm]
(ba) edge node[above] {$ g_0$} (aa)
(ba) edge node[above] {$ g_1 $} (ca)
(da) edge node[above] {$ g_2$} (ca)
(da) edge node[above] {$ g_3$} (ea);
\end{tikzpicture}.
\]
Assume we are given homotopy equivalences $\phi_{i,j}:G_iF_j\rightarrow F_jG_i$ for $i,j\in \Z_{\geq 0}$.  If the obstructions $z_{f_i,g_j}$ are null-homotopic for all $i,j\geq 0$, then $\Tot(F_\bullet)$ commutes with $\Tot(G_\bullet)$ up to homotopy.
\end{proposition}
\begin{proof}
The proof is very similar to the proof of Lemma \ref{lemma:FcommutesWithZigzag}.  The existence of $z_{f_i,g_j}$ for all $i,j$ implies that every possible square involving the $\phi_{i,j}$, $g_i$, and $f_j$, commutes up to homotopy.  By Lemma \ref{lemma:FcommutesWithZigzag}, this in particular gives us an equivalence $\Tot(G_\bullet)F_j\rightarrow F_j\Tot(G_\bullet)$ for all $j\in \Z_{\geq 0}$.  We intend to show that these equivalences fit into homotopy commutative squares:
\begin{equation}\label{eq:FzigzagSquare}
\begin{diagram}
\Tot(G_\bullet)\otimes F_j & \rTo^{\Id \otimes f} & \Tot(G_\bullet) \otimes F_{j\pm 1}\\
\dTo^{\simeq} && \dTo^{\simeq }\\
F_j\otimes \Tot(G_\bullet) & \rTo^{f\otimes \Id} & F_{j\pm 1}\otimes \Tot(G_\bullet)
\end{diagram}
\end{equation}
for all $j\in \Z_{\geq 0}$ odd.  Given this, Lemma \ref{lemma:FcommutesWithZigzag} implies that $\Tot(F_\bullet)$ commutes with $\Tot(G_\bullet)$, as claimed.  The proof that the square (\ref{eq:FzigzagSquare}) commutes up to homotopy is very similar to the proof that the square (\ref{eq:coneCommuteSquare}) commutes up to homotopy, and we will leave the details to the reader.
\end{proof}


The results in this section give sufficient conditions for interpolation complexes to commute with one another.  In the nicest possible situations, the obstructions vanish for degree reasons, and the following gives a nice criterion for commutativity.
\begin{corollary}
Let $\l_i$ and $\mu_i$ be invertible scalar objects with $\l_2\l_1\inv$ and $\mu_2\mu_1\inv$ small ($i=1,2$).  Let $F,G\in \AC$ be objects and $\a_i:\l_i\rightarrow F$, $\b_i:\mu_i\rightarrow G$ be scalars.  Suppose that $\a_i$ commute with $G$ and $\b_i$ commute with $F$ for $i=1,2$.  If the degree $-1$ homology groups of $\Hom^\bullet(\l_i\otimes \mu_j, F\otimes G)$ are zero for $i,j\in \{1,2\}$, then the then $\Cij{\a_2}{\a_{1}}(F)$ and $\Cij{\a_1}{\a_2}(F)$ commute with $\Cij{\b_2}{\b_{1}}(G)$ and $\Cij{\b_1}{\b_2}(G)$.\qed
\end{corollary}

\subsection{The self obstruction}
\label{subsec:selfObstruction}

There is another cycle that could potentially obstruct some of our constructions.  
\begin{proposition}\label{prop:wCycle}
If $\a\in\Homg(\one,F)$ commutes with $F$, then $w_\a = h_{\a}\circ (\a \l)$ is a cycle in $\Hom^{-1}(\l \l, FF)$.  Here, $\l$ is the degree of $\a$ and $h_\a$ is a chosen homotopy, as in Definition \ref{def:gCommutesWithF}. \qed
\end{proposition}
\begin{proposition}\label{prop:selfObstruction}
If $w_\a\simeq 0$, then $\Cone(\a)^{\otimes 2} \simeq (F\oplus \l[1])\otimes \Cone(\a)$.
\end{proposition}
\begin{proof}
If $w_\a=h_{\a}\circ (\a \l)$ is a boundary, then there exists $k_\a\in \Hom^{-2}(\l\l,FF)$ such that $d\circ \ell - \ell\circ d = h_{\a}\circ (\a \l)$.  We intend to show that $\a:\l\rightarrow F$ becomes null-homotopic upon tensoring on the left or right with $\Cone(\a)$.  Indeed, consider the element $m\in \Hom^{-1}(\l\Cone(\a), F\Cone(\a))$ represented by the non-horizontal arrows of the following diagram:
\[
\begin{tikzpicture}[baseline=-.12cm]
\node (aa) at (0,0){$\l\l$};
\node (ba) at (2.5,0){ $\l F$};
\node (ab) at (0,-2.5){$F\l $};
\node (bb) at (2.5,-2.5){$F F$ };
\path[->,>=stealth',shorten >=1pt,auto,node distance=1.8cm,font=\small]
(aa) edge node[left] {0} (ab)
(aa) edge node {$\l {\a}$} (ba)
(ab) edge node {$h_{\a}$} (bb)
(ba) edge node {$F{\a}$} (bb)
(aa) edge node[xshift=-.6cm,yshift=.5cm] {$-k_\a$} (bb);
\draw[frontline,->,>=stealth',shorten >=1pt,auto,node distance=1.8cm]
(ba) to node[xshift=-1.1cm,yshift=0cm] {$\tau_{\l, F}$} (ab);
\end{tikzpicture}
\]
This element satisfies $d\circ m+m\circ d = \a\otimes \Id_{\Cone(\a)}$.  This shows that $\a\otimes \Id_{\Cone(\a)}\simeq 0$.  On one hand, the cone on this map is  $\Cone(\a)^{\otimes 2}$.  On the other hand, the cone on this map is isomorphic to the cone on the zero map, which is $(F\oplus \l[1])\otimes \Cone(\a)$.  This completes the proof.  
\end{proof}

\printbibliography

\end{document}